\documentclass[12pt]{article}
\usepackage{graphics,graphicx}
\usepackage{amsmath}
\usepackage{amssymb,amsfonts}

\newcommand{\hh}{H_{PC}}

\newcommand{\thk}{\th^{k-1}}

\newcommand{\co}{c_{_{\! 1}}}
\newcommand{\ct}{c_{_{\! 2}}}
\newcommand{\glc}{\gl_{\crt}}
\newcommand{\thc}{\th_{\crt}}
\newcommand{\km}{\frac{1}{k-1}}

\newcommand{\thz}{\th_{_0}}
\newcommand{\sfrac}[2]{\mbox{$\frac{#1}{#2}$}}

\renewcommand{\and}{~~{\rm  and}~~}
\newcommand{\vb}{V\sm \VV}
\newcommand{\thto}{\th_{_{\!2}}}

\newcommand{\bd}{\bar{d}}
\newcommand{\bbr}{\bb[\, \,}
\newcommand{\thth}{\th_{_3}}
\newcommand{\gln}{\gl n}

\newcommand{\crt}{{\rm crt}}

\newcommand{\al}{\aleph}
\newcommand{\0}{\emptyset}

\newcommand{\gll}{\Lambda}
\newcommand{\gL}{\Lambda}

\newcommand{\thl}{\th_{_{\! \!\gl}}}

\newcommand{\dd}{\Delta}

\topmargin by -\headsep \textheight 8.9in \oddsidemargin 0pt
\evensidemargin \oddsidemargin \marginparwidth 0.5in \textwidth
6.5in
\newcommand{\bb}{\Big}

\newcommand{\beq}[1]{\begin{equation}\label{#1}}
\newcommand{\mn}[0]{\medskip\noindent}
\newtheorem{thm}{Theorem}[section]
\newtheorem{cor}[thm]{Corollary}
\newtheorem{lem}[thm]{Lemma}
\newtheorem{ex}[thm]{Example}
\numberwithin{equation}{section}

\newcommand{\bq}[2]{\begin{equation}\label{#1} #2 \end{equation}}

\newcommand{\C}{\mathcal{C}}
\renewcommand{\c}{\mathcal{C}}

\newcommand{\kk}{\frac{k}{k-1}}
\newcommand{\gd}{\delta}
\renewcommand{\th}{\theta}
\newcommand{\case}[4]{
\left\{ \begin{array}{ll} {#1} & \mbox{#2} \\ {#3}  & \mbox{#4}
\end{array} \right.}

\newcommand{\thi}{\th_{_i}}
\newcommand{\tho}{\th_{_{\! 1}}\hspace{-0.4mm}}

\newcommand{\tht}{\tilde{\th}}
\newcommand{\bn}{\bigskip\noindent}
\newcommand{\sm}{\setminus}
\newcommand{\ra}{\rightarrow}
\newcommand{\pr}{\Pr}
\newcommand{\sub}{\subseteq}
\newcommand{\qed}{\hfill$\square$}
\newcommand{\gs}{\sigma}
\newcommand{\raf}[1]{(\ref{#1})}
\newcommand{\pf}{\noindent {\bf Proof. }}
\newcommand{\eps}{\varepsilon}
\newcommand{\ga}{\alpha}
\newcommand{\gb}{\beta}
\newcommand{\gl}{\lambda}

\newcommand{\bean}{\begin{eqnarray*}}
\newcommand{\eean}{\end{eqnarray*}}

\newcommand{\gt}{\theta}

\newcommand{\poi}{{\rm Poi}}
\newcommand{\bin}{{\rm Bin}}

\newcommand{\eeq}{\end{equation}}
\newcommand{\old}[1]{}

\newcommand{\h}{{\cal H}}

\newcommand{\f}{{\cal F}}

\newcommand{\gls}{\Lambda}

\newcommand{\fnpk}{F_{_{\!{PC}}}(n,p\,;k)}

\newcommand{\reg}{\topmargin -3pt \advance \topmargin by
-\headheight \advance \topmargin by -\headsep \textheight 9.1in
\oddsidemargin -.1in \evensidemargin \oddsidemargin
\marginparwidth 0.5in \textwidth 6.7in}

\begin{document}

\old{ \vspace{6mm}
\begin{center}
{\Large \bf Poisson Cloning Model for Random Hypergraphs } \\[10mm]
{\it Jeong Han Kim \\ Microsoft Research \\ One Microsoft Way \\
  Redmond, WA98052  \\ jehkim@microsoft.com }
\end{center}

}

\renewcommand{\include}{\old}
\newcommand{\edc}{\end{document}}

\newcommand{\BB}{\tilde{B}}
\newcommand{\FF}{\tilde{F}}
\newcommand{\VV}{\tilde{V}}
\newcommand{\bt}{\tilde{b}}
\newcommand{\nt}{\tilde{N}}


\vspace{6mm}
\begin{center}
{\Large \bf Poisson Cloning Model for Random Graphs  } \\[10mm]
{\it Jeong Han Kim\footnote{This work was supported by the Korea Science and Engineering Foundation
 (KOSEF) grant funded by the Korea government(MOST) (No. R16-2007-075-01001-0).} \\ Yonsei University \\ Department of Mathematics \\ Seoul  120-749,
Korea \\ jehkim@yonsei.ac.kr }
\end{center}

\vspace{5mm} \mn {\bf Abstract.} In the random graph $G(n,p)$ with
$pn$ bounded, the degrees of the vertices are almost i.i.d Poisson
random variables with mean $\gl:= p(n-1)$. Motivated by this fact,
we introduce the Poisson cloning model $G_{PC} (n,p)$ for random
graphs in which the degrees are i.i.d Poisson random variables
with mean $\gl$. Then, we first establish a theorem that shows the
new model is equivalent to the classical model $G(n,p)$ in an
asymptotic sense. Next, we introduce a useful algorithm, called
the cut-off line algorithm,  to generate the random graph $G_{PC}
(n,p)$. The Poisson cloning model  $G_{PC}(n,p)$  equipped with
the cut-off line algorithm enables us to very precisely analyze
the sizes of the largest component and the $t$-core of $G(n,p)$.
This new approach to the problems yields not only elegant proofs
but also improved bounds that are essentially best possible.

\mn We also consider the Poisson cloning models for random
hypergraphs and random $k$-SAT problems. Then, the $t$-core
problem for random hypergraphs and the pure literal algorithm for
random $k$-SAT problems are analyzed.

\section{Introduction}
The notion of a random graph was first introduced  in 1947 by
Erd\H{o}s \cite{Erd47} to show the existence of a graph with a
certain Ramsey property. A decade later, the theory of the random
graph began with the paper entitled {\em On Random Graphs I }~ by
Erd\H{o}s and R\'enyi \cite{ER59}, and the theory had been
developed  by
 a series \cite{ER60,ER61, ER64, ER66, ER68} of papers of
them. Since then, the subject has become one of the most  active
research areas. Many
 researchers have
devoted themselves to studying various properties of random
graphs, such as the emergence of the giant component
\cite{ER60,BB84a,TL90c}, the connectivity \cite{ER59,ER61, BT85},
the existence of perfect matching \cite{ER64,ER66,ER68,BT85}, the
existence of Hamiltonian cycle(s) \cite{KS83,BB84b,BF85}, the
$k$-core problem \cite{BB84b, TL91d, PSW}, and the graph
invariants like the independence number \cite{BE76,Ma76} and the
chromatic number \cite{SS87,BB88,TL91b}. (The list of references
here is far from being exhaustive.)

There are two canonical models for random graphs, both of which
were originated in the simple model introduced  in  \cite{Erd47}.
In the binomial model $G(n,p)$ on a set $V$ of $n$ vertices, each
of ${n \choose 2}$ possible edges is in the graph with probability
$p$, independently of other edges. Thus, the probability of
$G(n,p)$ being a fixed graph $G$ with $m$ edges is $p^m
(1-p)^{{n\choose 2} -m}$. The uniform model $G(n,m)$ on $V$ is a
graph chosen uniformly at random from  the set of all graphs on
$V$ with $m$ edges. Hence, $G(n,m)$ becomes a fixed graph $G$ with
probability ${{n \choose 2} \choose m}^{-1}$, provided $G$ has $m$
edges. Most of asymptotic behaviors of the two models are almost
identical if their expected numbers of edges are the same. (See
Proposition 1.13 in \cite{JLR}.) The random graph process, in
which random edges are added one by one, is also extensively
studied.
For more about models and/or basics of random graphs, we recommend
two books with the identical title {\em Random Graphs} by Bollob\'as
\cite{RGB}, and by Janson, {\L}uczak and Ruci\'nski \cite{JLR}.

The phase transition phenomenon is among  most interesting topics of
random graphs. Specifically, the phase transition phenomena
regarding the emergences of the giant (connected) component and
the $t$-core problem have attracted much attention.
  In  their  monumental paper entitled {\em On the Evolution of Random
  Graphs} \cite{ER60},   Erd\H{o}s and
 R\'enyi proved that, for the size  $\ell_1 (n, p) $ of
 the largest component  of
 $G(n, p)$,
$$ \ell_1 (n, p) =~~\left\{ \begin{array}{ll}   O(\log n), ~~&
~{\rm if}~ ~ \limsup_{n\ra \infty} p(n-1) <1
 \\  (1+o(1)) \th_\gl n, & ~{\rm if}~~
\lim_{n\ra \infty} pn = \gl >1,
  \end{array}\right. $$
  where $\th_\gl$ is the positive solution of the equation
  $  1-\th-e^{-\gl\th}=0$.

Why does the size of the largest component change so dramatically
around $\gl=1$? It was Karp \cite{Kar}  who nicely explained the
reason. To find a component $C(v)$ of a fixed vertex $v$ of
$G(n,p)$, one may first expose the vertices that are adjacent to
$v$, and keep repeating the same procedure by taking each of those
adjacent vertices: Initially, $v$ is active and all other vertices
are neutral. At each step, we take an active vertex $w$ and expose
all neutral vertices adjacent to $w$. This can be done  by
checking if $\{w, w'\} \in G(n,p)$  or not for all neutral
vertices $w'$. Then, activate all neutral vertices  that are
adjacent to $w$. The vertex  $w$ is no longer active, and only
non-activated neutral vertices remain neutral. The process
terminates  when there is no more active vertex. Clearly, the
process will stop after finding all the vertices in the component
containing $v$. Provided  the number of neutral vertices does not
decrease so fast, the number of newly activated vertices has a
distribution  close to that of the binomial random variable
$\bin(n-1, p)$, where

\vspace{-3mm}
$$ \pr [\bin(n-1, p) = \ell] = {n-1 \choose \ell} p^\ell
(1-p)^{n-1-\ell}. $$

\vspace{-1mm}\noindent
 Particularly, the mean of the  number is
close to $pn$. If $pn \leq 1-\gd$ for a fixed $\gd>0$, then the
process is expected to die out quickly almost every time. Thus,
all $C(v)$'s  are expected to be small. If $pn \geq 1+\gd$, then
the process may survive forever with  positive probability. Hence,
$C(v)$ can be large with positive probability; as there are
many (actually $\Theta(n)$) trials, at least one of $C(v)$'s is expected to be large. Applying
this approach to the random directed graph, Karp was able to prove
a phase transition phenomenon for the size of the largest strong
component.

Notice that, when $pn=\Theta(1)$,  the distribution of
$\bin(n-1,p)$ is very close to the Poisson distribution with
mean $\gl:=p(n-1)$. Hence, we may further expect  that the
process described above could be approximated by the Galton-Watson
branching process defined by a Poisson random variable $\poi(\gl)$
with mean $\gl$, where
$$ \pr [\poi(\gl) =\ell] = e^{-\gl} \frac{\gl^\ell}{\ell!}. $$
Generally, the Galton-Watson branching process defined by a random
variable $X$ starts with a single unisexual organism. The organism
will give birth to $X_1$ children, where $X_1$ is a random
variable with  the same distribution as $X$. The same but
independent birth process continues from  each of the children and
the grandchildren and so on,   until no more descendant exists.
(For more information regarding Galton-Watson branching processes,
one may refer \cite{AN}.) For simplicity, we  say the
Poisson($\gl$) branching process  for the Galton-Watson branching
process defined by $\poi(\gl)$.

\mn {\bf The Poisson cloning model.} To convert  the above
observation to a rigorous proof, it is needed to overcome or
bypass two main obstacles. The first one is that the degrees of
vertices of $G(n,p)$ are not exactly i.i.d Poisson random
variables. Though they have the same distribution as
$\bin(n-1,p)$, they are not mutually independent.
 For example, the
sum of all  degrees must be even as  it is twice the number of
edges. This cannot be guaranteed if the degrees are independent.
The second one is that the number of neutral vertices keeps
decreasing. Even if both obstacles do not cause substantial
differences in many cases, one needs  at least to  keep tracking
small differences for rigorous proofs. Since these kinds of small
differences occur almost everywhere in the analysis, they
sometimes make rigorous analysis significantly difficult, if not
impossible. Furthermore, the fact that the number of neutral vertices decreases not only plays a crucial role but also yields a different result in the case that one wants know more precise behaviors.

As an approach to bypass the first obstacle, we introduce the
Poisson cloning model $G_{PC} (n,p)$ for random graphs  in which
the degrees are  i.i.d Poisson random variables with mean
$\gl=p(n-1)$. Moreover, the new model is  equivalent to the
classical model $G(n,p)$ in an asymptotic sense.
 Actually, defining the model is not extremely difficult:
First take i.i.d Poisson $\gl$ random variables $d(v)$ indexed by
all vertices $v$ in $V$. Then take $d(v)$ copies, or clones, of each
vertex $v$. If the sum of $d(v)$'s is even, then we generate a
uniform random perfect matching on the set of all clones. An edge
$\{v, w\}$ is in the random graph $G_{PC}(n,p)$ if a clone of $v$
is matched to a clone of $w$ in the random perfect matching.
The resulting graph may or may not a simple graph. If
the sum is odd, one may just take a graph with a self loop. Hence,
the graph is always not simple if the sum is odd. In the next section,
the Poisson cloning model is to be
  defined with
details.

 It is also possible to extend  the model to uniform
hypergraphs, where a $k$-uniform hypergraph on the vertex set $V$
is a collection of subsets of $V$ with size $k$. A graph is then a
$2$-uniform hypergraph. In the binomial model $H(n,p\,;k)$ for
random $k$-uniform hypergraphs, each of ${n \choose k}$ edges is
 in the hypergraph with probability $p$, independently of other edges.  The Poisson
cloning model for random $k$-uniform hypergraphs may be similarly defined and is denoted by
$H_{PC}
 (n,p\,;k)$.

The following theorem shows that the new model is essentially
equivalent to the binomial model.
\begin{thm} \label{pc}
Suppose $k\geq 2$ and $ p= \Theta (n^{1-k})$.  Then, for any
collection ${\cal H}$ of $k$-uniform simple hypergraphs,
 $$ c_{_1}\! \pr [
H_{PC} (n,p\,;k) \in {\mathcal{H}}] \leq   \pr [ H (n,p\,;k) \in
{\mathcal{H}}]
  \leq
c_{_2} \! \bb(\pr [ H_{PC} (n,p\,;k) \in
{\mathcal{H}}]^{{\frac{1}{k}}} + e^{-n}\bb),
 $$
where $$c_{_1}= k^{1/2}e^{\frac{p}{n}{k \choose 2}{n\choose
k}+\frac{p^2}{2}{n\choose k}}+O(n^{-1/2}), ~~c_{_2}=
\bb(\frac{k}{k-1}\bb) \bb(c_{_1}(k-1)\bb)^{1/k}+o(1),
$$ and  $o(1)$ goes to $0$ as $n$ goes to infinity.
\end{thm}

To overcome the second obstacle, we present an algorithm,  called
the cut-off line algorithm, that enables us to generate the
Poisson cloning model and analyze problems simultaneously. As a
consequence, the size of the largest component of $G_{PC} (n,p)$
can be described very precisely. It is also possible to analyze
the size of the $t$-core of the random hypergraph $H_{PC}
(n,p\,;k) $, where the $t$-core of a hypergraph  is the largest
subhypergraph with minimum degree at least $t$.

\mn {\bf The emergence of the giant component.} After the phase
transition result of Erd\H{o}s and
 R\'enyi, it has remained to determine the size of the largest
 component when $pn\ra 1$. Though Erd\H{o}s and
 R\'enyi suggested that the size $\ell_1 (n,p)$ of the largest component
 could be only    $O(\log
 n)$, $\Theta(n^{2/3})$, or $\Theta (n)$,
 Bollob\'as \cite{BB84a} showed  that $\ell_1 (n,p)$
 increases rather continuously by
estimating it quite accurately for $pn-1\geq n^{-1/3} \sqrt{\log
n\,}  /2 $. Later {\L}uczak \cite{TL90c} was able to estimate
$\ell_1 (n,p)$
 for $pn-1 \gg n^{-1/3}.$

 In statements in theorems and lemmas, etc.,  of this paper, we use the
following convention.

 \mn {\bf Convention:} When we say that a statement is true for all
$\ga$ in the range $a\ll \ga \ll b$, it actually means that there
is (small) constant $\eps>0$ so that the statement is true for
$\ga$ in the range $a/\eps \leq \ga \leq \eps b$.

 \begin{thm} \label{size} \cite{TL90c}
  {\rm  (Supercritical Phase)} Suppose $\gl=\gl(p,n)= 1+\eps$ with $\eps\gg n^{-1/3}
  $.
  Then, for large enough $n$, with probability at least
  $1-7 (\eps^3 n/8)^{-1/9}$,
  $$ \bb| \ell_1 (n, p) - \th_{\gl} n\bb| \leq \frac{n^{2/3}}{5} , $$
   and all
  other components are
  smaller than $n^{2/3}$.
  \end{thm}

Using estimations for the number of connected graphs with certain
numbers of  vertices and edges, and the first and second moment
methods, one may also obtain the following result for the
subcritical phase.
\begin{thm} {\rm  (Subcritical Phase)} Let $\gl(n,p)=1-\eps$ with $n^{-1/3}
  \ll \eps\ll 1$.
  Then, for any positive constant $\gd \leq 1/3$ and large
  enough $n$, with probability at least $1-(\frac{8}{\eps^3
  n})^{\gd/4}$,
$$ \bb| \ell_1 (n, p) -  \frac{2\log(\eps^3 n )}{\eps^{2}}
\bb| \leq \frac{\gd \log (\eps^3 n)}{\eps^{2}}  .
$$
  \end{thm}
There have been many results  regarding the structure of the
largest component too, for which readers may refer \cite{JLR,
TL90c, JKLP, LPW} and references therein.


For Poisson branching processes, a duality principle has been
known. A pair $(\mu, \gl)$ with $\mu< 1 < \gl$ is called a {\em
conjugate pair} if $\mu e^{-\mu} =\gl e^{-\gl}$. It is easy to see
that $\mu=(1-\thl)\gl$ for a conjugate pair $(\mu,\gl)$.  For a
conjugate pair $(\mu, \gl)$, the distribution of the
Poisson($\gl$) branching process conditioned that the process dies
out is exactly the same as that of  the Poisson($\mu$) branching
process. (See e.g. \cite{AS}, p164.)  A similar but a little bit coarse duality was
observed for the random graph $G(n,p)$ and $G(n^*,p)$ with
$\gl=\gl(n,p)>1 $ and $n^*=(1-\thl)n$.  Notice that $1-\thl$ is
the extinction probability for the Poisson($\gl$) branching
process. It has been known that the component sizes  of $G(n^*,p)$
and those of $G(n,p)$ excluding the largest component are the same
in an asymptotic sense (see \cite{AS}).

The Poisson cloning model $G_{PC}(n,p)$ equipped with the cut-off
line algorithm enables us to not only estimate $\ell_1 (n, p) $
more accurately but also establish a precise discrete duality principle:
 In the supercritical phase
$\gl:=\gl(n,p)=1+\eps$ with
 $n^{-1/3}
  \ll \eps\ll 1$, $G_{PC} (n,p)$ can be decomposed into three vertex
    disjoint graphs $C$, $S$ and $G$ whp (with high probability),
  where $C$ is a connected graph  of size
  about $\thl n$, $S$ is a smaller graph of size about $\eps^{-2} \ll \thl
  n$,   and $G$ has the same distribution as $G_{PC} (n^*,p^*)$ with
  $n^* \approx (1-\thl)n$ and $p^* \approx p$, which yields
  $\gl(n^*,p^*)\approx \mu:= (1-\thl)\gl$.
   In the subcritical phase $\gl=1-\eps$ with $n^{-1/3}\ll \eps \ll
1$, the largest component is of size
$$ \frac{\log (\eps^3 n) -2.5 \log \log (\eps^3 n)+O(1)}{-(\eps
+\log (1-\eps) ) } $$ whp. The precise statements are as follows.
We concentrate on the cases $\eps \ll 1 $ for which more careful
analysis is required. It is believed that the proofs are easily
modified for the cases of positive constants $\eps$.

\begin{thm} \label{mainpc}

\noindent {Supercritical Phase:} Let $ \gl:=\gl(n,p) =1+\eps$ with
 $n^{-1/3}
  \ll \eps\ll 1$,  $\mu:=(1-\thl)\gl$ and  $1\ll \ga \ll
(\eps^3 n)^{1/2}$. Then,
  with probability $1-e^{-\Omega (\ga^2)}$,
$G_{PC} (n,p)$ may be decomposed  into  three vertex disjoint
graphs $C$, $S$ and $G$, where $C$ is connected and
 $$ \thl n
  -\ga(n/\eps)^{1/2} \leq |C| \leq \thl n
  + \ga(n/\eps)^{1/2}, $$ and $ |S| \leq \frac{\ga^2}{\eps^2 }, $
and $G$ has the same distribution as $G_{PC}(n^*,p^*)$ for some
$n^*$ and $p^*$ satisfying
$$ (1-\thl)n - \ga (n/\eps)^{1/2} \leq n^* \leq (1-\thl)n + \ga
(n/\eps)^{1/2} , $$ and $$ \mu - \ga (\eps n)^{-1/2} \leq
\gl(n^*,p^*) \leq \mu + \ga (\eps n)^{-1/2}. $$

\mn Subcritical Phase: Suppose $ \gl:=\gl(n,p) =1-\eps$ with
 $n^{-1/3}
  \ll \eps\ll 1$. Then, the size $\ell_{1}^{PC}(n,p)$ of the largest component
of $G_{PC} (n,p)$ satisfies
$$ \pr\bb[ \ell^{PC}_1 (n,p)
\geq  \frac{\log (\eps^3 n) -2.5 \log \log (\eps^3 n)+c}{-(\eps
+\log (1-\eps) ) } \bb] \leq  2e^{-\Omega(c)}, $$ and
$$ \pr\bb[ \ell_1^{PC} (n,p)
\leq  \frac{\log (\eps^3 n) -2.5 \log \log (\eps^3 n)-c}{-(\eps
+\log (1-\eps) ) } \bb] \leq  2e^{-e^{\Omega(c)}},$$ for any
positive constant $c> 0$.

\mn Inside Window: Suppose $ \gl:=\gl(n,p) =1+\eps$ with
 $|\eps| = O(n^{1/3})$. Then, whp,
 $$\ell_{1}^{PC} (n,p) = \Theta (n^{2/3}). $$
 (All constants in $\Omega(\cdot)$'s do not depend on any of
$\eps$, $\ga$ and $c$.)
\end{thm}

By Theorem \ref{pc}, a corollary regarding  $G(n,p)$ follows.

\begin{cor} \label{main} {Supercritical region:}
Suppose $ \gl =\gl(n,p)  =1+\eps$ with
 $n^{-1/3}
  \ll \eps\ll 1$, and  $1\ll \ga \ll
(\eps^3 n)^{1/2}$. Then, in $G (n,p)$,
$$ \pr [\, |\ell_1 (n,p) -
\thl n| \geq \ga (n/\eps)^{1/2} ] \leq 2e^{-\Omega(\ga^2)}. $$
 Moreover, for the size $\ell_2 (n,p)$ of the second largest component and
 $\eps^*= 1-(1-\thl)\gl$,
$$ \pr\bb[ \ell_2 (n,p)
\geq  \frac{\log ((\eps^*)^3 n) -2.5 \log \log ((\eps^*)^3
n)+c}{-(\eps^* +\log (1-\eps^*) ) } \bb] \leq  2e^{-\Omega(c)},
$$ and
$$ \pr\bb[ \ell_2 (n,p)
\leq  \frac{\log ((\eps^*)^3 n) -2.5 \log \log ((\eps^*)^3
n)-c}{-(\eps^* +\log (1-\eps^*) ) } \bb] \leq
2e^{-e^{\Omega(c)}},$$ for any positive constant $c> 0$.

\mn {Subcritical region}: Suppose $\gl=1-\eps$ with
 $n^{-1/3}
  \ll \eps\ll 1$, then,
$$ \pr\bb[ \ell_1 (n,p)
\geq  \frac{\log (\eps^3 n) -2.5 \log \log (\eps^3 n)+c}{-(\eps
+\log (1-\eps) ) } \bb] \leq  2e^{-\Omega(c)}, $$ and
$$ \pr\bb[ \ell_1 (n,p)
\leq  \frac{\log (\eps^3 n) -2.5 \log \log (\eps^3 n)-c}{-(\eps
+\log (1-\eps) ) } \bb] \leq  2e^{-e^{\Omega(c)}}, $$ for any
positive constant $c> 0$.

\mn Inside Window: Suppose $ \gl:=\gl(n,p) =1+\eps$ with
 $|\eps| = O(n^{1/3})$. Then, whp,
 $$\ell_{1} (n,p) = \Theta (n^{2/3}). $$
\end{cor}

\mn {\bf The emergence of the $t$-core.} There are at least two
possible directions to extend  the problem of connected components.
Observing that the minimum degree in a component must be larger
than or equal to $1$, one may consider subgraphs with minimum
degree at least $t\geq2$. For a graph $G$, the $t$-core is the
largest subgraph with minimum degree at least $t$. As the minimum
degree of the union of two subgraphs is at least the smaller
minimum degree of the two, the $t$-core of a graph is unique. It
is also easy to see that the $t$-core must be an induced subgraph.
For this reason, the $t$-core of $G$ sometimes refers to its
vertex set. Denoted by $V_t (G)$ is (the vertex set of) the
$t$-core of $G$. As the $1$-core $V_1(G)$ is the set of all
non-isolated vertices, we consider the cases $t\geq 2$ throughout
this paper. If there is no subgraph with minimum degree $t$, the
$t$-core is defined to be empty.

Another direction is to consider the $t$-connectivity, where a
graph is $t$-connected if the graph remains connected after any
$t-1$ vertices are removed. Higher orders of connectivity have
been used to understand various structures of graphs. Clearly, if
a non-empty subgraph  is $t$-connected, then its minimum degree
must be $t$ or larger.

In 1984, Bollob\'as \cite{BB84b} initiated the study of $t$-core,
$t\geq 2$, and observed that, provided  $t\geq 3$ and $pn$ is
larger than a fixed constant, the $t$-core of $G(n,p)$ is
non-empty and $t$-connected whp. {\L}uczak \cite{TL91d} proved
that, for $t\geq 3$, there is an absolute constant $c$ such that
the $t$-core of $G(n,p)$
 is either empty, or larger than $cn$ and
$t$-connected,  whp. In particular, as far as the random graph
$G(n.p)$ is concerned, the $t$-core problem is the same as the
$t$-connectivity problem. Moreover, if $\gl(n,p)$ is less than
$1$, then the $t$-core of $G(n,p)$ is empty whp as the size of the
largest component is $O(n^{2/3})$ whp. As $p$ increases while $n$
is fixed, the probability of the $t$-core of $G(n,p)$ being
non-empty keeps increasing. Let $p_{_t} (n, \gd) $ be the infimum
of all $p$ that makes the probability larger than or equal to a
constant  $\gd$ in the range $0< \gd <1$. Then Bollob\'as's result
implies that $n p_{_t} (n,\gd) $ is bounded from above by a
constant. Though $n p_{_t} (n,\gd) $ may still have no limit value
as $n$ goes to infinity, it seems to be more natural to expect
that the limit exists. Furthermore, as it happens often in phase
transition phenomena, the limit, if exists, is also expected to be
independent of $\gd$. In other words, the phase transition is
expected to be sharp.

For $t=2$, the $2$-core of a graph $G$ is non-empty if and only if
$G$ contains a cycle. It is easy to see by the first moment method
that $G(n,p)$ with $p=o(1/n)$ does not contain a cycle whp. For a
constant $c$ in the range $0<c<1$, $G(n,p)$ with $pn=c$ may or may
not have a cycle with positive probability. Particularly,  the
phase transition for the existence of non-empty $2$-core is not
sharp. In the graph process $(G(n,m))_{m=0,1...}$, in which a
random edge is added one by one without repetition, Janson
\cite{SJ87} found the limiting distribution for the length of the
first cycle, especially he showed that the length  is bounded whp.
However, the expectation of the length is known to be $\Theta
(n^{1/6})$ due to Flajolet et al. \cite{FKP}. The two facts are
not contradicting each other, since there are random variables
that are bounded whp, but their expectations are  not. For
example, $\pr[X=1]=1-1/n$ and $\pr[X=n^2] =1/n$.

Bollob\'as \cite{BB84b} proved that, if $t\geq 5$ and $\gl
(n,p):=p(n-1) \geq \max\{ 67, 2t+6\}$, then $G(n,p)$ has a
non-empty $t$-core.  Chv\'atal \cite{VC} introduced the notion of
critical $\gl_t$, without proving existence,  satisfying, as $n$
goes to infinity,
$$ \pr
\Big[ \mbox{$G(n,p)$ has a non-empty $t$-core }\Big] \ra
\case{0}{if $\gl (n,p)  < \gl_{t}-\gd $}{1}{if $\gl (n,p)
>\gl_{t}+\gd$,}
$$
 for any constant $\gd>0$. He also proved
  $\gl_3 \geq 2.88$, if exists, and  claimed that $\gl_4 \geq
4.52$ and $\gl_5 \geq 6.06$ etc. could be proven  by the same
method. It is Pittel, Spencer and Wormald \cite{PSW} who proved a
more general theorem that implies that $\gl_t$ exists for fixed
$t\geq 3$. They identified the values too. We present a slightly
weaker version of the theorem.

For a Poisson random variable $\poi(\rho)$ with mean $\rho$, let
$P(\rho, i)=\pr[\poi(\rho)=i]$ and $Q(\rho, i)=\pr[\poi(\rho)\geq
i]$, i.e.,
$$ P(\rho,i) = e^{-\rho}\frac{\rho^i}{i!}, \and
Q(\rho, i) = \sum_{j=i}^{\infty} P(\rho,j)=
 e^{-\rho}\sum_{j=i}^{\infty} \frac{\rho^j}{j!},  $$
 and let $$
 \gl_t = \min_{\rho>0} \frac{\rho}{Q(\rho,t-1)}. $$

\begin{thm}
\label{PSW} Let $t\geq 3, \gl(n,p)=p(n-1)$. Then
$$ \pr \Big[ \mbox{$G(n,p)$ has a non-empty $t$-core }\Big] \ra \case{0}{if
$\gl (n,p) < \gl_{t}-n^{-\gd}$}{1}{if $\gl (n,p)
>\gl_{t}+n^{-\gd}$,}
$$ for any $\gd \in (0,1/2)$,  and the  $t$-core
when $\gl(n,p) > \gl_t+n^{-\gd}$ has  $(1+o(1)) Q(\thl \gl,t) n $
vertices, whp, where $\thl$ is the largest solution for the
equation
$$ \th - Q(\th\gl, t-1)=0. $$
\end{thm}

There have been much studies about the $t$-cores of various types
of random graphs and random hypergraphs. Fernholz and Ramachandran
\cite{FRa,FRb} have studied random graphs conditions on given
degree sequences. Cooper \cite{CC} found the critical values for
the $t$-cores of uniform multihypergraphs with given degree
sequences that include the random $k$-uniform hypergraph
$H(n,p\,;k)$. Molloy \cite{MM} has considered cores for random
hypergraphs and random satisfiability problems for Boolean
formulas.  Recently, Janson and M. J. Luczak \cite{JL} also gave
seemingly simpler proofs for $t$-core problems that cover the
result of Pittel, Spencer and Wormald. For more information and
techniques used above mentioned papers, readers may refer
\cite{JL}.


Using the Poisson cloning model for random hypergraphs together
with the cut-off line algorithm, we  completely analyze the
$t$-core problem for the random uniform hypergraph. We also
believe that the cut-off line algorithm can be used to analyze the
$t$-core problem for random hypergraphs conditioned on certain
degree sequences as in \cite{CC, FRa, FRb,JL}. As the $2$-core of
$G(n,p)$ behaves  quite differently from the other $t$-cores of
$H(n,p\,; k)$, we exclude the case $k=t=2$. The case requires much
more careful analysis and will be studied in a subsequent paper.

The critical value for the problem turns out to be the minimum
$\gl$ such that there is a positive solution for the equation
 \beq{crteq} \th- Q(\thk \gl, t-1)=0. \eeq
It is not difficult to check that  the minimum is
 \beq{rho}  \gl_{\crt} (k,t)
:=\min_{\rho>0} \frac{\rho}{ Q(\rho,t-1)^{k-1}}.\eeq For  $\gl>
\gl_{\crt} (k,t)$, let $\thl$ be the largest solution for the
equation $\th^{\frac{1}{k-1}}- Q(\th \gl, t-1)=0.$


\begin{thm}
\label{sizeofcore} Let  $k,t\geq 2 $, excluding $k=t=2$, and $\gs
\gg n^{-1/2}$.

\mn Subcritical Phase: If $\gl(n,p\,; k):= p {n-1\choose k-1}
=\gl_{\crt} -\gs$ is uniformly bounded from below by $0$ and
$i_{_0}(k,t)$ is the minimum $i$ such that ${i \choose k} \geq
ti/k$, then
$$ \pr  [ V_t (H(n,p\,;k)) \not=\emptyset ] \leq 2e^{-\Omega (\gs^2 n) }+
O ( n^{-(t-1-t/k) i_{_0}(k,t)}),  $$ and, for any
 $\gd>0$,
 \bq{subsize}{\pr  [ |V_t (H(n,p\,;k)) |\geq  \gd n  ] \leq 2e^{-\Omega (\gs^2 n)
}+2e^{-\Omega(\gd^{2k/ (k-1))} n)}.}

\mn Supercritical Phase: If $\gl= \gl(n,p\, ;k) =\gl_{\crt} +\gs $
is uniformly bounded from above, then, for all $\ga $ in the range
 $1\ll \ga  \ll \gs n^{1/2} $,
  $$ \pr [\, \, |\, |V_t(n,p\, ; k)|- Q (\thl \gl,t) n | \geq
 \ga  ( n/\gs)^{1/2}\,  ] = e^{-\Omega (\ga ^2 ) },$$
 and, for any  $i\geq t$ and the sets $V_t (i)$ (resp.  $W_t
 (i)$)
  of vertices of degree $i$ (resp.   larger than or equal to $i$) in the $t$-core,
$$ \pr \bbr \bb|  |V_t (i)| -  P(\thl \gl, i) n
\bb|\geq  \gd  n   \bb] \leq 2e^{-\Omega(\min\{ \gd^2 \gs n, \gs^2
n\} )},
$$
and
$$ \pr \bbr \bb|  |W_t (i)| -  Q(\thl \gl, i) n
\bb|\geq  \gd   n   \bb] \leq 2e^{-\Omega(\min\{ \gd^2 \gs n,
\gs^2 n\} )}.
$$
In particular, for $\rho_{_\crt}:= \th_{\gl_{\crt} (k,t)}
\gl_\crt(k,t)$,
$$   |V_t (i)| = (1+O(\gs^{1/2}))^i   P(\rho_{_\crt}, i) n
+O\bb( (n/\gs)^{1/2} \log n\bb),
$$
with probability $1- 2e^{-\Omega(\min\{\log^2 n , \gs^2 n\} )}$.
Moreover, if all $|V_t (i)|$, $i\geq t$, are given, each simple
graph with the degree sequence induced by $ |V_t (i)|$, ${i\geq
t}$, is equally likely to be the $t$-core.
\end{thm}

\old{
  \bq{degree}{ \pr \bb[\, |\, |V_t (l)| - P (\gt_\gl^{k-1}\gl;l)n| \geq
\ga (n/\gs)^{1/2} \bb] = e^{-\Omega (\ga ^2 ) } . } Generally, for
all $l\geq t$ with $P (\gt_\gl^{k-1}\gl;l)n\ra \infty$ and $1\ll
\ga \ll \min\{ \gs n^{1/2} , (P (\gt_\gl^{k-1}\gl;l)n)^{1/2}, (\gs
n)^{1/2}/l \}
  $,
 \bq{degree2}{ \pr \bb[\, \, \bb|\, |V_t (l)| - P (\gt_\gl^{k-1}\gl;l)n\bb| \geq
\ga \bb(1+ \gs^{-1/2}l P^{1/2} (\gt_\gl^{k-1}\gl;l)\bb) \bb(P
(\gt_\gl^{k-1}\gl;l) n\bb)^{1/2}  \bb] = e^{-\Omega (\ga ^2 ) } ,
} and for all $l\geq t$ with  $P (\gt_\gl^{k-1}\gl;l)n \leq  c$,
for some $c=c(n) >0$ uniformly bounded from above, and $ 1\ll
-\log c^r/r! \ll \min\{\gs^2 n, \gs n/l^2\} $,
$$ \pr [ \, |V_t(l)| \geq r ] \leq  ((1+o(1)) c)^r/r!.
$$
Moreover, if $|V_t (l)|$ are given for all $l\geq t$, then the
$t$-core of $H_{PC} (\gl )$
 is isomorphic to the
 cloning model conditioned on a degree sequence with $|V_t (l)|$
 vertices of degree
$l$. } 

As one might guess, we will prove a stronger theorem (Theorem
\ref{stop}) for the Poisson cloning model $H_{PC} (n,p\, ; k)$,
from which Theorem \ref{sizeofcore} easily follows.

\bn {\bf The pure literal algorithm for the random $k$-SAT
problem} Recently the satisfiability problem for Boolean formulas
has played a central role in the theory of computational
complexity. An instance of the problem is a formula given by a
conjunctive normal form (CNF), that is, a conjunction of
disjunctions. Each disjunction, or clause, is of the form
$(y_1\vee \cdots \vee y_k)$, where $y_i$'s are chosen among $2n$
literals consisting of $n$ Boolean variables and their negations.
Given a formula, the problem is whether there exists an assignment
of the $n$ variables satisfying the formula. When such an
assignment exists,  the formula is satisfiable. Otherwise, it is
unsatisfiable. A pair of literals $y,z$ is {\em strictly distinct}
if $y$ is neither $z$ nor the negation of $z$. When the input
formula is restricted to have only clauses of $k$ pairwise
strictly distinct literals, called {\em $k$-clauses}, the problem
is called the $k$-SAT problem.

It is known that the satisfiability problem is NP-complete
\cite{Coo71}, so that determining whether an arbitrary  formula is
satisfiable or not is regarded very difficult (assuming P $\neq$
NP), in the sense that it is at least as hard as any problem whose
solutions can be verified in polynomial time. Cook \cite{Coo71}
proved that even the $k$-SAT problem for $k\geq 3$ is NP-complete
too, while the $2$-SAT problem can be solved in polynomial time.
Given these facts, researchers have tried to find heuristic
algorithms that are able to determine, in polynomial time, the
satisfiability of most of $k$-SAT formulas, especially $3$-SAT
formulas. Among others, a number of heuristic algorithms have been
considered based on  Davis-Putnam algorithm \cite{DP}. Since one
has to define``most" before showing that his or her algorithm
works for most $k$-SAT formulas, random models for the $k$-SAT
problem have been introduced.

The most common models for the random $k$-SAT problems is the
uniform model $F_k (n,m)$ and $\fnpk$. Here,  $F (n,m\,; k)$ is
sampled uniformly at random from the set of all $k$-SAT formulas
with $n$ variables and  $m$ clauses. The other model  may be
constructed by selecting each $k$-clause with probability $p$
independently of all other clauses. The random formula $F(n,p\,; k
)$ is a conjunction of all selected clauses. Since there are $2^k
{n \choose k}$ $k$-clauses all together, the expected number of
clauses in the formula is $m_p:=2^k p {n \choose k}$. Other models
include one formed by selecting a uniform random $k$-clause $m$
times with replacement. It is known that theses three models are
essentially equivalent if $m_p =m$ and $m=\Theta(n)$.

Not surprisingly, the random $2$-SAT and the random $3$-SAT
problems have been most extensively studied and many results have
been published. For $k=2$, Chv\'{a}tal and Reed \cite{CR}, Goerdt
\cite{Go} and Fernandez de la Vega \cite{FV} independently proved
that the random $2$-SAT problem undergoes  a  phase transition at
$2pn =1$, that is, for $F(n,p)=F(n,p\,;k)$,
$$\lim_{n\ra \infty}
\pr[ F(n,p) ~\mbox{is satisfiable} ] = \case{1}{if $\limsup_{n\ra
\infty} 2pn <1 $}{ 0}{if $\liminf_{n\ra \infty} 2pn
>1$.} $$
Since there are $2n$ literals, $2pn$ is the right parameter.
Bollob\'as et al. \cite{BBCKW} took more sophisticated approaches
to determine the scaling window for the problem: Let $\rho \gg
n^{-1/3}$. Then
$$
\pr[ F (n,\sfrac{1-\rho}{2n}) ~\mbox{is satisfiable} ]
=1-\Theta(1/\rho^3), $$ and $$ \pr[ F(n,\sfrac{1+\rho}{2n})
~\mbox{is satisfiable} ] = e^{-\Theta(\rho^3)}
$$

 Though it is believed that the random $k$-SAT problem,
$k\geq 3$, undergoes a similar phase transition, it remains as a
conjecture. Only sharp transitions are known by a seminal result
of Friedgut \cite{FB}.

To find a satisfying assignment for a $k$-SAT formula, one may
apply the pure literal algorithm (PLA): A literal is {\em pure} in
the formula if it belongs to one or more clauses of the formula,
while its negation is in no clause. PLA keeps selecting a pure
literal, setting it  true, and removing clauses containing the
literal as they are already satisfied. It stops when there is no
more pure literal. We say that PLA {\em succeeds} if no clause
remains in the formula after it stops, and it {\em fails}
otherwise. Clearly, the formula is satisfiable if PLA succeeds.
The converse is not true, for example, $(y\vee z) \wedge
(\bar{y}\vee \bar{z})$ is satisfiable whereas no pure literal
exists. Broder, Frieze, and Upfal \cite{BFU} analyzed PLA for the
random $3$-SAT problem to show that, for $F_3 (n,p)$,  if
$\limsup\frac{m_{_p}}{n} < 1.63 $ then PLA succeeds whp, and if $
\liminf\frac{m_{_p}}{n} >1.7$ then it fails whp. Mitzenmacher
\cite{M} used the differential equation method introduced by
Wormald \cite{W} to claim, without rigorous proof, that the
threshold for PLA exists and it is the solution of certain
equations, which are somewhat complicated. That is, there is a
constant  $c_{_k}$, $k\geq 3 $, so that PLA succeeds whp if
$\limsup\frac{m_{_p}}{n} < c_{_k}$, and fails whp if $
\liminf\frac{m_{_p}}{n} > c_{_k}$. For more about upper bounds for
the random satisfiability problems, readers may refer
 \cite{BFU, DB,EF, FP, Jea,  Ka, KM,  Kea, Zi}. For more advanced algorithms
than PLA, which give various improved lower bounds, one may refer
\cite{A1, ASb, Dsur, CF, FS}. And a variation of the random
satisfiability called the $(2+p)$-SAT problem can be found in
\cite{AKKK, MZ1, MZ2}.

For rigorous analysis of PLA, we  consider  the Poisson cloning
model for the random $k$-SAT problems, $k\geq 2$. Since  a
$k$-clause can be regarded as a (hyper)edge consisting of $k$
vertices, the Poisson cloning model $F_{PC}(n,p\, ; k)$ can be
defined as $H_{PC} (2n,p\,;k)$ on $V=\{ x_{_1}, \bar{x}_{_1}, ...,
x_{_n}, \bar{x}_{_n}\}$ without  edges that contain both of a
variable and its negation. Then it is not difficult to establish
an asymptotic equivalence between $F_{PC}(n,p\, ;k)$ and $F(n,p\,;
k)$ from  Theorem \ref{pc}. The details can be found in the next
section.

\begin{thm} \label{equiv}
Suppose $ p= \Theta (n^{1-k})$.  Then, for any collection ${\cal
F}$ of $k$-SAT formulas,
 $$ c_{_1}^* \pr [
\fnpk \in \f] \leq    \pr [ F (n,p\,;k) \in \f]
 \leq
c_{_2}^* (\pr [ \fnpk \in \f ]^{{\frac{1}{k}}} + e^{-n}),
 $$
 where
$$c_{_{1}}^* = k^{1/2}e^{\frac{p}{n}{k \choose 2}{2n\choose k}
+ \frac{p^2}{2}{2n\choose k}}+o(1), ~~c_{_{2}}^* =
e^{\frac{p(1-1/k)}{2n}{k \choose 2}{2n\choose
k}}\bb(\frac{k}{k-1}\bb) \bb((k-1)c_{_{1}}^* \bb)^{1/k}+o(1).
$$
and $o(1)$ goes to $0$ as $n$ goes infinity.
\end{thm}

As mentioned above,  PLA undergoes a sharp phase transition.
Actually, for $k\geq 3$,  it turns out that the phase transition
is similar to that of the 2-core in the random hypergraphs. The
case $k=2$ is similar to the $2$-core problem of $G(n,p)$ and will
be studied in a subsequent paper. Let
$$ \gl_{\crt} (k)  :=\min_{\rho>0} \frac{\rho
    }{Q(\rho,1)^{k-1}}= \min_{\rho>0} \frac{\rho
    }{(1-e^{-\rho})^{k-1}}, $$
and, for $\gl \geq \gl_{\crt} (k)$, let $\thl$ be the largest
solution of the equation $\th^{\frac{1}{k-1}} -1+e^{-\th\gl} =0$.
Denote  by $X_R(n,p\,;k)$ is the set of variables in
$X:=\{x_{_1},..., x_{_n}\}$ whose truth values are not determined
by PLA. The remaining formula is denoted by $F_R(n,p\,;k)$. The
residual degree $d_R(z)$ of a literal $z$ of $X_R (n,p\,;k)$ is
the number of clauses in $F_R(n,p\,;k)$ containing $z$. It is easy
to see that $F_R(n,p\,;k)$ is independent of choices of pure
literals when PLA is carried out.

\begin{thm}
\label{plrthm} Let $\gl(n,p\, ;k)= p{2n-1 \choose k-1}$,
$k\geq 3$ and $\gs\gg n^{-1/2}$. 
If $\gl(n,p\, ;k)  <  \gl(k) -\gs $ is uniformly bounded from
below by $0$ and $i_{_0}(k)$ is the minimum such that $2^k {i
\choose k} \geq 2i/k$, then
$$ \pr[ X_R(n,p\,;k) \not=\0\, ] \leq 2 e^{-\Omega(\gs^2 n)}
+ O(n^{-(1-2/k)i_{_0}(k)}).$$ \mn Supercritical Phase: If $\gl=
\gl(n,p\, ;k) =\gl_{\crt} +\gs $ is uniformly bounded from above,
then, for all $\ga $ in the range
 $1\ll \ga  \ll \gs n^{1/2} $,
 $$  \pr [\, \, |\, |X_R(n,p\, ; k)|- (1-e^{-\thl\gl})^2 n | \geq
 \ga  ( n/\gs)^{1/2}\,  ] = e^{-\Omega (\ga ^2 ) },$$
 and, for any  $i,j\geq 1$ and the sets $X_R (i,j)$ (resp.  $Y_R
 (i,j)$)
  of variables $x\in X_R(n,p\,;k)$ with $d_R (x)=i, d_R(\bar{x})=j$
   (resp.   $d_R (x)\geq i, d_R(\bar{x})\geq j$),
$$ \pr \bbr \bb|  |X_R (i,j)| -  P(\thl \gl, i)P(\thl \gl, j) n
\bb|\geq  \gd  n   \bb] \leq 2e^{-\Omega(\min\{ \gd^2 \gs n, \gs^2
n\} )},
$$
and
$$ \pr \bbr \bb|  |Y_R (i,j)| -  Q(\thl \gl, i)Q(\thl\gl,j) n
\bb|\geq  \gd   n   \bb] \leq 2e^{-\Omega(\min\{ \gd^2 \gs n,
\gs^2 n\} )}.
$$
Moreover, if all $|X_R (i,j)|$, $i,j\geq 1$, are given, each
formula with the degree sequence induced by $ |X_R (i,j)|$,
${i,j\geq 1}$, is equally likely to be the residual formula $F_R
(n,p\,;k)$.
\end{thm}
Like the $t$-core problem,  we will prove a stronger theorem
(Theorem \ref{stopplr}) for the Poisson cloning model $F_{PC}
(n,p\, ; k)$, from which Theorem \ref{sizeofcore} easily follows.

\bn

 The rest of the paper is organized as follows.  In the next
section, the Poisson cloning model is defined with details.  The
cut-off line algorithm and the cut-off line lemma are presented in
Section \ref{grpm}. Section \ref{sldi} is for Chernoff type large
derivation inequalities that will be used in most of our proofs.
In Section \ref{ldi}, a generalized core is defined and  the main
lemma is presented. As the proof of Theorem \ref{mainpc} is more
sophisticated, Theorems \ref{sizeofcore} and \ref{plrthm} are
first  proved in Sections \ref{score} and \ref{sat}. Section
\ref{giant} is for the proof of Theorem \ref{mainpc}. Closing
remark follows in Section \ref{remarks}.

\section{The Poisson Cloning Model} 


\mn In this section, we define the Poisson cloning models
$G_{PC}(n,p)$
 for random graphs and  generally $H_{PC}(n,p\, ; k)$ for random hypergraphs.
Then,  Theorem \ref{pc} will be proven.

To construct $G_{{PC}}(n,p)$, we first take i.i.d Poisson
$\gl=p(n-1)$ random variables $d(v)$ indexed by vertices $v$ in
the set $V$ with $|V|=n$. Then, take $d(v)$ copies of each vertex
$v\in V$. The copies of $v$ are called {\em clones of $v$}, or
simply {\em $v$-clones}. Since the sum of Poisson random variables
is also a Poisson random variable, the total number $N_\gl :=
\sum_{v\in V} d(v)$ of clones is a Poisson $\gl n$ random
variable. It is sometimes convenient to take a reverse, but
equivalent, construction. We first take a Poisson $\gl n=2p {n
\choose 2}$ random variable $N_\gl$ and then take $N_\gl$
unlabelled clones. Each clone is to be independently labelled as
$v$-clone uniformly at random, in the sense that $v$ is chosen
uniformly at random from  $V$. It is well-known that the numbers
$d(v)$ of $v$-clones are i.i.d Poisson random variables with mean
$\gl$.

If $N_{\gl}$ is even, the multigraph $G_{PC} (n, p) $ is defined
by generating a (uniform) random perfect matching of those $N_\gl$
clones, and contracting clones of the same vertex. That is, if a
$v$-clone and a $w$-clone are matched, then  the edge $\{v, w \}$
is in $G_{PC}(n,p)$ with multiplicity. In case that $v=w$, it
produces a loop that contributes 2 in the degree of $v$. If
$N_\gl$ is odd, we may define $G_{PC}(n,p)$ to be any graph with a
special loop that, unlike other loops, contributes only 1 in the
degree of the corresponding vertex. In particular, if $N_\gl$ is
odd, $G_{PC}(n,p)$ is not a simple graph.

Strictly speaking, $G_{{PC}}(n,p)$ varies depending on how to
define it when $N_\gl$ is odd. However, if only simple graphs are
concerned, the case of $N_\gl$ being odd would not  matter. For
example, the probability that $G_{{PC}}(n,p)$ is a simple graph
with a component  larger than $0.1n$ does not depend on how
$G_{{PC}}(n,p)$ is defined when  $N_\gl$ is odd, as it is not a
simple graph anyway. Generally, for any collection $\mathcal{G}$
of simple graphs, the probability that $G_{PC} (n,p)$ is in
$\mathcal{G}$ is totally independent of how $G_{{PC}}(n,p)$ is
defined when $N_\gl$ is odd. Notice that properties of simple
graphs are actually mean  collections of simple graphs. Therefore,
when properties of simple graphs are concerned, it is not
necessary to describe $G_{{PC}}(n,p)$ for odd $N_{\gl}$.

Here are two specific ways to define $G_{PC}(n,p)$.

\vspace{-3mm}
\begin{ex}\label{ex1} One may keep matching two clones
chosen uniformly at random among all unmatched clones.
\end{ex}

\vspace{-7mm}
\begin{ex}\label{ex2} One may keep choosing his or her favorite
unmatched clone, and matching it to a clone selected  uniform at
random from all other unmatched clones.
\end{ex}

\vspace{-3mm} \noindent If $N_\gl$ is even, both examples  would
yield  uniform random perfect matchings. If $N_\gl$ odd, each of
them would yield a matching and an unmatched clone. We may create
the special loop consisting of the vertex for which the unmatched
clone is labelled. More specific ways to choose random clones will
be described in the next section,

Generally for $k\geq 3$, the Poisson cloning model $H_{PC}
(n,p\,;k)$ for $k$-uniform hypergraphs may  be  defined by the
same way: We take i.i.d Poisson $\gl=p{n-1 \choose k-1}$ random
variables $d(v)$, $v\in V$, and then take $d(v)$ clones of each
$v$. If $N_{\gl}:=\sum_{v\in V} d(v)$ is divisible by $k$, the
multihypergraph $H_{PC} (n, p;k) $ is defined by generating a
uniform random perfect matching consisting of $k$-tuples of
clones, and contracting clones of the same vertex. That is, if
$v_{_1}$-clone, $v_{_2}$-clone, ..., $v_{_k}$-clone are matched in
the perfect matching, then  the edge $\{v_{_1}, v_{_2}, ...,
v_{_k}\}$ is in $H_{PC}(n,p\, ;k)$ with multiplicity. If $N_\gl$
is not divisible by $k$, $H_{PC}(n,p\,;k)$ may be any hypergraph
with a special edge consisting of $N_{\gl} -k  \lfloor N_{\gl}/k
\rfloor$ vertices. In particular, $H_{PC}(n,p\, ;k)$ is not
$k$-uniform when $N_\gl$ is not divisible by $k$. Therefore, as
long as properties of simple $k$-uniform hypergraphs are
concerned, we do not have to describe $H_{PC}(n,p\, ;k)$ when
$N_\gl$ is not divisible by $k$.

We show that the Poisson cloning model $H_{PC}(n,p\,;k)$, $k\geq
2$, is contiguous to the classical model $H(n,p\,;k)$ when the
expected average degree is a constant.

\mn {\bf Theorem \ref{pc}} {\em {(Restated)} Suppose $k\geq 2$ and
$ p= \Theta (n^{1-k})$. Then, for any collection ${\cal H}$ of
$k$-uniform simple hypergraphs,
 $$ c_{_1} \! \pr [
H_{PC} (n,p\,;k) \in {\mathcal{H}}] \leq   \pr [ H (n,p\,;k) \in
{\mathcal{H}}]
  \leq
c_{_2} \! \bb(\pr [ H_{PC} (n,p\,;k) \in
{\mathcal{H}}]^{{\frac{1}{k}}} + e^{-n}\bb),
 $$
where $$c_{_1}= k^{1/2}e^{\frac{p}{n}{k \choose 2}{n\choose
k}+\frac{p^2}{2}{n\choose k}}+O(n^{-1/2}), ~~c_{_2}=
\bb(\frac{k}{k-1}\bb) \bb(c_{_1}(k-1)\bb)^{1/k}+o(1),
$$ and  $o(1)$ goes to $0$ as $n$ goes to infinity.
}

 \mn \pf We assume that  the random perfect
matching is generated by keeping choosing $k$ unlabelled clones
and labelling them uniformly at random, as any other way to
generate it is equivalent provided  $N_\gl$ is divisible by $k$.
Let $H$ be  a fixed simple $k$-uniform hypergraph with $m$ edges.
Then $H_{PC} (n,p\,;k) = H$ if and only if $N_\gl=km$ and the $km$
clones are labelled so that $H$ is yielded after contraction.
 The
first $k$ clones are labelled to be one of the $m$ edges with
probability $m\frac{k}{n}\frac{k-1}{n} \cdots \frac{1}{n}=m
\frac{k!}{n^k}$, and the second $k$ clones are labelled to be one
of the remaining $m-1$ edges with probability $(m-1)
\frac{k!}{n^k}$, and so on. That is,  $$ \pr[H_{PC} (n,p\,;k) = H]
= \pr[ N_\gl = km] m! \bb( \frac{k!}{n^k}\bb)^m.
$$

As $N_\gl$ is a Poisson random variable with mean $\gl n= pn{n-1
\choose k-1}=kp{n \choose k}$, we have
 \bean \pr[ N_\gl = km] &=& e^{-kp {n \choose k}} \frac{ (k p {n \choose
k})^{km}}{(km)!} \\ &=& (1+O(\sfrac{1}{m}))  e^{-kp {n \choose k}}
\frac{k^{km} p^{km} {n \choose k}^{km}}{(2\pi km)^{1/2}
(\frac{km}{e})^{km}}  \\
&=&  e^{O(\frac{1}{m})} (2\pi km)^{-1/2} \bb( e^{-p{n\choose k}}
\frac{p^m {n \choose k}^m}{(\frac{m}{e})^m}\bb)^{k} ,\eean unless
$m=0$.  Therefore,
$$ \pr[H (n,p\,;k) = H] = p^m (1-p)^{{n \choose k} -m}
= p^m e^{-p{n \choose k}-\frac{p^2}{2}{n \choose k} +pm +O(p^3n^k
+ p^2 m))}
$$
implies that
 \bean
\frac{\pr[ H_{PC} (n,p\,;k) = H] }{\pr[H (n,p\,;k) = H] } &=&
e^{-pm +\frac{p^2}{2}{n \choose k} +O(p^3 n^k + p^2 m+1/m)}
\\
& & \times ~(2\pi km)^{-1/2}  m! \bb(\frac{k!}{n^k}\bb)^m {n
\choose k}^m \bb(\frac{e}{m}\bb)^m \bb( e^{-p{n\choose k}}
\frac{p^m {n \choose k}^m}{ (\frac{m}{e})^m}\bb)^{k-1}   .\eean
Since
$$
(2\pi m)^{-1/2}  m! \bb(\frac{e}{m}\bb)^m = 1+O(\sfrac{1}{m})=
e^{O(1/m)}, \and \frac{k!}{n^k} {n \choose k} = e^{-
(1+\Theta(1/n)) {k \choose 2}/n },
$$
we conclude that
$$
\frac{\pr[ H_{PC} (n,p\,;k) = H] }{\pr[H (n,p\,;k) = H] } =
k^{-1/2}  \bb( e^{-p{n\choose k}} \frac{p^m {n \choose k}^m}{
(\frac{m}{e})^m}\bb)^{k-1} e^{ - (1+\Theta(\frac{1}{n})){k \choose
2}\frac{m}{n}-pm+\frac{p^2}{2}{n \choose k} +O(p^3n^k + p^2m+1/m)}
.$$ For the case  $m=0$, we replace $1/m$ in the last $O(\cdot)$
term by $\frac{1}{m+1}$. Then,  it is easy to see
 \beq{hpch}
\frac{\pr[ H_{PC} (n,p\,;k) = H] }{\pr[H (n,p\,;k) = H] } =
k^{-1/2}  \bb( e^{-p{n\choose k}} \frac{p^m {n \choose k}^m}{
(\frac{m}{e})^m}\bb)^{k-1} e^{ - (1+\Theta(\frac{1}{n})){k \choose
2}\frac{m}{n}-pm+\frac{p^2}{2}{n \choose k} +O(p^3n^k +
p^2m+\frac{1}{m+1})} ,\eeq for all $m$ in the range $0\leq m \leq
{n \choose k}$.
 Let $R_m = e^{-p{n\choose k}} \frac{p^m {n \choose k}^m}{
(\frac{m}{e})^m}$, or equivalently $R_m =e^{-\frac{\gl n}{k}}
(\frac{e\gln}{ km})^m$ by $\gl n = k p{n \choose k}$. Then
$\frac{R_{m+1}}{R_m} = (1+O(\frac{1}{m}))\frac{ \gln  }{k m}$.
This gives that $R_m$ has its maximum $1+O(\frac{1}{ n})$ when $m=
\frac{\gln}{k} +O(1)$, assuming $\gl =\Theta (1)$, or $p=\Theta
(n^{1-k})$. Moreover, it is not difficult to show that
$$ R_m = (1+O(\sfrac{1}{\gl n})) e^{-\frac{
 (1+O(\frac{m}{n}))(km-\gl n )^2 }{2\gln} }~~\mbox{if
~~$|m- \frac{\gln}{k} | \leq  n$}, \and R_m \leq e^{-\Omega {(|km-
\gln |)} }~~\mbox{otherwise.} $$ Hence
$$
\frac{\pr[ H_{PC} (n,p\,;k) = H] }{\pr[H (n,p\,;k) = H] } \leq
k^{-1/2} e^{-{k \choose 2} \frac{\gl}{k}- \frac{p\gln}{2k}}
+O(n^{-1/2})= k^{-1/2} e^{-\frac{p}{n}{n \choose k} {k \choose
2}-\frac{p^2}{2}{n \choose k} }+O(n^{-1/2})=:c_{_1}^{-1},$$ which
yields
$$ \pr [ H(n,p\,; k) =H ] \geq c_{_1} \pr [ H_{PC} (n,p\,; k) =H]
. $$ Thus,
 \bean  \pr [ H(n,p\,; k) \in \h] &=& \sum_{H\in\h}  \pr [ H(n,p\,; k) =H
 ]\\
&\geq & \sum_{H\in\h} c_{_1}\pr [ H_{PC} (n,p\,; k) =H] =c_{_1}
\pr [ H_{PC} (n,p\,; k) \in \h]. \eean

For the upper bound, take the minimum $m_{_1}\geq p{n \choose k}$
such that \beq{rm1} e^{-(m_{_1}-p{n \choose k})({k\choose 2}/n
+p)} R_{m_{_1}} \leq \bb(c_{_1} (k-1) \pr [ H_{PC}(n,p\,;k)\in
\h]\bb)^{1/k}+ e^{-n }. \eeq It is routine to check that
 $   m_{_1} =\Theta(n) $. Let $$\h_1 = \{ H\in \h: \mbox{the
 number of edges in $H$ is at least $p{n \choose k}$}\}. $$
 Then
  \bean \pr
[ H(n,p\,;k)\in \h_1] &=& \sum_{m\geq m_{_1}} \sum_{H\in \h \atop
|H|=m} \pr [ H(n,p\,;k)=H] + \sum_{m: \atop p{n \choose k}\leq m<
m_{_1}}\sum_{H\in \h \atop |H|=m} \pr [
H(n,p\,;k)=H]   \\
&\leq& \pr\bb[\, \, \bin \bb({n\choose k}, p\bb)\geq m_{_1} \bb] +
\sum_{m: \atop p{n \choose k}\leq m< m_{_1}}\sum_{H\in \h \atop
|H|=m} \pr [ H(n,p\,;k)=H] . \eean For $m> m_{_1}$,
$$ \frac{\pr\bb[\, \, \bin \bb({n\choose k},
p\bb)=m \bb]}{\pr\bb[\, \, \bin \bb({n\choose k}, p\bb)=m-1 \bb]}
\leq  (1+O(p)) \frac{ p{n \choose k}}{m} $$ implies that
  \bq{m1a}
{
 \pr\bb[\, \, \bin \bb({n\choose k}, p\bb)\geq m_{_1} \bb] \leq
(1/2+o(1)) (2\pi m_{_1})^{1/2} \pr\bb[\, \, \bin \bb({n\choose k},
p\bb)=m_{_1} \bb]
 }
and, if $m_{_1}- p{n \choose k} \gg n^{1/2}$,
  \bq{m1b}
   {
   \pr\bb[\, \,
\bin \bb({n\choose k}, p\bb)\geq m_{_1} \bb] \ll
 m_{_1}^{1/2} \pr\bb[\, \, \bin \bb({n\choose k},
p\bb)=m_{_1} \bb].
}
  Observe that   \bean  \pr\bb[\, \, \bin
\bb({n\choose k}, p\bb)=m_{_1} \bb] &=& { {n \choose k} \choose
m_{_1} } p^{m_{_1}} (1-p)^{{n \choose
k}-m_{_1}} \\
 &\leq&  (1+o(1)) e^{-p{n \choose k}} \frac{p^{m_{_1}} {n\choose
 k}^{m_{_1}}}{(2\pi m_{_1})^{1/2} (\frac{m_{_1}}{e})^{m_{_1}}}
 e^{-{m_{_1} \choose 2}/{n \choose k} + pm_{_1} - \frac{p^2}{2}{n
 \choose k} + O(1/n)} \\
 &=& (1+o(1)) (2\pi m_{_1})^{-1/2} R_{m_{_1}}
 e^{-{m_{_1} \choose 2}/{n \choose k} + pm_{_1} - \frac{p^2}{2}{n
 \choose k} + O(1/n)} ,
\eean as $R_{m_{_1}} = e^{-p{n \choose k}} \frac{p^{m_{_1}}
{n\choose
 k}^{m_{_1}}}{ (\frac{m_{_1}}{e})^{m_{_1}}}$. Since
 $$ \mbox{$-{m_{_1} \choose 2}/{n \choose k} + pm_{_1} - \frac{p^2}{2}{n
 \choose k} = -\frac{1}{2{n \choose k}} (m_{_1}- p{n\choose k})^2
 + O(1/n)$,} $$  \raf{rm1} together with \raf{m1a} and \raf{m1b} gives
 \beq{nkpm} \pr\bb[\, \, \bin \bb({n\choose k}, p\bb)\geq m_{_1} \bb]
 \leq  (1/2+o(1))\bb(  \bb(c_{_1} (k-1) \pr [ H_{PC}(n,p\,;k)\in
\h]\bb)^{1/k}+e^{-n }\bb). \eeq

 For $m$ in the range $p{n \choose k} \leq m <
 m_{_1}$,
 $$  e^{-(m_{_1}-p{n \choose k})({k\choose 2}/n +p)}
 R_m \geq \bb(c_{_1} (k-1) \pr [ H_{PC}(n,p\,;k)\in
\h]\bb)^{1/k}. $$ Provided  $H$ has  $m$ edges,  \raf{hpch} yields
 \bean
\pr[H (n,p\,;k) = H]  &=&  (1+o(1)) k^{1/2}  R_m^{1-k}  e^{ {k
\choose 2}\frac{m}{n}+pm-\frac{p^2}{2}{n \choose k} }\pr[ H_{PC}
(n,p\,;k) = H]  \\
& \leq & (1+o(1)) k^{1/2} e^{\frac{p}{n}{n \choose k}{k \choose
2}+\frac{p^2}{2} {n \choose k}} \bb(c_{_1} (k-1) \pr [
H_{PC}(n,p\,;k)\in \h]\bb)^{\frac{1-k}{k}} \\
& & ~~\times~ \pr[ H_{PC} (n,p\,;k) = H]  .\eean Hence
 \bean
\sum_{m: \atop p{n \choose k}\leq m< m_{_1}}\sum_{H\in \h \atop
|H|=m} \pr [ H(n,p\,;k)=H] &\leq&
 (c_{_1}+o(1)) \bb(c_{_1} (k-1) \pr [
H_{PC}(n,p\,;k)\in \h]\bb)^{\frac{1-k}{k}} \\
 & & ~\times  ~\pr[ H_{PC}
(n,p\,;k) \in\h_1]. \eean This together with  \raf{nkpm} gives
 \bean
 \pr[ H
(n,p\,;k) \in\h_1] &\leq & (1/2+o(1))\bb(  \bb(c_{_1} (k-1) \pr [
H_{PC}(n,p\,;k)\in \h]\bb)^{1/k}+e^{-\frac{p}{10}{n \choose k}
}\bb) \\
& & + (c_{_1}+o(1)) \bb(c_{_1} (k-1) \pr [ H_{PC}(n,p\,;k)\in
\h]\bb)^{\frac{1-k}{k}} \pr[ H_{PC} (n,p\,;k) \in\h_1] \eean

Similarly, for the maximum $m_{_2}<  p{n \choose k}$ such that
 $$ R_{m_{_2}} \leq \bb(c_{_1} (k-1) \pr [ H_{PC}(n,p\,;k)\in
\h]\bb)^{1/k}+ e^{-n },$$ and
  $$\h_2 = \{ H\in \h: \mbox{the
 number of edges in $H$ is less than  $p{n \choose k}$}\}, $$
 we have
 \bean
 \pr[ H
(n,p\,;k) \in\h_2] &\leq & (1/2+o(1))\bb(  \bb(c_{_1} (k-1) \pr [
H_{PC}(n,p\,;k)\in \h]\bb)^{1/k}+e^{-n
}\bb) \\
& & + (c_{_1}+o(1)) \bb(c_{_1} (k-1) \pr [ H_{PC}(n,p\,;k)\in
\h]\bb)^{\frac{1-k}{k}} \pr[ H_{PC} (n,p\,;k) \in\h_2].  \eean

As $\pr[ H (n,p\,;k) \in\h]= \pr[ H (n,p\,;k) \in\h_1] + \pr[ H
(n,p\,;k) \in\h_2]$, we finally have
 \bean
\pr[ H (n,p\,;k) \in\h]& \leq & (1+o(1))\bb(  \bb(c_{_1} (k-1) \pr
[ H_{PC}(n,p\,;k)\in \h]\bb)^{1/k}+e^{-n
}\bb) \\
& & + (c_{_1}+o(1)) \bb(c_{_1} (k-1) \pr [ H_{PC}(n,p\,;k)\in
\h]\bb)^{\frac{1-k}{k}} \pr[ H_{PC} (n,p\,;k) \in\h]\\
 &\leq& c_{_2}\bb( \pr
[ H_{PC}(n,p\,;k)\in \h]^{1/k}+e^{-n }\bb).   \eean \qed

\bn

Theorem \ref{pc} may be generalized in the case that there are
some small number of forbidden edges. For example, the random
$k$-SAT formula $F(n,p\,;k)$ may be regarded  as $H_{PC}
(2n,p\,;k)$ on $V=\{ x_{_1}, \bar{x}_{_1}, ..., x_{_n},
\bar{x}_{_n}\}$ without  edges that contain both of a variable and
its negation. Suppose there is a set $B$ of forbidden edges with
$|B|=\frac{\gb}{n} {n \choose k}$ for $\gb=O(1)$. Each edge not in
$B$ is in the random $k$-uniform hypergraph $H^{(B)} (n, p\,;k)$
with probability $p$ independently of all other edges.

\begin{thm} \label{pc2} Suppose $k\geq 2$, $ p= \Theta (n^{1-k})$ and $B$ is a
set of $\frac{\gb}{n} {n \choose k}$ with   $\gb=O(1)$. Then, for
any collection ${\cal H}$ of simple $k$-uniform hypergraphs
without  edges in $B$,
 $$ c_{_1}(\gb)  \pr [
H_{PC} (n,p\,;k) \in \h] \leq  \pr [ H^{(B)} (n,p\,;k) \in \h]
 \leq
 c_{_2} (\gb) \bb( \pr [ H_{PC} (n,p\,;k) \in \h]^{{\frac{1}{k}}}
+ e^{-\tiny{\frac{p}{10}}{n\choose k}}\bb),
 $$
where
$$c_{_1} (\gb) = c_{_1} e^{\frac{p(\gb+o(1))}{n}{n
\choose k}}, ~~c_{_2} (\gb)= c_{_2} e^{\frac{p(\gb+o(1))}{n}{n
\choose k}}.
$$
\end{thm}

\mn \pf The result follows since
$$ \frac{\pr [ H^{(B)} (n,p\,;k) =H]}{\pr [ H (n,p\,;k) =H]} = (1-p)^{-|B|}
=e^{\frac{p(\gb+o(1))}{n}{n \choose k}}
$$
for $H \in \h$ implies that
$$
\pr [ H^{(B)} (n,p\,;k)\in \h ]= e^{\frac{p(\gb+o(1))}{n}{n
\choose k}}\pr [ H (n,p\,;k)\in \h ]. $$ \qed

\bn

\bn

 For the random $k$-SAT
formulas, suppose  $\gb_{_k}$ satisfies
$$ \frac{\gb_{_k}}{2n}{2n \choose k} = {2n \choose k} - 2^k {n
\choose k}, ~~~{\rm or~equivalently,} ~~~~ \frac{\gb_{_k}}{2n}=1-
\frac{2^k {n \choose k}}{ {2n \choose k} } . $$ Then
$$\frac{2^k {n \choose k}}{ {2n \choose k} }
=e^{-\frac{1}{n} {k \choose 2}+ \frac{1}{2n}{k \choose 2}
+O(\frac{1}{n^2})}= 1-\frac{1}{2n} {k \choose 2}
+O\bb(\frac{1}{n^2}\bb)$$ implies that
$$
 \gb_{_k}= {k \choose 2}  + O\bb(\frac{1}{n}\bb). $$

\begin{cor} If $k\geq 2$ and $p= \Theta (n^{1-k})$, then, for any
collection  $\f$ of $k$-SAT formulas,
$$ c_{_{1}}^* \pr [
F_{PC} (n,p\,;k) \in \f] \leq   \pr [ F (n,p\,;k) \in \f]
 \leq
c_{_{2} }^* \bb( \pr [ F_{PC} (n,p\,;k) \in \f]^{{\frac{1}{k}}} +
e^{-n}\bb),
 $$
where
$$c_{_{1}}^* = k^{1/2}e^{\frac{p}{n}{k \choose 2}{2n\choose k}
+ \frac{p^2}{2}{2n\choose k}}+o(1), ~~c_{_{2}}^* =
e^{\frac{p(1-1/k)}{2n}{k \choose 2}{2n\choose
k}}\bb(\frac{k}{k-1}\bb) \bb((k-1)c_{_{1}}^* \bb)^{1/k}+o(1).
$$
\end{cor}
\qed

\section{The Poisson $\gl$-Cell and the Cut-Off Line Algorithm}
\label{grpm}

\newcommand{\vp}{\xi}
\mn To generate a uniform random perfect matching of $N_\gl$
clones, we may keep matching $k$ unmatched clones uniformly at
random (cf. Example \ref{ex1}).  Another way is to choose the
first clone as we like and match it to $k-1$ clones selected
uniformly at random among all other unmatched clones (cf. Example
\ref{ex2}). As there are many ways to choose the first clone,  we
may take  a way that makes the given problem easier to analyze.
Formally,  a sequence $\mathcal{S}=(S_i)$  of choice functions
determines how to choose the first clone at each step, that is,
$S_i$ tells which unmatched clone is to be  the first clone to
form  the $i^{\rm th}$ edge in the random perfect matching. A
choice function may be deterministic or random. If less than $k$
colones remain unmatched, the edge consisting of those clones will
be added. The clone chosen by $S_i$ is called the $i^{\rm th}$
chosen clone, or simply a chosen clone.

We  also  present  a more specific way to select the $k-1$ random
clones to be matched to the chosen clone. The way presented here
will be useful to solve problems mentioned in the introduction.
First,
 independently assign to each clone a uniform random real number
between $0$ and $\gl=p{n-1 \choose k-1}$. For the sake of
convenience,  a clone is called the largest, the smallest, etc. if
so is the number assigned to it. In addition, for $ 0\leq \th \leq
1$,  a clone is called {\em $\th\gl$-large} (resp. {\em
$\th\gl$-small}) if its assigned number is larger than or equal to
(resp. smaller than) $\th\gl$. To visualize the labelled clones
with assigned numbers, one may consider $n$ horizontal line
segments from $(0,j)$ to $(\gl, j)$, $j=0, ..., n-1$ in the
two-dimensional plane $\mathbb{R}^2$. The $v_{_j}$-clone with
assigned number $x$ can be regarded as the point $(x,j)$ in the
corresponding line segment. Then, each line segment with the
points corresponding to clones with assigned numbers is an
independent Poisson arrival process with density $1$, up to time
$\gl$. The set of these Poisson arrival processes is called a {\em
Poisson $(\gl,n)$-cell}, or simply  a {\em $\gl$-cell}.

We will consider  sequences of choice functions that choose an
unmatched clone without changing the joint distribution of the
numbers assigned to all other unmatched clones. Such a  choice
function is  called {\em oblivious}. A sequence of oblivious
choice functions are also called {\em oblivious}. The choice
function that chooses the largest unmatched clone is not
oblivious, as the numbers assigned to the other clones must be
smaller than the largest assigned number. For an instance of an
oblivious choice function, one may consider
 the choice function that chooses a $v$-clone for a
vertex $v$ with fewer than $3$ unmatched clones. For a more
general example, let a vertex $v$ and its clones be  {\em
$t$-light} if there are fewer than $t$ unmatched $v$-clones.

\begin{ex}\label{ex3} Suppose there is  an order of all clones that is independent
of the assigned numbers. The sequence of the choice functions that
choose the first $t$-light clone is oblivious.
\end{ex}

   A cut-off line algorithm is determined by a  sequence
   of oblivious choice functions. Once a clone is obliviously chosen,
the largest $k-1$ clones among all
 unmatched clones are  to be matched  to the chosen clone.
 This may be further implemented by
moving the cut-off line to the left until $k-1$ vertices are
found: Initially, the cut-off line of the $\gl$-cell is the
vertical line in $\mathbb{R}^2 $ containing the point $(\gl, 0)$.
The initial cut-off value, or cut-off number, is $\gl$.
 At the first step, once the chosen clone is given,
 move the cut-off line to the left until exactly $k-1$
unmatched clones, excluding the chosen clone,  are on or in the
right side of the line. These $k-1$ clones together with the
chosen clone form the first edge in the random perfect matching.
The new cut-off value $\gL_1$ is to be the assigned number to the
$(k-1)^{\rm th}$ largest clone. Here we assumed that no two
distinct clones are assigned the same number as  the probability
of such an event  is $0$. The new cut-off line is, of course, the
vertical line containing $(\gL_1, 0)$. Repeating this procedure,
one may obtain the $i^{\rm th}$ cut-off value $\gL_i$ and the
corresponding cut-off line.

Notice that, after the $i^{\rm th}$ step ends with the cut-off
value $\gL_{i}$, all numbers assigned to unmatched clones are
i.i.d uniform random numbers between $0$ to $\gL_i$, as the choice
functions are oblivious. Let  $U_i$ be the number of unmatched
clones after step $i$. That is, $U_i=N_\gl -ik$.  Since the
$(i+1)^{\rm th}$ choice function tells how to choose  the first
clone to form the $(i+1)^{\rm th}$ edge without changing the
distribution of the assigned numbers, the distribution of
$\gL_{i+1}$ is the distribution of the $(k-1)^{\rm th} $ largest
number among $U_i -1$ independent uniform random numbers between
$0$ and $\gL_i$. Let $1-T_{j}$ be the random variable representing
the largest number among $j$ independent uniform random numbers
between $0$ and $1$. Or equivalently in distribution sense,
$T_{j}$ is the random variable representing the smallest number
among the random numbers.   Then the largest number among the
$U_i-1$ random numbers  has the same distribution as $\gll_i
(1-T_{U_i-1})$. Repeating this $k-1$ times, we have
$$ \gL_{i+1} = \gL_{i} (1-T_{U_i-1})(1-T_{U_i-2})\cdots
(1-T_{U_i-k+1}),$$ and hence \begin{align*}  \gL_{i+1} &= \gL_{i}
(1-T_{U_i-1})
\cdots (1-T_{U_i-k+1}) \\ &= \gL_{i-1}(1-T_{U_{i-1}-1})
 \cdots (1-T_{U_{i-1}-k+1})\cdot
(1-T_{U_i-1})
 \cdots (1-T_{U_i-k+1})
\\
&= \gl \prod_{j=N_\gl -1 \atop k \nmid N_\gl - j }^{N_\gl -
(i+1)k+1} \bb(1-T_{j} \bb) . \end{align*}

It is crucial   to observe  that, once $N_\gl$ is given, all $T_i$
are mutually independent random variables. This makes the random
variable $\gL_{i}$ highly concentrated near its mean, which
enables us to develop theories as if $\gL_{i}$ were a constant.
The cut-off value  $\gL_i$ will provide enough information to
resolve some otherwise difficult problems.

In the next section, we will prove the following  slightly general
lemma regarding the concentration of $\prod (1-T_j)$.

\begin{lem}
\label{cut-off}  For positive integer $k$, let $T_j$'s be mutually
independent, $j= N, N-1, ..., N-lk$ with $ N -lk \gg 1 $, and let
$R$ be a non-empty subset of $\{0, 1, ..., k-1\}$ with $|R|=r$.
Then, denoting $\gt_{_i} = (1-ik/N)^{1/k}$, we have, for $\eps\leq
0.1 $,
$$ \pr\bb[\,  \, \,   \max_{i:1\leq i \leq l} \bb| \prod_{j=N \atop
N-j \in_{k} R }^{\theta_{_i}^{k} N} \bb(1-T_{j} \bb)
-\theta_{_i}^{r} \bb|\geq \eps  \bb] \leq 10 e^{-\frac{1+o(1)}{7}
\min\{\eps \th_{_l}^k N, \, \frac{\eps^2 k \,\th_{_l}^k
N}{2(1-\th_{_l})}
 \}}.
$$
In particular, if $\th_{_l} = \Omega (1)$, then
$$ \pr\bb[\,  \, \,   \max_{i:1\leq i \leq l} \bb| \prod_{j=N \atop
N-j \in_{k} R }^{\theta_{_i}^{k} N} \bb(1-T_{j} \bb)
-\theta_{_i}^{r} \bb|\geq \eps  \bb] \leq 10 e^{-\Omega (
\min\{\eps  N, \, \frac{\eps^2  N}{2(1-\th_{_l})}
 \})}.
$$
\end{lem}

\old{For $\th$ in the range $0\leq \th \leq 1$, let $\gL(\th)$ be
the cut-off value when $(1-\th^{k}) \gl n $ or more clones are
matched for the first time. Conversely, let $N(\th)$ be the number
of matched clones until the cut-off line reaches  $\th^{k-1}
\gl$.}

 The cut-off line lemma follows from Lemma \ref{cut-off}:
 For $\th$ in the range $0\leq \th \leq 1$, let $\gL(\th)$ be the
cut-off value when $(1-\th^{\kk}) \gl n $ or more clones are
matched for the first time. Conversely, let $N(\th)$ be the number
of matched clones until the cut-off line reaches  $\th \gl$.

\begin{lem} (Cut-off Line Lemma)
\label{cf2} Let $k\geq 2$ and $\gl >0$ be fixed.
 Then, for $\tho <1 $ uniformly bounded below from $0$
 and $0<\dd\leq n$,
$$ \pr \bb[  \max_{\th:\tho \leq \th \leq 1 } |\gL (\th) -\th \gl|  \geq
\sfrac{\dd}{n}
 \bb]
  \leq  2e^{-\Omega(\min\{\dd, \frac{\dd^2}{(1-\tho)n}\})},$$
and $$ \pr \bb[  \max_{\th:\tho \leq \th \leq 1 } |N(\th)
-(1-\th^\kk) \gl n|
 \geq \dd  \bb]
  \leq  2e^{-\Omega(\min\{\dd, \frac{\dd^2}{(1-\tho)n}\})}.
    $$
\end{lem}

\newcommand{\xit}{\xi_{_\th}}
\pf  Suppose $N_\gl=\gl n + h$  is given. As $N_\gl$ is a Poisson
$\gl n$ random variable,
$$ \pr \bbr  | N_\gl - \gl n | \geq c  \min\bb\{ n, \frac{ \dd}{1-\tho}\bb\}\bb]
\leq 2 e^{-\Omega( \min\{n, \frac{\dd^2}{(1-\tho)^2 n }\})} \leq 2
e^{-\Omega( \min\{\dd, \frac{\dd^2}{(1-\tho) n }\})},
$$
where $c$ is a (small) constant to be specified later. Hence it is
enough to consider $h$ in the range $|h|\leq c \min\{n, \frac{
  \dd}{1-\tho}\}$.

For $\tho\leq \th\leq 1$, let $\xit $ be the solution of the
equation $ (1-\th^\kk ) \gl n = (1-\xit^\kk) N_\gl$. Then
$$ \xit =
\bb(1-\frac{(1-\th^\kk) \gl n}{N_\gl}\bb)^{\frac{k-1}{k}} = \bb(
\th^\kk + \frac{(1-\th^\kk) h}{\gln +h}\bb)^{\frac{k-1}{k}}= \th +
O\bb(\frac{(1-\th) |h|}{\gl n +h }\bb)$$ implies that $\xit$ is
uniformly bounded from below by $0$ (for small enough $c$). Lemma
\ref{cut-off} gives
$$  \pr \bb[ \max_{\th:\tho \leq \th \leq 1 }| \gL (\th) -
\xit \gl | \geq \sfrac{\dd}{2n}  \bb| N_\gl=\gl n + h \bb] \leq 2
e^{-\Omega(\min\{\dd, \frac{\dd^2}{(1-\xi_{\tho})n}\})}.
$$
\old{ 
$$ \xit = \th \bb(1+ \frac{(1-\th^k) h }{\th^k N_\gl}\bb)^{1/k}
 =  \th +
 O\bb(\frac{ (1-\th)|h |}{ n}\bb), 
~~ ~~  \xit^{k-1}  = \th^{k-1}  + O\bb(\frac{(1-\th)|h|}{ n}\bb)$$
and  $$ 1-\xi_{\tho}  = 1-\tho +O\bb(\frac{(1-\tho)|h|}{
n}\bb)=1-\tho+ O(\dd/n).$$} 
 Taking small enough $c$, we also have $| \xit -\th
 |\leq \frac{\dd}{2\gl n}$ and
 \bean
 |\gL (\th) - \th  \gl |
 &\leq &  |\gL (\th) - \xit  \gl | + \gl | \xit -\th
 |\\
 &\leq&  |\gL (\th) - \xit  \gl | +\frac{\dd}{2n}.
 \eean
Therefore, if $|\gL (\th) - \th  \gl |\geq \frac{\dd}{n}$ then
$|\gL (\th) - \xit  \gl |\geq \frac{\dd}{2n}$ and hence the
probability that such $\th$ exists in the range $\tho\leq \th\leq
1$ is at most
$$ 2 e^{-\Omega(\min\{\dd, \frac{\dd^2}{(1-\xi_{\tho})n}\})} \leq 2
e^{-\Omega(\min\{\dd , \frac{\dd^2}{(1-\tho )n+\dd/2\gl }\})} \leq
2 e^{-\Omega(\min\{\dd , \frac{\dd^2}{(1-\tho )n }\})}.
$$

For the second inequality, it is enough to observe that
 $ |N(\th) - (1-\th^\kk)
\gl n  | \geq \dd $   implies that  $\gL(\th + \Omega(\dd/n) )
\leq  \th \gl$ or $\gL(\th - \Omega(\dd/n) )\geq  \th \gl$.

\qed

For the Poisson $\gl$-cell conditioned on $N_\gl=N$, a similar
lemma may be obtained.

\begin{lem} (Cut-off Line Lemma for $N$ clones)
\label{cf3} Let $k\geq 2$, $\gl >0$ be fixed.
 Then, for the Poisson $\gl$-cell conditioned  on $N_\gl=N$, and for
 $\tho <1 $ uniformly bounded below from $0$
 and $0<\dd\leq N$,
$$ \pr \bb[  \max_{\th:\tho \leq \th \leq 1 } |\gL (\th) -\th \gl|  \geq
\sfrac{\dd}{N}
 \bb]
  \leq  2e^{-\Omega(\min\{\dd, \frac{\dd^2}{(1-\tho)N}\})},$$
and $$ \pr \bb[  \max_{\th:\tho \leq \th \leq 1 } |N(\th)
-(1-\th^\kk) N|
 \geq \dd  \bb]
  \leq  2e^{-\Omega(\min\{\dd, \frac{\dd^2}{(1-\tho)N}\})}.
    $$
\end{lem}

\section{Large Deviation Inequalities }
\label{sldi}

In this section, a generalized  Chernoff bound and an inequality
for random processes are to be shown.  Let $X_1, ..., X_m$ be a
sequence of random variables such that the distribution of $X_i$
is determined if all the values of $X_1, ..., X_{i-1}$ are known.
For example, $X_i = \gls (\th_i)$ with $1\geq \th_1 \geq \cdots
\geq \th_m \geq 0$ in a Poisson $\gl$-cell. If the upper and/or
lower bounds are known for the conditional means $E[X_i| X_1, ...
, X_{i-1}]$ and for the conditional second and third moments, then
Chernoff type large deviation inequalities may be obtained not
only for $\sum_{j=1}^m X_j$ but $\min_{1\leq i\leq m} \sum_{j=1}^i
X_j$ and/or $\max_{1\leq i\leq m} \sum_{j=1}^i X_j$. Large
deviation inequalities for such minimums or maximums are
especially useful in various situations. Lemma \ref{cut-off} can
be shown using such inequalities too.

\begin{lem}\label{basic}
Let $ X_1, ..., X_m $ be a sequence of   random variables. Suppose
  \beq{mean} E[X_i| X_{1}, ..., X_{i-1} ] \leq
\mu_i , \eeq and there are positive constants $a_{_i}$, $b_{_i}$,
and $\vp_{0}$ so that
 \beq{second}  E[ (X_i-\mu_i)^2|X_{1}, ..., X_{i-1} ]
\leq a_i, \eeq
 and
  \beq{thirdco}
E[ (X_i - \mu_i)^3 e^{\vp (X_i-\mu_i) }| X_{1}, ..., X_{i-1} ]
\leq b_i ~~~\mbox{for all $0\leq \vp \leq \vp_{_0} $}.\eeq
 Then for any $\ga$ with $0<\ga
 \leq \vp_{_0} (\sum_{i=1}^m a_i)^{1/2}$,
$$ \pr \bb[ \sum_{i=1}^m X_i  \geq  \sum_{i=1}^{m} \mu_i
+\ga \bb(\sum_{i=1}^m a_i\bb)^{1/2} \bb] \leq \exp\bb(
-\frac{\ga^2}{2}\bb(1- \frac{\ga \sum_{i=1}^m b_i}{3 (\sum_{i=1}^m
a_i)^{3/2} }\bb)\bb). $$ Similarly, \beq{mean2} E[X_i| X_{1}, ...,
X_{i-1} ] \geq \mu_i \eeq together with \raf{second} and
\beq{thirdcos} E[ (X_i - \mu_i)^3 e^{\vp (X_i-\mu_i) }| X_{1},
..., X_{i-1} ] \geq  b_i ~~~\mbox{for all $-\vp_{_0}\leq \vp< 0
$}\eeq
  implies that
$$ \pr \bb[ \sum_{i=1}^m X_i  \leq \sum_{i=1}^m \mu_i
-\ga \bb(\sum_{i=1}^m a_i\bb)^{1/2} \bb] \leq \exp\bb(
-\frac{\ga^2}{2}\bb(1-\frac{\ga \sum_{i=1}^m b_i}{3 (\sum_{i=1}^m
a_i)^{3/2} }\bb)\bb). $$
\end{lem}

\mn \pf  We first   show  that
$$ E[ e^{{\vp \sum_{j=1}^{i} (X_j-\mu_j)} }] \leq e^{\vp^2\sum_{j=1}^i
 a_{_j}/2 +  |\vp|^3 \sum_{j=1}^i b_{_j} /6}
   ~~~\mbox{for $i=0,...,n$}, $$
using induction. As
 \bean E[e^{{\vp \sum_{j=1}^{i}
(X_j-\mu_j)} }] &=& E\bb[E[e^{{\vp
\sum_{j=1}^{i} (X_j-\mu_j)} }|X_{1},..., X_{i-1}]\bb] \\
& =& E\bb[e^{{\vp \sum_{j=1}^{i-1} (X_j-\mu_j)} }E[e^{{\vp
(X_i-\mu_i)} }|X_{1},..., X_{i-1}]\bb], \eean it is enough to show
$$E[e^{{\vp (X_i-\mu_i)} }|X_{1},..., X_{i-1}] \leq e^{ \vp^2
a_{_i}/2 +  \vp^3 b_{_i} /6}. $$
 For $0< \vp\leq \vp_{_0} $,     Taylor theorem
gives
 \bean E[ e^{\vp (X_i -\mu_i)}|X_{1},..., X_{i-1}] \! \! \!  \! &= &  \! \! \! \! 1+
\frac{\vp^2 E[(X_i -\mu_i)^2|X_{1},..., X_{i-1}]}{2} +\frac{\vp^3
E[(X_i -\mu_i)^3
 e^{\vp^* (X_i
-\mu_i)}|X_{1},..., X_{i-1}]}{6}  \\
\! \! \!  \! &\leq & \! \! \! \!  1+ \frac{\vp^2 a_{_i}}{2} +
 \frac{\vp^3 b_{_i}}{6}\leq e^{
\vp^2 a_{_i}/2 +  \vp^3 b_{_i} /6}
   ,
 \eean for some $\vp^*$ between $0$ and $\vp$.

Let $\vp= \ga (\sum_{i=1}^m a_i )^{-1/2}
    \leq \vp_{_0}$. Then
$$E[ e^{\vp \sum_{j=1}^m (X_j -\mu_j) }]
\leq \exp\bb(\frac{\ga^2}{2}
    +\frac{\ga^3 \sum_{i=1}^m b_i}{6
(\sum_{i=1}^m a_i)^{3/2} }\bb ),$$ and $$
 \pr \bb[ \sum_{j=1}^m X_j - \sum_{j=1}^m \mu_j  \geq \ga \bb(\sum_{i=1}^m a_i \bb)^{1/2}\bb]
 \leq E[
   e^{\vp (\sum_{j=1}^m (X_j -\mu_j)-\ga (\sum_{i=1}^m a_i )^{1/2})}]
 \leq \exp\bb(
    -\frac{ \ga^2}{2} +\frac{\ga^3 \sum_{i=1}^m b_i}{6
(\sum_{i=1}^m a_i)^{3/2} } \bb ).
 $$
 Similarly, \raf{mean2} together with \raf{second},
\raf{thirdco} and  $\vp= -\ga (\sum_{i=1}^m a_i
 )^{-1/2}$ together with \raf{thirdcos} gives
$$
 \pr \bb[ \sum_{j=1}^m X_j - \sum_{j=1}^m \mu_j  \leq -\ga
 \bb(\sum_{i=1}^m a_i \bb)^{\frac{1}{2}}\bb]
 \leq E[
   e^{\vp (\sum_{j=1}^m (X_j -\mu_j)+\ga (\sum_{i=1}^m a_i )^{\frac{1}{2}})}]
 \leq \exp\bb(
    -\frac{ \ga^2}{2} +\frac{\ga^3 \sum_{i=1}^m b_i}{6
(\sum_{i=1}^m a_i)^{3/2} } \bb ).
$$

 \qed

\bn

 As it  is sometimes tedious   to point out the
value of $\ga$ and to check the required bounds  for it,  the
following forms of inequalities are often more convenient.

\begin{cor} \label{uni}(Generalized Chernoff bound)  If   $ \gd\vp_{_0} \sum b_i \leq  \sum a_{_i} $
for some $0<\gd\leq 1$, then \raf{mean}-\raf{thirdco} imply
$$ \pr \bb[ \sum_{i=1}^m X_i   \geq  \sum_{i=1}^m \mu_i+ R  \bb] \leq
e^{-\frac{1}{3}\min \{\gd  \vp_{_0}R , \, \,  R^2/ \sum_{i=1}^m
a_i \}},
$$ for all $R>0$. Similarly, If   $ -\gd\vp_{_0} \sum b_i \leq  \sum
a_{_i} $ for some $0<\gd\leq 1$, then \raf{second}, \raf{mean2}
and  \raf{thirdcos}  yield
$$ \pr \bb[ \sum_{i=1}^m X_i   \leq \sum_{i=1}^m \mu_i- R \bb]
\leq e^{-\frac{1}{3} \min \{  \gd \vp_{_0}R, \, \,  R^2/
\sum_{i=1}^m a_i
                        \}} $$
for all $R>0$.
\end{cor}

\noindent\pf For $R \leq \gd  \vp_{_0} \sum_{i=1}^m a_i$,  Lemma
\ref{basic}  with $\ga= R(\sum_{i=1}^m a_{_i})^{-1/2}$ gives
$$ \pr \bb[ \sum_{i=1}^m X_i  \geq   \sum_{i=1}^m \mu_i+ R  \bb] \leq
\exp\bb( - \frac{R^2}{ 3\sum_{i=1}^m a_i } \bb).
$$
 If $R \geq \gd \vp_{_0}\sum_{i=1}^m a_i$, one may replace
$a_{_i}$ by $a_{_i}^* \geq a_{_i}$ satisfying $\sum a_{_i}^* =
R/(\gd \vp_{0} ) $ and obtain
$$ \pr \bb[ \sum_{i=1}^m X_i - \sum_{i=1}^m \mu_i  \geq  R  \bb]
=\pr \bb[ \sum_{i=1}^m X_i - \sum_{i=1}^m \mu_i  \geq  (\gd
\vp_{_0} R)^{1/2} \bb(\frac{R}{\gd\vp_{_0}} \bb)^{1/2} \bb] =
\exp( -\gd\vp_{_0}  R /3 ),
$$
as $$ \frac{(\gd \vp_{_0} R)^{1/2} \sum b_i}{3(\sum a^*_i)^{3/2}}=
     \frac{(\gd \vp_{_0} R)^{1/2} \sum b_i}{3(R/(\gd\vp_{_0}))^{3/2}}
     = \frac{\gd^2 \vp_{_0}^2 \sum b_i}{3R} \leq \frac{\gd \vp_{_0} \sum
      a_i}{3R} \leq \frac{1}{3}.
$$

\qed

We now try to obtain inequalities for the maximum and the minimum
of the random process.  Let $X_\th$, $\th\geq 0$, be random
variables, which are possibly set-valued. Here $\th$ may be real
numbers as well as non-negative integers. For $\th\geq 0$, suppose
$\Gamma(\th)$ is a random variable depending on
$\{X_{\th'}\}_{\th' \leq \th}$ and $\th$, and
$$\psi=\psi (\{X_{\th'}\}_{\th'\leq \tho}; \thz,\tho), ~~{\rm and} ~~
\psi_\th= \psi_\th ( \{X_{\th'}\}_{\th'\leq \tho};\thz,  \th,
\tho). $$ The random variables $\psi $ and $\psi_\th$ are to be
used to bound $\Gamma(\th)$.

\begin{ex}\label{ex4} Let $X_1,X_2, ...$ be i.i.d Bernoulli random variables with
mean $p$ and $S_i = \sum_{j=1}^i X_j$. Set $\Gamma(i) = |S_i -
ip|$ and
  $$\psi = \Gamma(n)~~{\rm and}~~ \psi_i =| S_n -S_{i} - (n-i) p|.
  $$
Then, since
$$ S_i -ip  = S_n - np  - ( S_n -S_{i} - (n-i)
p), $$  we have
$$ \Gamma(i) \leq \psi + \psi_i.
$$
\end{ex}

\begin{ex}\label{ex5} ~Consider  the $\gl$-cell on $n$ vertices defined in the previous
section.  Let $v_\th $ be the vertex that has its largest clone at
$(1-\th)\gl$.
 If such a vertex does not exist, $v_\th$ is defined to be $\al$,
 assuming $\al \not\in V$.
 As there is no
possibility that two distinct clones are assigned the same number,
$v_\th$ is well-defined. Let $X_\th=v_\th$ and $\VV(\th)$ be the
set of vertices that contain no clone larger than or equal to
$(1-\th)\gl$. That is, $\VV(\th) = V\sm \{ v_{\th'}: 0\leq \th'
\leq \th\}$. Clearly, $ E[ |\VV(\th)|] = e^{-\th\gl} n. $
Observing that, for $\thz\leq \th\leq \tho$,
$$ e^{-(\tho-\th)\gl}\bb| |\VV(\th)| - e^{-\th\gl} n\bb| \leq
\bb| |\VV(\tho)|-  e^{-\tho\gl }n \bb|+  \bb| |\VV(\tho)| -
e^{-(\tho-\th)\gl} |\VV(\th)|\bb|, $$ we may set $\Gamma(\th) = |
|\VV(\th)| - e^{-\th\gl} n|$, $$\psi = e^{(\tho-\thz
)\gl}\Gamma(\tho), \and \psi_\th = e^{(\tho-\thz )\gl}\bb|
|\VV(\tho)| - e^{-(\tho-\th)\gl} |\VV(\th)|\bb|. $$ \end{ex}

\mn

 We bound the probability $\max_{\thz\leq \th \leq \tho}
\Gamma(\th) \geq R$ and $\min_{\thz\leq \th \leq \tho} \Gamma(\th)
\leq R$ under some conditions.

\begin{lem} \label{gcb2}  Let $0\leq \thz< \tho $, $R=R_1+R_2$, $R_1,R_2>0$
and $\Phi_\th$ be events depending on $\{X_{\th'}\}_{\th'\leq
\th}$. If
$$\Gamma(\th) \leq  \psi + \psi_\th, ~~\forall~\thz\leq \th \leq \tho ,$$
then
\begin{align*} \pr \bb[ \max_{\thz\leq \th \leq \tho} \Gamma(\th) \geq R
\bb] &\leq \pr \bb[ \psi \geq R_1\bb ]+ \pr\bb[
\bigcup_{\th:\thz\leq \th \leq \tho} \overline{\Phi}_\th \bb]
\\ & + \max_{\th: \thz\leq \th \leq \tho} \max_{\{X_{\th'} \}_{\th' \leq \th}}
 1(\Phi_{\th}) \pr
\bb[ \psi_\th \geq R_2 \bb| \{X_{\th'} \}_{\th' \leq \th}\bb].
\end{align*}
Similarly, if
$$\Gamma(\th) \geq  \psi  +
\psi_\th,  ~~\forall~\thz\leq \th \leq \tho , $$
 then  \begin{align*}  \pr
\bb[ \min_{\thz\leq \th \leq \tho} \Gamma(\th) \leq - R \bb] &\leq
\pr \bb[ \psi \leq  - R_1\bb ]+ \pr\bb[
\bigcup_{\th:\thz\leq \th \leq \tho} \overline{\Phi}_\th \bb]  \\
& + \max_{\th: \thz\leq \th \leq \tho} \max_{\{X_{\th'} \}_{\th'
\leq \th}}
 1(\Phi_\th) \pr \bb[ \psi_\th
\leq  - R_2 \bb| \{X_{\th'} \}_{\th' \leq \th}\bb].
\end{align*}
\end{lem}


\mn {\bf Example \ref{ex4}} {\em  (continued) As
$$ \pr [ \psi \geq R_1] \leq e^{-\Omega (\min\{R_1,
\frac{R_1^2}{p(1-p)n}\})} $$ and $$ \pr [ \psi_i  \geq R_2| X_1,
..., X_i ] = \pr [ \psi_i  \geq R_2 ] \leq e^{-\Omega (\min\{R_2,
\frac{R_2^2}{p(1-p)(n-i)}\})},
$$
Lemma \ref{gcb2} for  $R_1=R_2=R/2$ and $\Phi_\th=\emptyset$ gives
$$ \pr [ \max_{i:0\leq i \leq n} |S_i -pi|  \geq R]
\leq  e^{-\Omega (\min\{R, \frac{R^2}{p(1-p)n}\})}. $$ }

\noindent {\bf Example \ref{ex5}}  {\em (continued) Since
$$ |\VV(\th)| =\sum_{v\in V} 1(\mbox{$v$ has no $(1-\th)\gl$-large
clone}) $$ is a sum of i.i.d Bernoulli random variables with mean
$e^{-\th\gl}$,
$$ \pr\bb[ \bb| |\VV(\th)|- e^{-\th\gl} n \bb| \geq R \bb]
\leq 2e^{-\Omega (\min\{R, \frac{R^2}{\th n}\})}, $$ especially $$
 \pr\bb[ \psi  \geq R/2 \bb]
\leq 2e^{-\Omega (\min\{R, \frac{R^2}{\tho n}\})}. $$ Once $\{
X_{\th'}\}_{\th'\leq \th}$ is given,   $\VV(\th)$ is determined
and
$$ \VV(\tho) = \sum_{v\in \VV(\th)} 1(\mbox{$v$ has no $(1-\tho)\gl$-large
clone}) $$ is a sum of i.i.d Bernoulli random variables with mean
$e^{-(\tho-\th)\gl}$. Thus,
  $$ \pr\bb[ \psi_\th  \geq R/2 \bb| \{
X_{\th'}\}_{\th'\leq \th}  \bb] \leq 2e^{-\Omega (\min\{R,
\frac{R^2}{(\tho-\th)|\VV(\th)|}\})} \leq 2e^{-\Omega (\min\{R,
\frac{R^2}{\tho n}\})},
$$
and Lemma \ref{gcb2} for $\thz=0$ and $\Phi_\th= \emptyset$ yields
$$ \pr \bb[ \max_{\th:0\leq \th \leq \tho} \bb ||\VV(\th)|-e^{-\th \gl }n
\bb|\geq R\bb]\leq 2e^{-\Omega (\min\{R, \frac{R^2}{\tho n}\})}.
$$ }

The proof of Lemma \ref{gcb2} follows.

\mn {\bf Proof of Lemma \ref{gcb2}}  Let $\tau$ be the first time
$\th$ in the range $\thz\leq \th \leq \tho$ when $\Gamma(\th) \geq
R$.
If no such $\th$ exists, $\tau=\infty$ and $\psi_\tau =-\infty$.
Observe that
 \bean
 \pr \bb[ \tau < \infty  \bb] &\leq &
 \pr \bb[ \tau < \infty, ~\psi \leq R_1, \bigcap_{\th:\thz\leq \th \leq \tho} \Phi_\th \bb] +
\pr [ \psi \geq R_1] + \pr\bb[ \bigcup_{\th:\thz\leq \th \leq
\tho} \overline{\Phi}_\th \bb]\\
&\leq& \pr \bb[  ~\psi_\tau \geq R_2, \bigcap_{\th:\thz\leq \th
\leq \tho} \Phi_\th \bb]+ \pr [ \psi \geq R_1] + \pr\bb[
\bigcup_{\th:\thz\leq \th \leq \tho} \overline{\Phi}_\th \bb].
\eean Considering the conditional probability on  $\{
X_\th\}_{\th\leq \tau}$, we have  \bean
  \pr \bb[ ~\psi_\tau \geq R_2, \bigcap_{\th:\thz\leq
\th \leq \tho} \Phi_\th \bb] &=&   E\bb[  \pr \bb[ \psi_\tau  \geq
R_2, \bigcap_{\th:\thz\leq \th \leq \tho} \Phi_\th  \bb| \{
X_\th\}_{\th \leq \tau}\bb] \bb] \\
&\leq & E\bb[  1(\Phi_\tau)  \pr \bb[ \psi_\tau \geq R_2 \bb| \{
X_\th\}_{\th \leq \tau}\bb] \bb]. \eean

As
$$1(\Phi_\tau)  \pr \bb[ \psi_\tau \geq R_2 \bb| \{
X_\th\}_{\th \leq \tau}\bb] \leq \max_{\th: \thz\leq \th \leq
\tho} \max_{\{X_{\th'} \}_{\th' \leq \th}}
 1 (\Phi_\th )  \pr
\bb[ \psi_\th \geq R_2 \bb| \{X_{\th'} \}_{\th' \leq \th}\bb], $$
the desired inequality follows.

 Applying the same argument  for $-\Gamma (\th)$, the second part also follows.

\qed

\mn

\bn
 {\bf Proof of Lemma \ref{cut-off} }  As
 $$\prod_{j=N  \atop N-j \in_k R}^{\gt_{_i}^k N} (1-T_{j})
   = \exp\bb(\sum_{j=N  \atop N-j \in_k R}^{\gt_{_i}^k N} \log
   (1-T_j)\bb), $$
   we  show a high concentration for
 $\log (1-T_j)$.
Since $ \pr [ \exists ~j, ~ T_j \geq 1/2] \leq
\sum_{j=N}^{\gt_{_l}^k N} 2^{-j} \leq 2^{-\gt_{_l}^k N +1},$ and
$$ -x -x^2 \leq \log (1-x) \leq -x  ~~\forall x : 0\leq x \leq 1/2, $$
with probability at least $1- 2^{-\th_{_l}^k N+1}$, we have
$$ -T_j - T_j^2 \leq  \log
   (1-T_j) \leq - T_j, ~~~\mbox{for all $j$}. $$
 Thus, it is enough to show that
 both of $E[S^*_i]$ and $E[T^*_i]$ are very close to
$\th_{_i}^r$, and
 $$ T^*_i  := \sum_{j=N  \atop N-j \in_k R}^{\gt_{_i}^k N}
  T_j,  \and  S^*_i:= \sum_{j=N  \atop N-j \in_k R}^{\gt_{_i}^k N}
S_j $$ are highly concentrated.  That is, we will show that
$$ \pr [ \max_i  |S^*_i -E[S^*_i]| \geq \eps ]
\leq 4 e^{-\frac{1+o(1)}{6} \min\{\frac{\eps^2 k \, \th_{l}^k
N}{2(1-\th_{l})}, \, \eps \th_{_l}^k N \}},$$ and
$$\pr [ \max_i
|T^*_i -E[T^*_i]| \geq \eps ]\leq 4 e^{-\frac{1+o(1)}{6}
\min\{\frac{\eps^2 k \, \th_{l}^k N}{2(1-\th_{_l})}, \, \eps
\th_{l}^k N \}},  $$ together with
$$ E[S^*_i] =-r \log \theta_{_i} + o\bb(
\bb(\frac{1-\gt_{_i}}{\theta_{_i}^k N}\bb)^{1/2} \bb), ~~{\rm and}
~~E[T^*_i]= -r \log \theta_{_i} + o\bb(
\bb(\frac{1-\gt_{_i}}{\theta_{_i}^k N}\bb)^{1/2} \bb). $$ The
$o\bb( \bb(\frac{1-\gt_{_i}}{\theta_{_i}^k N}\bb)^{1/2} \bb)$
terms do not matter,  since the desired inequality is trivial
unless $\eps = \Omega((\frac{1-\th_{_l}}{\th_{_l}^k N})^{1/2})$.
If $\eps = \Omega((\frac{1-\th_{_l}}{\th_{_l}^k N})^{1/2})$, then
the above concentration inequalities for $0.95 \eps$ give
$$ \pr \bb[ \max_i \bb|\sum_{j=N  \atop N-j \in_k R}^{\gt_{_i}^k N}
\log (1-T_j)   - r\log \thi \bb| \geq 0.95\eps \bb] \leq 8
e^{-\frac{1+o(1)}{7} \min\{\frac{\eps^2 k \, \th_{l}^k
N}{2(1-\th_{_l})}, \, \eps \th_{l}^k N  \}},
$$
which along with $\eps\leq 0.1$ yields
  \bean
\pr\bb[\,  \, \,   \max_{i:1\leq i \leq l} \bb| \prod_{j=N \atop
N-j \in_{k} R }^{\theta_{_i}^{k} N} \bb(1-T_{j} \bb)
-\theta_{_i}^{r} \bb|\geq \eps  \bb] &\leq&  8
e^{-\frac{1+o(1)}{7} \min\{\frac{\eps^2 k \, \th_{l}^k
N}{2(1-\th_{_l})}, \, \eps \th_{l}^k N  \}}+ 2^{-\gt_{_l}^k N +1}\\
&\leq & 10 e^{-\frac{1+o(1)}{7} \min\{\frac{\eps^2 k \, \th_{l}^k
N}{2(1-\th_{_l})}, \, \eps \th_{l}^k N  \}}.\eean

For the concentration inequalities without the maximum, it is
enough to check the hypotheses of the generalized Chernoff bound.
First, it is routine to check that
$$ E[ T_j^h] = \frac{ h!}{(j+1)(j+2)\cdots (j+h)} \leq
\frac{h!}{j^h},  $$  for positive integer $h$. Thus,   for $j\gg
1$,
   \beq{var}  Var[S_j] \leq    E[S_j^2]
\leq \frac{3}{j^2} =:a_{_j} , ~~{\rm and}~ \bb|\sum_{j=N  \atop
N-j \in_k R}^{\gt_{_l}^k N} \frac{3}{j^2} -
\frac{3r(1-\gt_{_l}^k)}{k \gt_{_l}^{k} N}\bb|\leq
\frac{4k(1-\gt_{_l}^k)}{\gt_{_l}^{2k}
N^2}=o\bb(\frac{1-\gt_{_l}}{\th_{_l}^{k} N}\bb). \eeq Furthermore,
for $0< \vp \leq \vp_{_0}:=0.1\gt_{_l}^k N$,
$$  E[S_j^3 e^{\vp(S_j - E[S_j])}]
\leq  E[S_j^3 e^{\vp S_j }] \leq    \frac{10}{j^3}=:b_{_j},$$ and,
for $- \xi_0 \leq  \vp < 0 $,
$$  E[S_j^3 e^{\vp(S_j - E[S_j])}]
\leq E[S_j^3 e^{\vp E[S_j]}] \leq
 \frac{10}{j^3} .$$
Since
   \beq{tm}  \sum_{j=N  \atop N-j \in_k R}^{\gt_{_l}^k N}
\frac{10}{j^3} \leq  \frac{6r(1-\gt_{_l}^{2k})}{k\gt_{_l}^{2k}
N^2}, \eeq we have that $ \vp_{_0} \sum  b_j \leq \sum_{j=1}^l
a_{_i}$ and hence, for $ \eps>0$, the generalized Chernoff bound
(Corollary \ref{uni}) gives
$$  \pr \bb[ | S^*_i  -E[ S_i^*]|\geq \eps \bb]
\leq 2e^{-\frac{1+o(1)}{3} \min\{\frac{\eps^2 k \, \th_{l}^k
N}{1-\th_{_l}}, \, \eps \th_{l}^k N \}}
$$ and
$$
\pr \bb[ | S^*_{i}-S^*_{l}  -E[S^*_{i}-S^*_{l} ]|\geq \eps \bb]
\leq 2e^{-\frac{1+o(1)}{3} \min\{\frac{\eps^2 k \, \th_{l}^k
N}{1-\th_{_l}}, \, \eps \th_{l}^k N \}}.
$$
As
$$ |S^*_i -E[S^*_i]| \leq  |S^*_l -E[S^*_l] | +|S^*_{i}-S^*_{l}  -E[S^*_{i}-S^*_{l}
]|, $$ and $S^*_{i}-S^*_{l}$ is independent of $S_1,..., S_i$, we
may apply Lemma \ref{gcb2} with $\eps_1=\eps_2=\eps/2$ to obtain
$$\pr [ \max_{i:1\leq i\leq l}  |S^*_i -E[S^*_i]| \geq
\eps ] \leq 4 e^{-\frac{1+o(1)}{6} \min\{\frac{\eps^2 k \,
\th_{l}^k N}{2(1-\th_{_l})}, \, \eps \th_{l}^k N \}}.
$$
Similarly,
$$\pr [ \max_{i:1\leq i\leq l}  |T^*_i -E[T^*_i]| \geq
\eps ] \leq 4e^{-\frac{1+o(1)}{6} \min\{\frac{\eps^2 k \,
\th_{l}^k N}{2(1-\th_{_l})}, \, \eps \th_{l}^k N \}}.
$$

For the expectations, since
$$ \sum_{j=N  \atop N-j \in_k R}^{\gt_{i}^k N}
 E[ S_j]= \sum_{j=N  \atop N-j \in_k R}^{\gt_{i}^k N}
\frac{1+O(1/j)}{j+1} \leq \frac{r}{k}  \int_{\theta_{i}^k
N+1}^{N}\frac{dx}{x}+ O\bb(\frac{1-\gt_{i}^k}{\theta_i^k N}\bb)
=-r\log \theta_{i} + o\bb( \bb( \frac{1-\gt_{i}^k}{\theta_{i}^k
N}\bb)^{1/2} \bb), $$ and
$$ \sum_{j=N  \atop N-j \in_k R}^{\gt_{i}^k N}
 E[ S_j] \geq \sum_{j=N  \atop N-j \in_k R}^{\gt_{i}^k N} E[ T_j]=
 \sum_{j=N  \atop N-j \in_k R}^{\gt_{i}^k N}
  \frac{1}{j+1} \geq \frac{r}{k}  \int_{\theta_i^k
N+2}^{N+1}\frac{dx}{x}= -r \log \theta_i + o\bb(
\bb(\frac{1-\gt_{i}^k}{\theta_i^k N}\bb)^{1/2} \bb),
$$
we have
$$ E[S_i^*] =-r \log \theta_i + o\bb(
\bb(\frac{1-\gt_{i}^k}{\theta_i^k N}\bb)^{1/2} \bb), ~~{\rm and}
~~E[T_i^*]= -r \log \theta_i + o\bb(
\bb(\frac{1-\gt_{i}^k}{\theta_i^k N}\bb)^{1/2} \bb)$$ as desired.

\old{ For $1\ll
 \ga\leq 0.1 (\gt_{_l}^k N )^{1/2} $, we may take each bound larger so that
 the bound in \raf{var}  and \raf{tm} become $\frac{3r}{k \gt_{_l}^{k} N}+O\bb(
\frac{1}{\gt_{_l}^{2k} N^2}\bb)$ and $ \frac{6}{\gt_{_l}^{2k}
N^2},$ respectively. Thus we have
 $$ \pr\bb[\,  \, \,  \bb| \prod_{j=N
\atop N-j \in_{k} R }^{\theta^{k} N} \bb(1-T_{j} \bb) -\theta^{r}
\bb|\geq \ga (\theta^k N)^{-1/2}\bb] \leq e^{-\Omega(\ga^2)}.
$$}
\qed

\section{Generalized Core-Processes and Main Lemma}
\label{ldi}

In this section, we introduce generalized cores and the main
lemma. The main lemma will be  crucial  in  the proofs of theorems
mentioned in the introduction.

We start with a few terminology.  A {\em generalized degree} is an
ordered pair $(d_1, d_2)$ of
 non-negative integers. The inequality between two generalized degrees
is determined by the inequality between the first coordinates and
the reverse inequality between the second coordinates. That is,
 $(d_1, d_2)\geq (d'_1, d'_2)$ if and only if
  $d_1 \geq d_1' $ and $ d_2 \leq  d'_2$.
A {\em property} for generalized degrees is simply a set of
generalized degrees. A property $P$ is {\em increasing}  if
generalize degrees larger than an element  in $P$ are also in $P$.
 When a property $P$ depends only on
the first coordinate of generalized degrees, it is  a property for
degrees. For the $t$-core problem, we will use $P_{t-{\rm core}}=
\{ (d_1,d_2): d_1\geq t\}$. To estimate the size of the largest
component, we will set $P_{\rm comp} = \{ (d_1,d_2): d_2  =0 \}$.

Given the Poisson $\gl$-cell on the set $V$ of $n$ vertices and
$\th$ in the range $0\leq \th\leq 1$,  let $d_v (\th)$ be the
number of $v$-clones smaller than $\th\gl$. Similarly, $\bd_v
(\th) $ is the number of $v$-clones larger  than or equal to
$\th\gl$. Then, $D_v (\th) := (d_v (\th), \bd_v (\th))$ are i.i.d
random variables. In particular, for any property $P$, the events
$D_v (\th) \in P$ are independent and occur with the same
probability, say $p(\th, \gl; P)$, or simply $p(\th)$.

For an increasing property $P$, the $P$-process is defined as
follows. Construct the Poisson $\gl$-cell as described in Section
\ref{grpm}, where $\gl= p{n-1 \choose k-1}$. The vertex set $V=\{
v_{_0}, ..., v_{_{n-1}} \}$  will be regarded as an ordered set so
that the $i^{\rm th}$ vertex is $v_{_{i-1}}$. The $P$-process, or
generalized core-process generated by $P$, is a generalization of
Example \ref{ex2} for which choice functions choose $t$-light
clones.

\mn {\bf The  $P$-process}: Initially, the cut-off value $\gL=
\gl$. Activate all vertices $v$  with $D_v (1)\not \in P$.
 All clones of the activated vertices are
activated too. Put activated  clones in a stack in an arbitrary
order. However, this does not mean that the clones are removed
from the $\gl$-cell.

\mn (a)   If the stack is empty, go to (b). If the stack  is
nonempty, choose the first clone in the stack and move the cut-off
line to the left  until the largest $k-1$ unmatched clones,
excluding the chosen clone, are found. (So, the cut-off value
$\gL$ keeps decreasing.) Then, match the $k-1$ clones to the
chosen clone. Remove all matched clones from the stack  and
repeat. A vertex $v$ that has not been activated is to be
activated as soon as $D_v (\gL/\gl)\not\in P$. This can  be done
even before all $k-1$ clones are found.  Its unmatched clones are
to be activated too and put into the stack immediately. Clones
found while moving the cut-off line are also in the stack until
they are matched.

\mn (b) Activate the first vertex in $V$ that has not been
activated. Its clones are activated too. Put those clones into the
stack. Then, go to (a).

\mn Clones in the stack are called active.
The steps carried by the instruction described in (b) are called
{\em free steps} as we are free to choose any clone.

\medskip
 When the cut-off line is at $\th\gl$,  all $\th\gl$-large
clones are matched or will be matched at the end of the step and
all vertices $v$ with $D_v (\th) \not\in P$ have been activated.
All other vertices can have been activated only by free steps.
 Let $V(\th)=V_{P} (\th)$ be the set  of vertices $v$ with $D_v
(\th) \in P$, and let  $M(\th)=M_{P}(\th)$ be the number of $\th
\gl$-large clones plus the number of $\th \gl$-small clones of
vertices $v$ not in $V(\th)$. That is,
  \bq{newm} {M(\th) = \sum_{v\in V} \bd_v (\th) + d_v (\th ) 1
(v\not \in V (\th)) = \sum_{v\in V} \bd_v (\th) + d_v (\th) 1 (D_v
(\th ) \not\in P). }

Recalling that $N(\th)$ is the number of matched clones until the
cut-off line reaches $\th \gl$, the number $A(\th)$  of active
clones (when the cut-off value $\gL$ is) at $\th \gl$ is at least
as large as $M(\th)-N(\th)$. On the other hand, the difference
$A(\th)-(M(\th)-N(\th))$ is at most the number $F(\th)$ of clones
activated in free steps until $\th\gl$, i.e.,
 \beq{mainob} M(\th)-N(\th) \leq A( \th) \leq  M(\th)-N(\th)  + F(\th).
 \eeq

As the cut-off lemma gives a concentration inequality for
$N(\th)$,
$$ \pr \bb[  \max_{\th:\tho \leq \th \leq 1 } |N(\th)
-(1-\th^\kk) \gl n|
 \geq \dd  \bb]
  \leq 2 e^{-\Omega(\min\{\dd, \frac{\dd^2}{(1-\tho)n}\})},
    $$
a concentration inequality for $M(\th)$ will be enough  to obtain
a similar inequality  for $B(\th):=M(\th) -N(\th)$. More
precisely, we will show that, under  appropriate hypotheses,
$$ \pr \bb[\max_{\th:\tho\leq \th \leq 1} \bb| M (\th) - (\gl-q(\th))
n \bb| \geq \dd \bb] \leq 2e^{-\Omega (\min\{ \dd,
\frac{\dd^2}{(1-\tho) n}\})},$$
 where
$$q(\th) =q(\th,\gl;P)= E\bb[ d_v
(\th) 1(D_v (\th) \in P )\bb]. $$  As $d_v (\th)$'s and $D_v
(\th)$'s are identically distributed, $q(\th)$ does not depend on
$v$. Recall also $p(\th)= \pr [ D_v(\th) \in P]$.

As we will see later, $B(\th)$ is very close to $A(\th)$. Hence, a
concentration inequality for $B(\th)$ is crucial.

\begin{lem} \label{mainl} (Main lemma)  In the $P$-process,
if $\tho< 1$ uniformly bounded from below by $0$,
$1-p(\tho)=O(1-\tho)$ and $p(\tho)=\Omega (1)$, then, for all
$\dd$ in the range $0< \dd \leq n$,
$$ \pr \bb[\max_{\th:\tho\leq \th \leq 1} \bb| |V(\th)| - p(\th)
n \bb| \geq \dd \bb] \leq 2e^{-\Omega (\min\{ \dd, \frac{
\dd^2}{(1-\tho) n}\})}, $$ and
$$ \pr \bb[\max_{\th:\tho\leq \th \leq 1} \bb| B(\th)
 - (\gl\th^\kk -q(\th)) n \bb| \geq \dd \bb]
\leq 2e^{-\Omega (\min\{ \dd, \frac{\dd^2}{(1-\tho) n}\})}.$$
\end{lem}


\mn {\bf Remark.} If $\dd \gg n^{1/2} \log n$, the proof of the
main lemma is much easier: Without the max, the two concentration
inequalities follow from the generalized Chernoff bound. Since the
bounds for the probabilities are much less than $1/n$ and there
are $O(n)$ meaningful $\th$'s, the first moment method gives the
inequalities. This is already enough to prove Theorems
\ref{mainpc},  \ref{sizeofcore}  and \ref{plrthm} provided $|\gl-
\gl_{\crt}| \gg n^{1/2} \log n$.   The full strength of the lemma
is needed when the $\log n$ factor is missing.

\mn

 For the proof, we first show that a concentration inequality
for
$$ |V(\th)|=\sum_{v\in V} 1( D_v (\th) \in P ),$$
which is a sum of i.i.d Bernoulli random variables with mean
$p(\th)$. More generally, we have

\begin{lem}\label{sobr}
Suppose $X_i$'s are  i.i.d Bernoulli random variables  with mean
$p$. Then, for $\dd>0$,
$$ \pr\bb[ \, \, \bb|\sum_{i=1}^m X_i -pm\bb|\geq \dd\bb ] \leq 2
e^{-\Omega (\min\{ \dd, \frac{\dd^2}{p(1-p) m}\})}. $$
\end{lem}

\pf Since
$$ E[X_i ] = p , ~~~
E[(X_i -p)^2] = p(1-p), $$ and, for $\rho$ with $|\xi|\leq
\xi_{_0}:=\log 2$,
$$
\bb| E[(X_i-p)^3 e^{\xi(1(X_i (\th) =0)-p)}] \bb| \leq 2p(1-p) ,
$$
we may   set $a_i= p(1-p) $ and $b_i= 2p(1-p)$. Applying the
generalized Chernoff bound, we have the desired inequality.

 \qed

\begin{cor}\label{vth} For $\th$ in the range $\tho\leq \th\leq 1$
and with the same hypotheses as in the main lemma,
$$ \pr \bb[ \bb||V_P (\th)| - p(\th)  n\bb| \geq \dd \bb]
\leq 2e^{-\Omega (\min\{\dd , \frac{ \dd^2}{(1-\tho) n}\})}.$$
\end{cor}
As in Example \ref{ex5}, Lemma \ref{gcb2}  yields a concentration
inequality for all of $V(\th)$'s:

\begin{lem} With the same hypotheses as in the main lemma,
$$ \pr \bb[\max_{\th:\tho\leq \th \leq 1} \bb| |V(\th)| - p(\th) n \bb|
\geq \dd \bb] \leq 2e^{-\Omega (\min\{ \dd, \frac{ \dd^2}{(1-\tho)
n}\})}.$$
\end{lem}

\pf   Observing that, for $\tho\leq \th\leq 1$,
$$  \sfrac{p(\tho)}{p(\th)}\bb| |V (\th)| - p(\th) n\bb|
\leq \bb| |V (\tho)|- p(\tho)n \bb| + \bb| |V(\tho)| -
 \sfrac{p(\tho)}{p(\th)}|V(\th)|\bb|, $$
we  set $\Gamma (\th) = | |V (\th)| - p(\th) n|,$   $$
 \psi =  \sfrac{1}{p(\tho)}\Gamma  (\tho),   \and \psi_\th = \sfrac{p(\th)}{p(\tho)}
 \bb| |V (\tho)| -   \sfrac{p(\tho)}{p(\th)} |V (\th)|\bb|. $$
Clearly,  $ \Gamma  (\th ) \leq \psi + \psi_\th$. Corollary
\ref{vth} gives \beq{psi1} \pr [ \psi \geq \dd/2] \leq
2e^{-\Omega(\min\{ \dd , \frac{ \dd^2}{(1-\tho) n}\})}. \eeq
 Suppose $\{X_{\th'}:=V(\th') \}_{\th\leq \th'\leq 1}$ is given, especially
$V (\th)$ is given. Then, since  $P$ is increasing, we may write
$|V(\tho)|$ as
$$ |V (\tho)|=\sum_{v\in  V(\th)}   1( D_v (\tho) \in
P) ,
$$ with
$$ \pr [  D_v(\tho) \in
P |\{X_{\th'} \}_{\th'\leq \th} ]= \pr [ D_v (\tho) \in P | v\in V
(\th) ] = \frac{p(\tho)}{p(\th)}=:p(\tho, \th).$$   Lemma
\ref{sobr} then gives
 \beq{psi2} \pr [ \psi_{\th} \geq \dd/2 ]
\leq 2e^{-\Omega(
     \min \{ p(\tho, \th) \dd , \frac{p(\tho, \th)\dd^2}{(1-p(\tho,\th)) |V
     (\th)|}\} )}\leq 2e^{-\Omega(\min\{p(\tho\!) \dd , \frac{p(\tho \!)
\dd^2}{(1-p(\tho\!)) n}\})} \leq 2e^{-\Omega(\min\{\dd , \frac{
\dd^2}{(1-\tho) n}\})}. \eeq Lemma \ref{gcb2}  together with
\raf{psi1} and \raf{psi2} yields the desired inequality.

\qed

\bn

We now estimate $M(\th)$. First, since $\sum_{v\in V} \bd_v(\th)$
is a Poisson random variable with mean $(1-\th)\gl n $,
  \beq{firstsum} \pr\bbr  \bb|\sum_{v\in V} \bd_v(\th) - (1-\th)\gl n \bb| \geq \dd/2\bb]
\leq 2 e^{-\min\{ \dd, \frac{\dd^2}{(1-\th) n}\}}. \eeq For the
second sum in \raf{newm}, observe that
$$ \sum_{v\in V} d_v (\th) 1(v\not\in V (\th))
=  \sum_{v\in V}     d_v(\th)  1(D_v(\th) \not \in P)$$ is a sum
of i.i.d random variables with
$$ E[  d_v(\th) 1(D_v(\th)
\not \in P)]=E[  d_v(\th) -  d_v (\th) 1(D_i (\th) \in P ) \bb] =
\th\gl- q(\th).
$$
  Moreover, since $P$ is an
increasing property and $d_v(\th) $ is a Poisson $\th\gl $ random
variable, FKG inequality (see e.g.  Chapter 6 of \cite{AS}) gives
$$
E[  (d_v(\th) 1(D_v(\th) \not \in P))^i] \leq E[   d_v(\th)^i ]
\pr [D_i(\th) \not \in P ] = O(1-p(\th)),
$$
for all fixed $i$, e.g. $i=1,2,3$. Thus one may take $\xi_{_0}=1$
and $a_i, b_i=\Theta(1-p(\th))$ to satisfy all the conditions to
apply the generalized Chernoff bound and to obtain
 $$ \pr \bbr  \bb|\sum_{v\in V} d_v (\th) 1(v\not\in V (\th))
   - (\th\gl-q(\th)) n \bb| \geq \dd \bb] \leq
2e^{-\Omega(\min\{ \dd, \frac{\dd^2}{(1-p(\th))n}\})}.
 $$
 This together with \raf{firstsum} implies
that if  $1-p(\th) = O(1-\th)$, then
   \beq{smth}
 \pr \bbr  \bb|M(\th)
   - (\gl-q(\th)) n \bb| \geq \dd \bb] \leq
2e^{-\Omega(\min\{ \dd, \frac{\dd^2}{(1-\th)n} \})}.
 \eeq

As mentioned, the following lemma  is enough to prove the main
lemma.

\begin{lem} \label{amth} With the same hypotheses as in the main lemma,
$$ \pr \bb[\max_{\th:\tho\leq \th \leq 1} \bb| M (\th) - (\gl-q(\th))
 n \bb| \geq \dd \bb]
\leq 2e^{-\Omega (\min\{ \dd, \frac{\dd^2}{(1-\tho) n}\})}.$$
\end{lem}

\pf Clearly,
 $$| M(\th) -
 (\gl-q(\th)) n|\leq |M(\tho)-  (\gl-q(\tho)) n|+
 |M(\tho)-M (\th) - (q(\th)-q(\tho)) n|. $$
Let $\Gamma (\th) =| M (\th) -
 (\gl-q(\th)) n|$,
 $$ \psi = \Gamma (\tho),~~ \psi_{\th} =|M (\tho)-M (\th) - (q(\th)-q(\tho))
 n|, $$
and $\Phi_\th$ is the event $|V (\th) - p(\th)n| \leq
\frac{p(\tho)\dd}{4\gl}$. Then, \raf{smth} gives
$$ \pr [ \psi \geq \dd/2 ] \leq 2e^{-\Omega(\{\dd,
\frac{\dd^2}{(1-\tho) n}\})}. $$ For $\psi_\th$, suppose
$\{X_{\th'} :=(V(\th'), M (\th'))\}_{\th\leq \th' \leq 1}$ is
given. Using $$ M(\th) = \sum_{v\in V} \bd_v (\th) + d_v (\th ) 1
(v\not \in V (\th))=  \sum_{v\in V}   d_v (1) -  d_v(\th)1 (v \in
V (\th)) ,$$  we obtain
 \bean M( \tho) - M (\th) &=& \sum_{v\in  V} d_v (\th)
 1(v\in  V (\th))
- d_v (\tho) 1(v\in V (\tho)) .\eean
 Once $V(\th)$ is given, the
distributions of $$Y_v := d_v (\th)1(v \in V (\th)) -
 d_v (\tho)  1(v \in V (\tho))$$  depend on
neither $\{M(\th')\}_{\th\leq \th'\leq 1} $ nor $\{V(\th')\}_{\th<
\th'\leq 1} $ and hence, for $v\in  V (\th)$,
 \bean
 E\bb[  Y_v
\bb| \{ X_{\th'}\}_{\th\leq \th' \leq 1}\bb] & =& E\bb[d_v (\th) -
d_v (\tho)  1(v\in V (\tho))\bb| v\in  V (\th) \bb] \\
&= &
 \frac{q(\th) - q(\tho)}{p(\th)} .
\eean  
If $v \not\in V (\th)$, $Y_i=0$   since $P$ is increasing.

Also, for $v \in V(\th)$, we may write
 $$ Y_v =
d_v (\th) - d_v (\tho) + d_v (\tho) 1(v\not\in V (\tho))
$$
and \bean
 E\bbr Y_v^2 \bb|v \in V (\th)\bb]
&\leq &
 2 E\bb[  ( d_v (\th) - d_v (\tho) )^2 \bb|v \in V (\th) \bb]
  +  2E\bbr
d_v (\tho)^2 1(v\not\in V (\tho)) \bb| v\in V (\th) \bb].
 \eean
 First, for $j=1,2,3$,
$$
 E\bbr  ( d_v (\th) - d_v (\tho) )^j \bb|
v\in  V (\th) \bb]
  \leq
p(\th)^{-1} E\bbr (  d_v (\th) - d_v (\tho) )^j \bb] = O(\th-\tho)
= O(1-\tho) $$ for $p(\th)\geq p(\tho) =\Omega(1)$ and $ d_v (\th)
- d_v (\tho) $ is a Poisson  random variable with mean
$(\th-\tho)\gl =O(\th-\tho)$. For the second term, FKG inequality
gives
 \bean
E\bbr
  \bb(
d_v (\tho) 1(v\not\in V (\tho))\bb)^j \bb| v\in  V (\th) \bb]
&\leq & p(\th)^{-1} E\bbr
d_v (\tho)^j 1(v\not\in  V (\tho))\bb] \\
 &\leq&
 p(\tho)^{-1} E[
d_v (\tho)^j]\pr [v\not\in V (\tho)]\\
&=& O(1-p(\tho)) = O(1-\tho),
 \eean
 for $j=1,2,3$.
Therefore,
   $$ E\bbr \bb( Y_i - E[Y_i] \bb)^2 \bb| \{ X_{\th'}\}_{\th\leq \th' \leq 1}\bb] \leq
 E\bb[  Y_i^2 \bb|  \{ X_{\th'}\}_{\th\leq \th' \leq 1}\bb]= O(1-\tho).
$$
Similarly, for $\xi$ in the range $|\xi|\leq \xi_{_0}=1$, it is
not hard to show
 $$ \bb| E\bbr (Y_i-E[Y_i])^3 e^{\xi(Y_i-E[Y_i])}  \bb| \{ X_{\th'}\}_{\th\leq \th' \leq 1}\bb]
 \bb|
 = O(1-\tho). $$
Applying the generalized Chernoff bound, we have
$$ \pr \bbr \bb| \sum_{v\in V} Y_v -
\frac{q(\th)-q(\tho)}{p(\th)} | V (\th) | \bb| \geq \dd/4  \bb| \{
X_{\th'}\}_{\th\leq \th' \leq 1}\bb] \leq 2e^{-\Omega(\min \{\dd,
\frac{\dd^2}{1-\tho}\})}.
$$

Finally, as    the event $\Phi_\th$ guarantees
$$ \frac{q(\th)-q(\tho)}{p(\th)}  \bb| | V (\th) | -
p(\th)n\bb| \leq \dd/4  $$ for  $p(\tho) \leq p(\th)$ and $q(\th)
\leq \gl$, we  have
  \bean  & & 1(\Phi_\th) \pr \bbr \bb| \sum_{v\in v} Y_v -
(q(\th)-q(\tho))n \bb| \geq \dd/2 \bb| \{ X_{\th'}\}_{\th\leq \th'
\leq 1}\bb]  \\ & & \leq \pr \bbr \bb| \sum_{v\in V} Y_v -
\frac{q(\th)-q(\tho)}{p(\th)} | V (\th) | \bb| \geq \dd/4 \bb|
\{ X_{\th'}\}_{\th\leq \th' \leq 1}\bb]\\
& & \leq 2e^{-\Omega(\min \{\dd, \frac{\dd^2}{1-\tho}\})}.  \eean
Lemma \ref{gcb2}  yields the desired inequality.

\qed

\section{Cores of Random Hypergraphs}
\label{score}

This section is for the proof of  Theorem \ref{sizeofcore}. Let
$\gl>0$ and $H(\gl)=\hh (n,p)$, where $\gl =p{n-1 \choose r-1}$.
Let  the property $P= \{ (d_1, d_2) : d_1 \geq t \}$. Then
$$ p(\th) = Q(\th\gl, t) , \and q(\th) = \th\gl Q(\th \gl, t-1)
. $$ The main lemma gives

\begin{cor}\label{core}
For $ \tho \leq1$  uniformly bounded from below by $0$ and $\dd$
in the range $0< \dd \leq n$,
$$ \pr \bb[\, \,  \max_{\th: \tho\leq \th \leq 1} \bb| |V(\th)| - Q(\th\gl,t) n \bb|
\geq \dd  \bb] \leq  2 e^{-\Omega (\min\{ \dd , \frac{\dd^2}{
n}\})},
$$ and
$$ \pr \bbr \max_{\th:\tho\leq \th\leq 1}  \bb| B(\th)  -(\th^{\frac{1}{k-1}} -
Q(\th\gl, t-1))\th\gl  n \bb| \geq \dd \bb] \leq 2 e^{-\Omega
(\min\{ \dd, \frac{\dd^2}{ n}\})}.$$
\end{cor}

\mn

\mn {\em Subcritical Region}:  For $\gl= \gl_\crt -\gs$, $\gs \gg
n^{-1/2}$, and $\tho= \gd /\gl_{\crt} $ with $\gd= 0.1  $, it is
easy to see that there is a constant $c>0$ such that
$$
(\th^\km  - Q(\th\gl, t-1))\th \gl  n \geq c \gs n , ~~\mbox{for
all $\th$ in the range $\tho\leq \th \leq 1$}.$$ Let $\tau$ be the
first time the number $A(\th)$ of active clones at $\th\gl$
becomes $0$. Then the second part of Corollary  \ref{core} gives
  \begin{align*} \pr [ \tau \geq \tho] &\leq  \pr[ B(\th)=0~~\mbox{for some $\th$
with $\tho \leq \th \leq 1$}] \\ &\leq
 \pr \bbr
\max_{\th:\tho\leq \th\leq 1}  \bb| B(\th)  -(\th^\km - Q(\th\gl,
t-1))\th\gl  n \bb| \geq c\gs n  \bb] \\
&\leq  2 e^{-\Omega (\gs^2 n)}.\end{align*}
  As $\tho \gl \leq \tho \gl_{\crt} =\gd$, and hence $Q(\tho \gl, t) \leq
\gd/2$ for $t\geq 2$, the first part of Corollary \ref{core}
yields
$$ \pr [ |V_t(H_{PC} (n,p\,;k))| \geq \gd n ] \leq \pr [ \tau \geq \tho] + \pr [
|V(\tho)| \geq \gd n] \leq 2e^{-\Omega(\gs^2 n)}.  $$  Therefore,
Theorem \ref{pc} implies that
$$ \pr [ |V_t( H(n,p\,;k))|  \geq \gd n] \leq 2e^{-\Omega(\gs^2 n)}.$$

 To complete the proof, we  observe that the $t$-core of size $i$
has at least $t i /k$ edges. Let $Z_i$ be the number  of subgraphs
on $i$ vertices with
 at least $t i /k$ edges, $i=i_{_0}, ..., \gd n$, where $i_{_0} = i_{_0} (k,t)$ is  the least
 $i$ such that ${i\choose k} \geq ti/k$. Then, in $H(n,p\,;k)$,
$$ E[Z_i ] \leq {n \choose i} {{i \choose k}\choose  ti/k }
 p^{ti/k} \leq \frac{n^i}{i!} \frac{i^{ti}}{ (ti/k ) !}
  p^{ti/k}=: L_i,
 $$
where $ti/k$ actually means $\lceil ti/k \rceil$. Observe that
$$ \frac{L_{i+k}}{L_{i}} =O\bb(\frac{n^k}{i^k}
\frac{i^{kt}}{i^t} n^{-(k-1)t}\bb)
=O\bb(\bb(\frac{i}{n}\bb)^{(k-1)t-k }\bb)=O(\gd^{(k-1)(t-1)-1}).
$$
That is, $L_{i+k}/ L_i$ exponentially decreases. Since
$$ L_i = O(n^i n^{-i(k-1)t/k} ) = O(  n^{-i(t-1-t/k)}),$$ for
$i=i_{_0}, ..., i_{_0}+k-1$, it follows that
$$ \pr [ V_t( H(n,p\,;k))\not=\emptyset] \leq 2e^{-\Omega(\gs^2 n)}+O(  n^{-i_{_0}(t-1-t/k)}),
$$
as desired.

\qed

\bn {\em  Supercritical Region}: We will prove the following
theorem.
\begin{thm} \label{stop} If $\gl:=p{n-1 \choose k-1} = \glc + \gs$ with
 $\gs\gg n^{-1/2}$ and $0<\gd\leq 1$, then,
 with
probability $1-2e^{-\Omega(\min\{ \gd^2 \gs n, \gs^2 n\} )}$,
$V_t=V_t(H_{PC}(n,p\, ; k))$ satisfies
  $$ Q(\thl \gl, t) n - \gd n \leq |V_t|\leq   Q(\thl \gl, t) n + \gd n,
  $$
and  the degrees of vertices of  the $t$-core are i.i.d
$t$-truncated Poisson random variables with parameter $\gL_t:=\thl
\gl + \gb$ for some $\gb$ with $|\gb| \leq \gd$. Moreover, the
distribution of the $t$-core is the same as that of the
$t$-truncated Poisson cloning model with parameters $|V_t|$ and
$\gL_t$.
\end{thm}

Recall that $\thl$ is the largest solution for the equation
$$ \th^\km - Q(\th \gl, t-1)=0. $$

\pf First, it is not hard to check that  there are constant
$\co,\ct>0$ such that, for $\th$ in the range $\thl \leq \th \leq
1$,
$$ \th^\km - Q(\th \gl, t-1) \geq \co \gs^{1/2} (\th-\thl), $$
and, for $\th$ in the range $\thl- \ct \gs^{1/2}\leq \th \leq
\thl$,
$$ \th^\km - Q(\th \gl, t-1) \leq - \co \gs^{1/2} (\thl-\th ). $$

Let $\tau$ be the largest $\th$ with $A(\th) =0$. Then  $V(\tau)$
is the $t$-core of $H_{PC}(n,p\,; k)$. For $\tho=\thl+ \gd$ and
$\thto=\thl-\min\{ \gd, \ct \gs^{1/2}\}$ with $0< \gd \le 1 $,
Corollary  \ref{core} gives
  \begin{align*}  \pr [\tau \geq \tho] &\leq
    \pr[ B(\th)=0~~\mbox{for some $\th$
with $\tho \leq \th \leq 1$}] \\
&\leq  \pr \bbr \max_{\th:\tho\leq \th\leq 1}
 \bb| B(\th)  -(\th^\km  -
Q(\th\gl, t-1))\th\gl  n \bb| \geq \co \gs^{1/2} \gd n  \bb] \\
& \leq  2 e^{-\Omega (\gd^2 \gs  n)},\end{align*}
and
 \begin{align*}  \pr [\tau < \thto] &\leq \pr [ B(\thto) >0] \\
 &\leq \pr \bbr \bb| B(\thto)  -(\thto^\km  -
Q(\thto\gl, t-1))\thto \gl  n \bb| \geq \co \gs^{1/2} \min\{\gd,
\ct \gs^{1/2}\} n\bb] \\ &\leq  2 e^{-\Omega (\min\{\gd^2 \gs n,
\gs^2 n\} )}. \end{align*}

Since  $\frac{d}{d\th} Q(\th\gl, t)= \gl P(\th\gl,t-1)\leq \gl$,
we have
 $$ Q(\tho \gl, t)  \leq  Q(\thl \gl, t) +  \gl \gd, \and
Q(\thto  \gl, t)  \geq  Q(\thl  \gl, t) -  \gl \gd, $$ and
Corollary  \ref{core} implies that
$$ \pr[ V(\tho )- Q(\thl \gl, t) n \geq 2 \gl \gd n ] \leq
2e^{-\Omega (\gd^2 n)}, $$ and $$
 \pr[ V(\thto )- Q(\thl \gl, t) n \leq -2 \gl \gd n ] \leq
2e^{-\Omega (\gd^2 n)}. $$ Therefore,
$$ \pr[ |\tau - \thl| > \gd ] \leq \pr [ \tau \geq \tho] + \pr [ \tau \leq \thto]
\leq 2e^{-\Omega (\min\{\gd^2 \gs n,\, \gs^2 n\})},$$
 and,  replacing $\gd$ by
$\frac{\gd}{2\gl}$,
 \begin{align*} \pr [ | V(\tau)- Q(\thl \gl, t) n | \geq \gd n   ]
&\leq \pr [ \tau \geq \tho] + \pr [ \tau \leq \thto] + 2e^{-\Omega
(\gd^2 n)} \\
&\leq 2e^{-\Omega (\min\{\gd^2 \gs n,\, \gs^2 n\})}.
\end{align*}

Clearly, once $V(\tau)$  and $\gL_t:= \tau \gl$ are given,  the
residual degrees $d_v (\tau)$, $v \in V(\tau)$, are i.i.d
$t$-truncated Poisson random variables with parameter $\gL_t $.
\qed

\bn

Once $V_t$ and $\gL_t$ are given, $|V_t (i)|$, $i\geq t$, is the
sum of i.i.d Bernoulli random variables with mean $ p_{_i} (\gL_t)
:= \sfrac{P(\gL_t, i)}{Q(\gL_t, t)}. $ Similarly, the size of $W_t
(i) =\cup_{j\geq  i} V_t (j)$ is the sum of i.i.d Bernoulli random
variables with mean $ q_{_i} (\gL_t) := \sfrac{Q(\gL_t,
i)}{Q(\gL_t, t)}. $ Applying the generalized Chernoff bound (Lemma
\ref{uni}), we have
$$ \pr \bbr \bb|  |V_t (i)| - p_{_i}\!(\gL_t)  |V_t| \bb|\geq  \gd   |V_t|   \bb| V_t, \gL_t \bb]
\leq 2e^{-\Omega(\gd^2  |V_t|) },
$$
and
$$ \pr \bbr \bb|  |W_t (i)| - q_{_i}\!(\gL_t)  |V_t| \bb|\geq  \gd   |V_t|   \bb| V_t, \gL_t \bb]
\leq 2e^{-\Omega(\gd^2  |V_t|) }.
$$
Combining these with Lemma \ref{stop} and using
$$ |P(\rho, i) - P(\rho', i)|\leq |\rho-\rho'|, \and
|Q(\rho, i) - Q(\rho', i)|\leq  |\rho-\rho'|, $$   we obtain, for
any $i$,
$$ \pr \bbr \bb|  |V_t (i)| -  P(\thl \gl, i) n
\bb|\geq  \gd  n   \bb] \leq 2e^{-\Omega(\min\{ \gd^2 \gs n, \gs^2
n\} )},
$$
and
$$ \pr \bbr \bb|  |W_t (i)| -  Q(\thl \gl, i) n
\bb|\geq  \gd   n   \bb] \leq 2e^{-\Omega(\min\{ \gd^2 \gs n,
\gs^2 n\} )}.
$$

In particular, as $ \thl= \thc +\Theta(\gs^{1/2}) , $ for
uniformly bounded $\gs$,
 it follows that, for $\gl= \gl_{\crt} +\gs$,
$$   |V_t (i)| = (1+O(\gs^{1/2}))^i   P(\th_{\crt}\gl_{\crt} , i) n
+O\bb((n/\gs)^{1/2}\log n\bb),
$$
with probability $1- 2e^{-\Omega(\min\{\log^2 n , \gs^2 n\} )}$.

The last part of Theorem \ref{sizeofcore} does not follow from
Theorem \ref{stop} as $H(n,p\,;k)$ and $H_{PC} (n,p\,;k)$ do not
have the same distribution. We may directly  prove it instead.

For a hypergraph $H$, let $\tilde{H}$ be the hypergraph obtained
from $H$ by removing all edges in the $t$-core. Then, $\tilde{H}$
 has no edge that is entirely in $V_t(H)$, otherwise, the $t$-core
 becomes larger. Thus, the union of $\tilde{H}$ and any
 simple hypergraph on $V_t(H)$ is also simple. Therefore,
 conditioned on $\tilde{H}=\tilde{H}(n,p\,\;k)$, two hypergraphs $H_1$ and $H_2$
 on $V_t
 (H(n,p\,;k))$ with the same degree sequence
  are equally likely
 to be the $t$-core of $H(n,p\,;k)$: Notice that
 $$ \pr [ H(n,p\,;k)= \tilde{H} \cup H_1] = p^{\tilde{m}+m_{_1}}
 (1-p)^{{n \choose k} -\tilde{m}-m_{_1}}, $$
 and
  $$ \pr [ H(n,p\,;k)= \tilde{H} \cup H_2] = p^{\tilde{m}+m_{_2}}
 (1-p)^{{n \choose k} -\tilde{m}-m_{_2}}, $$
where $\tilde{m}$, $m_{_1}$ and $m_{_2}$ are the numbers of edges
in $\tilde{H}$, $H_1$ and $H_2$, respectively. Clearly,
$m_{_1}=m_{_2}$ as the degree sequences of $H_1$ and $H_2$ are the
same.

\section{The pure literal algorithm for the random $k$-SAT}
\label{sat}  To analyze the pure literal algorithm for a random
$k$-SAT, $k\geq 3$, it is necessary to consider a pair of degrees.
We consider this in a generalized framework. Starting with some
terminology, a  sequence of generalized degrees  is larger than or
equal to another with the same length if so is each pair of
corresponding generalized degrees. A property for sequences of
generalized degrees  is a set of generalized degree sequences. A
property $P$ is increasing if sequences of generalized degree
larger than an element  in $P$ are also in $P$.

Given a Poisson $\gl$-cell on the set $V$ of $n$ vertices, let
$\c=\{ C_i\} $ be an equipartition of the vertex set $V$. Then,
$D_i (\th) := (d_v (\th), \bd_v (\th))_{v\in C_i}$ are i.i.d
random variables. In particular, for any property $P$, the events
$D_i (\th) \in P$ are independent and occur with the same
probability, say $p(\th, \gl; P,\c)$, or simply $p(\th)$. For the
pure literal algorithm of  the random $k$-SAT problem, the
property is the set of pairs of generalized degrees $(d_1,d_2)$,
$(d'_1, d'_2)$ with $d_1, d_1' \geq 1$.

For an increasing property $P$ and an equipartition $\C=\{C_i
\}_{i=1,..., m}$, the $(P, \C)$-process is defined as follows.
Construct the Poisson $\gl$-cell as described in Section
\ref{grpm}, where $\gl= p{n-1 \choose k-1}$. The $(P, \C)$-process
is a generalization of the $P$-process.

\mn {\bf The  $(P, \C)$-process}: Initially, the cut-off value
$\gL= \gl$. Activate all vertices $v\in C_i$  with $D_i (1)\not
\in P$.
 All clones of the activated vertices are
activated too. Put activated  clones in a stack in an arbitrary
order. However, this does not mean that the clones are removed
from the $\gl$-cell.

\mn (a)   If the stack is empty, go to (b). If the stack  is
nonempty, choose the first clone in the stack and move the cut-off
line to the left  until the largest $k-1$ unmatched clones,
excluding the chosen clone, are found. (So, the cut-off value
$\gL$ keeps decreasing.) Then, match the $k-1$ clones to the
chosen clone. Remove all matched clones from the stack  and
repeat. A vertex in $C_i$ that has not been activated is to be
activated as soon as $D_i (\gL/\gl)\not\in P$. This can  be done
even before all $k-1$ clones are found.  Its unmatched clones are
to be activated too and put into the stack immediately. Clones
found while moving the cut-off line are also in the stack until
they are matched.

\mn (b) Activate all vertices in the  first $C_i$ no vertex of
which has not been  activated. Its clones are activated too. Put
those clones into the stack. Then, go to (a).

\mn Clones in the stack are called active.
The steps carried by the instruction described in (b) are called
free steps as it is free to artificially activate a vertex.

When the cut-off line is at $\th\gl$, all $\th\gl$-large clones
are matched or will be matched at the end of the step and all
vertices in $C_i$ with $D_i (\th) \not\in P$ are activated. All
other vertices can be activated only by free steps. Let
$V(\th)=V_{(P,\c)} (\th)$ be  the union of $C_i$ with $D_i (\th)
\in P$, and let $M(\th)=M_{(P,\c)}(\th)$ be the number of $\th
\gl$-large clones plus the number of $\th \gl$-small clones of
$v\not\in V(\th)$. That is,
$$ M(\th) = \sum_{v\in V} \bd_v (\th) + d_v (\th ) 1 (v\not \in V
(\th))= \sum_{i=1}^m \sum_{v\in C_i} \bd_v (\th) + d_v (\th) 1
(D_i (\th ) \not\in P). $$

Recalling that $N(\th)$ is the number of matched clones until the
cut-off line reaches $\th \gl$, the number $A(\th)$  of active
clones (when the cut-off $\gL$ is) at $\th \gl$ is at least as
large as $M(\th)-N(\th)$. On the other hand, the difference
$A(\th)-(M(\th)-N(\th))$ is at most the number $F(\th)$ of clones
activated in free steps until $\th\gl$, i.e.,
$$ M(\th)-N(\th) \leq A( \th) \leq  M(\th)-N(\th)  + F(\th). $$

As the cut-off lemma gives a concentration inequality for
$N(\th)$,
$$ \pr \bb[  \max_{\th:\tho \leq \th \leq 1 } |N(\th)
-(1-\th^\kk) \gl n|
 \geq \dd  \bb]
  \leq 2 e^{-\Omega(\min\{\dd, \frac{\dd^2}{(1-\tho)n}\})},
    $$
a concentration inequality for $M(\th)$ will be enough  to obtain
a similar inequality  for $B(\th):=M(\th) -N(\th)$. More
precisely, we will show that, under  appropriate hypotheses,
$$ \pr \bb[\max_{\th:\tho\leq \th \leq 1} \bb| M (\th) - (\gl-q(\th)) n \bb| \leq \dd \bb]
\leq 2e^{-\Omega (\min\{ \dd, \frac{\dd^2}{(1-\tho) n}\})},$$
 where
$$q(\th) =q(\th,\gl;P,\c)= E\bb[\bb( \frac{1}{|C_i|}\sum_{v\in C_i} d_v
(\th)\bb) 1(D_i (\th) \in P )\bb]. $$  As $D_i (\th)$ are
identically distributed, $q(\th)$ does not depend on $i$. Recall
also $p(\th)=p(\th,\gl;P,\c) = \pr [ D_i (\th) \in P]$. Here is a
generalization of the main lemma. Its proof is quite similar to
that of the main lemma and it is presented in the Appendix.

\begin{lem} \label{gmain} (Main lemma: generalized version)  In the
$(P, \c)$-process described above,
if $\tho< 1$ uniformly bounded from below by $0$, and
$|C_1|=O(1)$, $1-p(\tho)=O(1-\tho)$ and $p(\tho)=\Omega (1)$,
then, for all $\dd$ in the range $0< \dd \leq n$,
$$ \pr \bb[\max_{\th:\tho\leq \th \leq 1} \bb| |V(\th)| - p(\th) n \bb| \leq \dd \bb]
\leq 2e^{-\Omega (\min\{ \dd, \frac{ \dd^2}{(1-\tho) n}\})}, $$
and
$$ \pr \bb[\max_{\th:\tho\leq \th \leq 1} \bb| B(\th)
 - (\gl\th^\kk -q(\th)) n \bb| \leq \dd \bb]
\leq 2e^{-\Omega (\min\{ \dd, \frac{\dd^2}{(1-\tho) n}\})}.$$
\end{lem}

Let $\gl>0$ and $F(\gl)=F_{PC}  (n,p\, ; k)$, where $\gl =p  {2n-1
\choose k-1}$. As mentioned, we take the property $P= \{ ((d_1,
d_2), (d_1', d_2')) : d_1, d_1'  \geq 1 \}$ and  $C_i = \{x_i,
\bar{x}_i \}$. Then
$$ p(\th) = Q^2(\th\gl, 1) = (1-e^{-\th\gl})^2  , \and q(\th) =
\th\gl(1-e^{-\th\gl})
. $$ Let $X(\th)$ be the set of variables $x$ with both of $d_x
(\th)$ and $ d_{\bar{x}} (\th) $ larger than $0$. Then the main
lemma and $|V(\th)| = 2|X(\th)|$ give

\begin{cor}\label{plr}
For $ \tho \leq 1$  uniformly bounded from below by $0$ and $\dd$
in the range $0< \dd \leq n$,
$$ \pr \bb[\, \,  \max_{\th: \tho\leq \th \leq 1}  \bb| |X(\th)| - (1-e^{-\th\gl})^2  n \bb|
\geq \dd  \bb] \leq  2 e^{-\Omega (\min\{ \dd , \frac{\dd^2}{
n}\})},
$$
and
$$ \pr \bbr \max_{\th:\tho\leq \th\leq 1}  \bb| B(\th)  -2(\th^{\frac{1}{k-1}} -
(1-e^{-\th\gl}))\th\gl  n \bb| \geq \dd \bb] \leq 2 e^{-\Omega
(\min\{ \dd, \frac{\dd^2}{ n}\})}.$$
\end{cor}

\bn

As
 $1-e^{-\th\gl}=Q(\th\gl, t-1)$ with $t=2$, a similar  argument  used in the
 previous section may be applied to prove Theorem \ref{plrthm}.

 \mn {\em Subcritical Region}:  For $\gl= \gl_\crt -\gs$, $\gs
\gg n^{-1/2}$, let  $\tho= \gd /\gl_{\crt} $ with $\gd=0.1$, and
let $\tau$ be the first time the number $A(\th)$ of active clones
at $\th\gl$ becomes $0$. Since
$$
2(\th^\km  - (1-e^{\th\gl}))\th \gl  n \geq c \gs n , ~~\mbox{for
all $\th$ in the range $\tho\leq \th \leq 1$},$$ for a constant
$c>0$, the second part of Corollary \ref{plr} gives
  \bean \pr [ \tau \geq \tho] &=&  \pr[ B(\th)=0~~\mbox{for some $\th$
with $\tho \leq \th \leq 1$}] \\ &\leq&
 \pr \bbr
\max_{\th:\tho\leq \th\leq 1}  \bb| B(\th)  -2(\th^\km - (1-e^{-\th\gl}))\th\gl  n \bb| \geq c\gs n  \bb] \\
&\leq&  2 e^{-\Omega (\gs^2 n)}.\eean
  As $(1-e^{-\tho\gl})^2 \leq
\gd^2$, this and the first part of Corollary  \ref{core} yield
$$ \pr [ |X_R (F_{PC} (n,p\, ; k))| \geq \gd n ] \leq \pr [ \tau \geq \tho] + \pr [
|X(\tho)| \geq \gd n] \leq 2e^{-\Omega(\gs^2 n)}.  $$  Therefore,
Theorem \ref{equiv} implies that
$$ \pr [ |X_R ( F(n,p\, ; k))|  \geq \gd n] \leq 2e^{-\Omega(\gs^2 n)}.$$

 To complete the proof, we  observe that the residual formula  on
 $i$ variables
has at least $2 i /k$ clauses. Let $Z_i$ be number of subformulas
on $i$ variables  with
 at least $2 i /k$ clauses, $i=i_{_0}, ..., \gd n$, where
  $i_{_0} = i_{_0} (k)$ is the least
 $i$ such that $2^k {i\choose k} \geq 2i/k$. Then, in
 $F(n,p\,;k)$,
 $$ E[ Z_i] \leq {n \choose i} {2^k {i \choose k}\choose  2i/k }
 p^{2i/k} \leq \frac{n^i}{i!} \frac{(2i)^{2i}}{ (2i/k ) !}
  p^{2i/k}=: L_i,
$$
where $2i/k$ actually means $\lceil 2i/k \rceil$. Observe that
$$ \frac{L_{i+k}}{L_{i}} =O\bb(\frac{n^k}{i^k}
\frac{(2i)^{2k}}{i^2} n^{-2(k-1)}\bb)
=O\bb(\bb(\frac{i}{n}\bb)^{k-2}\bb)=O(\gd^{k-2}).
$$
That is,  $L_{i+k}/ L_i$ exponentially decreases for $k\geq 3$.
Since
$$ L_i = O(n^i n^{-2i(k-1)/k} ) = O(  n^{-i(1-2/k)}),$$
for $i=i_{_0}, ..., i_{_0}+k-1$, we have
$$ \pr [ | X_R( F(n,p\,;k))|\not=\0] \leq 2e^{-\Omega(\gs^2 n)}+O(  n^{-i_{_0}(1-2/k)}),
$$
as desired.

\qed

\bn {\em  Supercritical Region}: Applying the same argument used
to prove Theorem \ref{stop} in the previous section, we may easily
obtain
\begin{thm} \label{stopplr} If $\gl= p{2n-1 \choose k-1} \geq \glc + \gs$ with
$\gs \gg n^{-1/2}$ and $0<\gd\leq 1$, then, with probability
$1-2e^{-\Omega(\min\{ \gd^2 \gs n, \gs^2 n\} )}$, $|X_R(n,p\,;k)|$
satisfies
  $$ (1-e^{-\thl\gl})^2 n - \gd n \leq |X_R(n,p\,;k)| \leq
  (1-e^{-\thl\gl})^2 n + \gd n,
  $$
and  the degrees of literals  of $X_R(n,p\,;k)$   are i.i.d
$1$-truncated Poisson random variables with parameter $\gL_R:=\thl
\gl + \gb$ for some $\gb$ with $|\gb| \leq \gd$. Moreover, the
distribution of the residual formula on $X_R(n,p\,;k)$ is  the
same as that of the $1$-truncated Poisson cloning model with
parameters $|X_R(n,p\,;k)|$ and $\gL_R$.
\end{thm}

\pf The proof is almost identical to that of Theorem \ref{stop}
with $t=2$. \qed

\mn

 Once $X_R:=X_R(n,p\,;k)$ and $\gL_R$ are given, $|X_R
(i,j)|$, $i,j\geq 1$, is the sum of i.i.d Bernoulli random
variables with mean $ p_{_{i,j}} (\gL_t) := \sfrac{P(\gL_t,
i)P(\gL_t, j)}{(1-e^{-\gL_R})^2}. $ Applying Lemma \ref{sobr}, the
concentration inequality for a sum of i.i.d Bernoulli random
variables, we now have
$$ \pr \bbr \bb|  |X_R (i,j)| - p_{_{i,j}}\!(\gL_t)  |X_R| \bb|\geq  \gd   |X_R|
  \bb| X_R, \gL_R \bb]
\leq 2e^{-\Omega(\gd^2  |X_R|) }.
$$
Combining this with Lemma \raf{stop} and using
$$ |P(\rho, i) - P(\rho', i)|\leq |\rho-\rho'|, $$   we obtain, for
any $i$,
$$ \pr \bbr \bb|  |X_R (i,j)| -  P(\thl \gl, i)P(\thl \gl, j) n
\bb|\geq  \gd  n   \bb] \leq 2e^{-\Omega(\min\{ \gd^2 \gs n, \gs^2
n\} )}.
$$
Similarly,
$$ \pr \bbr \bb|  |Y_R (i,j)| -  Q(\thl \gl, i)Q(\thl\gl,j) n
\bb|\geq  \gd   n   \bb] \leq 2e^{-\Omega(\min\{ \gd^2 \gs n,
\gs^2 n\} )}.
$$

The last statement of Theorem \ref{plrthm} follows by the same
argument used in the previous section.

\section{The Emergence Of the Giant Component}
\label{giant}

\renewcommand{\AA}{\tilde{A}}

In this section, we prove  Theorem \ref{mainpc}. Let the property
$P= \{ (d_1, d_2) : d_{_2} =0 \}.$ Then
$$p(\th) = e^{-(1-\th)\gl}
, \and q(\th) = \th\gl e^{-(1-\th)\gl}
, $$ and the main lemma gives

\begin{cor}\label{mainlem}
For $ \tho \leq1 $  uniformly bounded from below  by $0$ and $\dd$
in the range $0< \dd \leq n$,
$$ \pr \bb[\, \,  \max_{\th: \tho \leq \th \leq 1 } \bb| |V (\th)| - e^{-(1-\th)\gl}  n \bb|
\geq \dd  \bb] \leq  2 e^{-\Omega (\min\{ \dd , \frac{\dd^2}{(1-
\tho ) n}\})},
$$
and
$$ \pr \bbr \max_{\th: \tho \leq \th \leq 1 }  \bb| B  (\th)
-(\th-e^{-(1-\th)\gl}) \th\gl n
 \bb| \geq \dd \bb] \leq 2 e^{-\Omega (\min\{ \dd,
\frac{\dd^2}{(1-\tho ) n}\})}.$$
\end{cor}

 To estimate $A(\th)$, it is now enough to estimate
 $F(\th)$ by \raf{mainob}.
Once good estimations for $F(\th)$  are established, we may take
similar (but slightly more complicated) approaches used in Section
\ref{score}. We  consider  an (imaginary) secondary stack with
parameter $\rho$, or simply $\rho$-secondary stack. Initially, the
secondary stack with parameter $\rho$ consists of the first $\rho
n$ vertices $v_{_0}, ..., v_{_{\rho n-1}}$ of $V$. The set of
those $\rho n$ vertices is denoted by $V_\rho$. Whenever the
primary stack is empty,  the first vertex in the secondary stack
that has not been activated is to be activated. Its clones are
activated too and put into the primary stack. The activated vertex
as well as vertices activated by other means is no longer in the
secondary stack. If the secondary stack is empty, go back to the
regular procedure. This does not change the $P$-process at all,
but will be used just for the analysis. Let $\tau_\rho$ be the
largest $\tau$ such that, at $\tau \gl$, the primary stack becomes
empty after  the secondary stack is empty. Thus, once the cut-off
line reaches  $\tau_{\rho} \gl$, no active clones are provided
from the secondary stack.  Denoted by
 $C(\rho)$  is the union of the components containing any
vertex in $V_\rho$.

The following lemma is useful to predict how large $\tau_{\rho}$
is.
\begin{lem}\label{secq}  Suppose $ 0< \gd, \rho < 1$ and $\tho,
\thto\leq 1$ are uniformly bounded from below by $0$. Then
$$\pr [ \tau_{_ {\!\rho}} \geq
\tho] \leq \pr [ \min_{\th:\tho \leq \th\leq 1 }  B(\th) \leq -
(1-\gd)\tho \gl e^{-(1-\tho)\gl} \rho n ] + 2e^{-\Omega (
\gd^2\rho n)} ,
$$
and conversely,
$$ \pr [ \tau_{_{\!\rho}} \leq
\thto ] \leq \pr [ B (\thto) \geq  - (1+\gd)\thto \gl
e^{-(1-\thto)\gl} \rho n] + 2e^{-\Omega (\gd^2  \rho n)}. $$
\end{lem}

\bn \pf For simplicity we will write $\tau$ and $W$ for
$\tau_\rho$ and $V_\rho$, respectively.  Since, at $\tau\gl $, the
primary stack is empty for the first time after no
 vertex is left in the  secondary stack, $C(\rho)$
  is exactly $ W \cup  V(\tau) $. And, all
clones of vertices in $W \cup V(\tau)  $ must have been matched.
Thus,
$$ M(\tau ) +N (W) \geq N(\tau),
~~ {\rm or ~ equivalently} ~~ B(\tau_{  }) \geq - N(W) ,
$$ where, in general,  $N(V')$ is the number of clones of $v\in V'$.
The inequality can be strict when there are $\tau\gl$-large clones
of $v\in W$.  However, clones of a vertex $v\in W$ that has no
$\tau\gl$-large clone are not counted in $M(\tau ) $. Thus,
$$ M(\tau)  +N(W(\tau)) = N(\tau) ,
~~ {\rm or ~ equivalently} ~~ B(\tau_{  }) = - N(W(\tau  )),
$$
where  $W(\th)$  is the set of vertices in $W$ that have no
$\th\gl$-large clone. Thus, $\tau$ is the unique  $\th$ such that
$B(\th) = - N(W(\th  ))$ and $B(\th') >  - N(W(\th'  ))$ for all
$\th'> \th$.

If $\tau_{ } \geq \tho$, then $B(\th) = - N(W(\th ))$ for some
$\th$ with $\tho \leq \th \leq 1$. As $ W(\tho)\sub W(\th)$ for
such $\th$, we have
 $$ B(\th) \leq  - N(W(\tho  ))
, ~~\mbox{ for some $\th$ in the range $\tho \leq \th \leq 1$}.$$
This  implies that
$$ \pr [ \tau \geq \tho] \leq \pr [N(W(\tho)) < (1-\gd)\tho\gl e^{-(1-\tho)\gl} \rho n
] + \pr[\min_{\th:\tho\leq \th\leq 1} B(\th) \leq - (1-\gd)\tho\gl
e^{-(1-\tho)\gl} \rho n ].
$$
For $\th$,  in general,
$$ N(W(\th)) = \sum_{v\in W}  d_v (\th) 1( \bd_v (\th) =0)$$
is a sum of i.i.d.  random variables with mean $\th\gl
e^{-(1-\th)\gl} $, it is easy to check by the generalized Chernoff
bound that
  \beq{nwth}
 \pr [ |N(W(\th)) -\th\gl e^{-(1-\th)\gl} \rho n |
   \geq  \gd \th\gl e^{-(1-\th)\gl} \rho n
 ]
  \leq 2e^{-\Omega ( \gd^2 \rho n)}. \eeq

Conversely, if $\tau_{  }\leq \thto$, then $ B(\thto) \geq
-N(W(\thto)),$ which together with \raf{nwth}  yields
 \bean \pr [ \tau \leq \thto] &\leq & \pr [N(W(\thto)) >
 (1+\gd)\thto \gl e^{-(1-\thto)\gl} \rho n] +
 \pr[ B(\thto) \geq - (1+\gd)\thto\gl e^{-(1-\thto)\gl}
\rho n] \\
&\leq & \pr[ B(\thto) \geq - (1+\gd)\thto\gl e^{-(1-\thto)\gl}
\rho n]+  2e^{-\Omega ( \gd^2 \rho n)}. \eean \qed

\mn

\newcommand{\ttau}{\tilde{\tau}}
 Since we are
interested in $\th$  close to $1$, it is convenient to define $\AA
(\th) = A(1-\th)$, $ \VV (\th) = V(1-\th)$, $\BB (\th) =
B(1-\th)$,  $\FF (\th) = F(1-\th)$ and $\ttau_{_{\!\rho}}=
1-\tau_{_{\!\rho}}$. Then Corollary \ref{mainlem} and Lemma
\ref{secq} can be written as

\begin{cor}\label{mainlem2}
For $ \tho >0  $  uniformly bounded from above by $1$ and $\dd$ in
the range $0< \dd \leq n$,
$$ \pr \bb[\, \,  \max_{\th: 0\leq \th \leq \tho } \bb| |\VV (\th)| - e^{-\th\gl}  n \bb|
\geq \dd  \bb] \leq  2 e^{-\Omega (\min\{ \dd , \frac{\dd^2}{ \tho
n}\})},
$$
and, for $\bt (\th)= (1-\th) ( 1-\th- e^{-\th\gl})\gln $,
$$ \pr \bbr \max_{\th: 0\leq \th \leq \tho }  \bb| \BB  (\th)  -\bt (\th)
 \bb| \geq \dd \bb] \leq 2 e^{-\Omega (\min\{ \dd,
\frac{\dd^2}{\tho  n}\})}.$$
\end{cor}

\begin{lem}\label{secq2}  Suppose $ 0< \gd, \rho < 1$ and $\tho, \thto\geq 0$ are
uniformly bounded from above  by $1$. Then
$$\pr [ \ttau_{_ {\!\rho}} \leq
\tho] \leq \pr [ \min_{\th:0 \leq \th\leq \tho }  \BB(\th) \leq -
(1-\gd)(1-\tho ) \gl e^{-\tho\gl} \rho n ] + 2e^{-\Omega (
\gd^2\rho n)} ,
$$
and conversely,
$$ \pr [ \ttau_{_{\!\rho}} \geq
\thto ] \leq \pr [ \BB  (\thto) \geq  - (1+\gd)(1-\thto) \gl
e^{-\thto\gl} \rho n] + 2e^{-\Omega (\gd^2  \rho n)}. $$ In
particular,  if $0< \tho, \thto \ll \gd$, then
$$\pr [ \ttau_{_ {\!\rho}} \leq
\tho] \leq \pr [ \min_{\th:0 \leq \th\leq \tho }  \BB(\th) \leq -
(1-\gd) \gl  \rho n ] + 2e^{-\Omega ( \gd^2\rho n)} ,
$$
and
$$ \pr [ \ttau_{_{\!\rho}} \geq
\thto ] \leq \pr [ \BB  (\thto) \geq  - (1+\gd) \gl  \rho n] +
2e^{-\Omega (\gd^2  \rho n)}. $$
\end{lem}

\bn {\bf Proof of Theorem \ref{main} }

\mn {\em Supercritical Region}: Suppose  $\gl=1+\eps$ with $\eps
\gg n^{-1/3}$ and $1\ll \ga \ll (\eps^3 n)^{1/2}$. Three phases
are to be  considered
 based upon on the values of $\th$.
 Let $\tho=\frac{\ga^2 }{\thl^2 n}$, $\thto = \thl - \ga(\thl n)^{-1/2}$, and $\th_{_3}=
\thl + \ga(\thl n)^{-1/2}$. (Recall $\thl$ is the larger solution
of the equation $1-\th -e^{-\th\gl} =0$.)

\mn To bound the size of the largest component, it is  enough to
show that all of the following events occur with probability
$1-e^{-\Omega (\ga^2)}$.

\mn
 (i) For $\rho = \tho\thl=\ga^2 (\thl n)^{-1}$, we have
 $\ttau_{\rho}
 \geq \tho$,
  especially $ \FF(\tho) \leq N(V_{\rho}).$

  \mn
  (ii) For $\th$ in the range $\tho\leq \th \leq \thto$, all
$\BB(\th)$ are positive.

\mn (iii)  For some $\th$ between $\thto$ and $\thth$,
$\AA(\th)=0$.

\mn

Once (i) and (ii) occur, as $\AA(\th) \geq \BB(\th)$, the vertices
activated between  $(1-\tho)\gl$ and $(1-\thto)\gl$ are all in the
same component, say $C_1$. Excluding the vertices in $V_\rho$, all
vertices that have a $(1-\thto)\gl$-large clone but no
$(1-\tho)\gl$-large clone  belong to $C_1$. That is,  $
\VV(\tho)\sm \VV(\thto)\sub C_1 \cup V_{\rho}$. Corollary
\ref{mainlem2} for $\dd=\ga (\thl n)^{1/2}/2$ gives
 $$ \pr [ |\VV(\tho)\sm \VV(\thto)| \leq e^{-\tho \gl }
 (1-e^{-(\thto-\tho)\gl} ) n - \ga (\thl n)^{1/2}] \leq
e^{-\Omega(\ga^2) }. $$ As $\tho\ll \ga (\thl n)^{-1/2} \ll
\thl\leq 1$,  $\rho n\ll \ga (\thl n)^{1/2}$, and
$\thto=\thl-\ga(\thl n)^{-1/2}$, we have
$$|C_1| \geq |\VV(\tho)\sm \VV(\thto)| - | V_{\rho}|
\geq (1-e^{-\thl\gl})n -O(\ga (n/\thl)^{1/2})
= \thl n -O(\ga  (n/\thl)^{1/2}), $$ with probability
$1-e^{-\Omega(\ga^2)}$.

Since (i) and (ii) imply $\ttau_{\rho}\geq \thto$, if (iii) occurs
in addition,  then $C_1\sub (V\sm \VV (\thth)) \cup V_\rho$.
Corollary  \ref{mainlem2} then yields $ |V\sm \VV(\thth) |\leq
(1-e^{-\thth\gl}) n + \ga (\thl n)^{1/2}$ with probability
$1-e^{-\Omega(\ga^2)}$. As $$ \rho n\ll \ga (\thl n)^{1/2}, \and
1-e^{-\thth\gl} = 1-e^{-\thl \gl}  + O(\ga(\thl n)^{-1/2})= \thl
+O(\ga(\thl n)^{-1/2}) , $$  $C_1$ is of size at most $\thl n+
O((n/\thl)^{1/2})$ with the desired probability. Replacing $\ga$
by $\gd \ga$ for an appropriate constant $\gd>0$, the bounds for
$|C_1|$  follows as desired in Theorem \ref{mainpc}.

The union $S$ of components constructed before $C_1$ is a subset
of $(\vb(\tho))\cup V_\rho$. Corollary \ref{mainlem2} with
$\dd=0.1 \tho n$ and $1-e^{-x} \leq x$ give
$$ |\vb (\tho)| \leq \tho\gl n   + 0.1 \tho n \leq
 1.1\tho \gl n  , $$
with probability $1-2e^{-\Omega (\tho n) }\geq 1-2e^{-\Omega
(\ga^2) }$. Thus,
$$|S|\leq 1.1\tho \gl n + \rho n \leq\sfrac{3\ga^2\gl}{\thl^2 }=\Theta
(\sfrac{\ga^2 }{\eps^2 }) ,$$ with probability $1- e^{-\Omega
(\ga^2)}$.

Clearly, the vertex set $V_R$ of the residual graph without $C_1$
and $S$ satisfies
$$ V(\thth) \sm V_\rho \sub V_R \sub V(\th_{_2}), $$
and hence
$$ (1-\thl)n - \ga (n/\eps)^{1/2} \leq |V_R|  \leq (1-\thl)n + \ga
(n/\eps)^{1/2} , $$ by replacing $\ga$ by $0.1 \ga$ if necessary.
Furthermore, the cut-off value $\gl^*$  when the construction of
$C_1$ is concluded is between $(1-\thth) \gl $ and $(1-\th_{_2})
\gl$, as desired. (Recall, $\thl=(2+O(\eps))\eps$.)

\old{
constant $\gd>0$, we have any component constructed before $C_1$
is of size at most $\eps^{-2} \log (\eps^3 n)$ with probability
$1-e^{-\Omega (\gd  \log (\eps^3 n))}$.

\old{
\log (\eps^3 n))^{1/2} \ll (\eps^3 n)^{1/2}$, we know  that  every
component constructed before $C_1$ is of size at most $\gb \log
(\eps^3 n)/\eps^2$ with probability $1-e^{-\Omega (\gb \log
(\eps^3 n))}$, by replacing $\gb$ by $\gd \gb$ for an appropriate
constant $\gd>0$. } 

 After $C_1$ has been constructed, the cut-off
line is less than or equal to $(1-\thto)\gl $ due to (ii). Since
$\thto= (1+o(1))\thl $ and $(1-\thto)\gl \leq 1-\Theta(\eps)$, the
process is already in the subcritical region and one may appeal to
the result for the subcritical region.}

For the proofs of (i), (ii) and (iii), we observe that there is  a
positive constant $c<1$ such that
 \beq{lower} \bt(\th)= (1-\th) (1-\th -e^{-\th\gl}) \geq c \th(\thl-\th)n, \eeq
for $0\leq \th \le \thl$, and
 \beq{upper} \bt(\th)=(1-\th) (1-\th -e^{-\th\gl}) \leq - c \thl (\th-\thl), \eeq
for $\th \geq \thl  $ with $\th$ uniformly bounded from above by
$1$.

\bn{\bf Proof of (i)} Since $\tho\ll \thl <1 $, $\rho = \tho \thl
$ and $b (\th) \geq 0$ for all $\th \leq \thl$, Lemma \ref{secq2}
for $\gd=0.5$ gives
  \bean \pr [
\ttau_{{\rho }}\leq \tho ]  &\leq & \pr [ \min_{\th:0\leq \th\leq
\tho}  \BB(\th) \leq - 0.5  \rho  \gl n ] + 2e^{-\Omega (
\rho n)} \\
&\leq & \pr [ \min_{\th:0\leq \th\leq \tho}  \BB(\th) -\bt(\th)
\leq - 0.5 \rho\gl n ] + 2e^{-\Omega ( \ga^2)} . \eean And
Corollary  \ref{mainlem2} yields
  $$\pr [ \min_{\th:0\leq \th\leq
\tho}  \BB(\th) -\bt(\th) \leq - 0.5 \rho  \gl n ] \leq
2e^{-\Omega ( \frac{\rho^2 n^2}{\tho n})}
\leq 2 e^{-\Omega (\ga^2)}.
$$

  \qed

\bn {\bf Proof of (ii)} If $\th$ is in the range $\tho^*:=\ga
(\thl n)^{-1/2} \leq \th \leq \thto$, then
$$ \bt(\th) \geq c \th (\thl -\th)\gl n \geq 0.9  c \ga (\thl n)^{-1/2}
\thl \gl n = 0.9 c\ga \gl  (\thl n)^{1/2} , $$ (see \raf{lower}).
Thus, Corollary \ref{mainlem2} yields
$$ \pr [ \min_{\th:\tho^{\! *} \leq \th\leq \thto}  \BB(\th) \leq 0]
\leq  \pr [ \min_{\th:\tho^{\! *} \leq \th\leq \thto}  \BB(\th)
-\bt(\th) \leq   -0.9 c\ga \gl  (\thl n)^{1/2}      ] \leq
2e^{-\Omega (\ga^2 )}.$$ For $\th$ between $\tho$ and $\tho^{\!
*}$, let $\tht_i= i \tho$, $i=1, ..., i_{_{\!\!f}}$, for the least
integer $i_{_{\!\! f}}$ with $i_{_{\!\! f}}\tho\geq \tho^*$. Then,
since $\bt(\th) \geq \bt(\tht_i) \geq (1+o(1))c \tht_i \thl\gl n $
for $\tht_{_i} \leq \th\leq \tht_{_{i+1}}$, applying
  Corollary \ref{mainlem2}, we have
  \bean \pr [ \min_{\th:\tht_i \leq \th \leq \tht_{i+1}} \BB(\th) \leq 0]
&\leq&   \pr [ \min_{\th:\tht_i \leq \th\leq \tht_{i+1}}  \BB(\th)
-\bt(\th) \leq   -0.9 c\tht_i  \thl \gl n       ]  \\ &\leq&
2e^{-\Omega ( \tht_i \thl^2 n )}=2e^{-\Omega (i \ga^2 )}.\eean
Thus,
$$ \pr [ \min_{\th:\thz \leq \th \leq \tho} \BB(\th) \leq 0]
\leq \sum_{i=1}^\infty 2e^{-\Omega (i \ga^2 )}=2e^{-\Omega ( \ga^2
)}.
$$

\qed

\bn {\bf Proof of (iii)} If (iii) does not occur, then (ii)
implies that $F(\thth) = F(\tho)$. Also,
$$ 0\leq \AA(\thth) \leq \BB(\thth) + \FF(\thth) = \BB(\thth) +\FF(\tho)$$
gives $\BB (\thth) \geq -\FF(\tho)$. By (i), this  yields $\BB
(\thth) \geq -N(V_{\rho}). $ As $N(V_{\rho})$ is a Poisson $ \gl
\rho n$ random variable, $$\pr[ N(V_{\rho}) \geq 1.1 \gl \rho n]
\leq 2e^{-\Omega(\rho n)}  \leq 2e^{-\Omega(\ga^2)}, $$ and hence
$$ \pr[ \BB (\thth) \geq -N(V_{\rho}) ]
\leq e^{-\Omega(\ga^2)} + \pr[ \BB (\thth) \geq -1.1 \gl \rho n ].
$$ As $\bt(\thth) \leq -  c \ga (\thl n)^{-1/2} \thl\gln = - c\ga
\gl (\thl n)^{1/2}$ (see \raf{upper}), $\rho n  = \ga^2/\thl \ll
\ga (\thl n)^{1/2}$, and
$$- \bt (\thth) - 1.1 \gl \rho n \geq    0.9 c \ga \gl (\thl
n)^{1/2}, $$ Corollary \ref{mainlem2} gives
$$
\pr[ \BB (\thth) \geq -1.1 \gl \rho n ] \leq \pr[ \BB
(\thth)-\bt(\thth)\geq 0.9c \ga  \gl (\thl n)^{1/2}   ] \leq
e^{-\Omega(\ga^2)}.
$$

\qed

\bn

\bn {\em Subcritical Region}: (Upper Bound) Suppose $\gl = 1-\eps$
with $\eps \gg n^{-1/3}$.  We take $\rho=\eps^2 / \log (\eps^3
n)$. Until the secondary stack with  $\rho=\eps^2/ \log (\eps^3
n)$ becomes empty, each free step  can be regarded as the start of
a branching process in which the number of children is the number
of newly activated clones. The branching process ends just before
the next free step. We call the $i^{\rm th}$ branching process for
the branching process initiated by  the $i^{\rm th}$ free step.
The whole process ends when no vertex is left in the secondary
stack at the conclusion of  a branching process. As the numbers of
newly activated clones are stochastically bounded by independent
Poisson $(1-\eps)$ random variables, we know that the $i^{\rm th}$
branching process dies out, say with $D_i$ descendants, for all
$i$. Let $D(1-\eps)$ be the number of descendants for the
branching process with (independent) Poisson $(1-\eps)$ children.
Then
$$ \pr [ D_i > k ] \leq \pr [ D(1-\eps) >  k] =
e^{\eps} \sum_{\ell\geq k+1} \frac{ \ell^{\ell-1}e^{-\ell}}{\ell
!} e^{(\eps+\log (1-\eps)) (\ell-1)}.
$$

As at most $D_i$  clones are involved in the $i^{\rm th}$
branching process, the size  of the component containing the
vertex activated by the $i^{\rm th}$ free step is at most $D_i$.
Observing there are at most $\rho n$ possible $i$,  we have the
following lemma.
\begin{lem} \label{up1} Suppose $\gl=1-\eps$ with $\eps \gg n^{-1/3}$. Then,
for the secondary stack $S_\rho$ with $\rho=\eps^2/ \log (\eps^3
n)$,
$$ \pr [ \, \exists \, v\in S_\rho, ~ |C_v| > k ] \leq
\frac{\eps^2e^{\eps}  n }{\log (\eps^3 n)}  \sum_{\ell\geq k+1}
\frac{ \ell^{\ell-1}e^{-\ell}}{\ell !} e^{(\eps+\log (1-\eps))
(\ell-1)}.
$$
In particular,
$$ \pr [ \, \exists \, v\in S_\rho, ~ |C_v| > k ] = O\bb( \frac{
e^{-\eps^2k}}{ k^{3/2} \log (\eps^3 n)}\bb). $$
\end{lem}

 We will now show  that the cut-off line decreases
fast enough.

\old{
 will be
enough to show that all components are small enough as desired. To
show that there is a component with the desired size, we take
$\rho= (\eps^3 n)^{-\gd}\eps^2 $ for a constant $\gd $ in the
range $0<\gd <0.1$. It will be shown that a component containing
one of $\rho n$ vertices is of size large enough as desired. } 

\begin{lem}\label{up2}
Suppose $ \eps\leq 0.01$ and $\rho= a \eps^2$ with $a\ll 1$. Then
$$ \pr [ (1-\ttau_\rho)\gl  \geq 1-\mbox{$(1+ \frac{a}{2}) \eps $} ] \leq
2e^{-\Omega( a\eps^3 n)}.
$$
\end{lem}

\pf  For $\gd=0.01$ and  $\thz=0.7a \eps$, Lemma \ref{secq2} gives
$$ \pr [ \ttau_{\rho} \leq \thz] \leq
\pr [ \min_{\th:0 \leq \th\leq \thz }  \BB(\th) \leq -
(1-\gd)(1-\thz ) \gl e^{-\thz\gl} \rho n ] + 2e^{-\Omega (
\gd^2\rho n)}. $$
%
Using
$$(1-\gd)(1-\thz ) \gl e^{-\thz\gl}
\rho n \geq 0.9a \eps^2 \gl n , $$ and
$$\bt(\th)\geq \bt(\thz) \geq - (1-\thz)(\eps\thz+\thz^2/2)
\gln\geq -0.8 a \eps^2 \gl n $$ 
for $0\leq \th \leq \thz$, we have, by Corolallry \ref{mainlem2},
that
   \bean \pr [ \min_{\th:0 \leq \th\leq \thz } \BB(\th) \leq -
(1-\gd)(1-\thz ) \gl e^{-\thz\gl} \rho n ] &\leq&
 \pr [ \min_{\th:0\leq \th\leq \thz} \BB(\th) \leq - 0.9
a\eps^2 \gl n  ]  \\
& \leq & \pr [ \min_{\th:0\leq \th\leq \thz} \BB(\th)-\bt(\th)
\leq - 0.1 a\eps^2 \gl  n ] \\
&\leq & 2e^{-\Omega( a\eps^3 n)}. \eean Since $\ttau_\rho > \thz$
implies that
$$(1-\ttau_\rho)\gl  < (1-\thz) (1-\eps) < 1-\mbox{$(1+ \frac{a}{2}) \eps $}, $$
the desired bound follows.
 \qed

 \bn

Applying Lemmas \ref{up1} and \ref{up2} iteratively for $a\ll 1$
 until the cut-off value is less than $0.99$,
we have
$$ \pr [ \, \exists \, v, ~ |C_v| > k ] \leq
\sum_{i=0}^{\infty} \frac{\eps_{_i}^2e^{\eps_{_i}} n }{\log
(\eps_{_i}^3 n)} \sum_{\ell\geq k+1} \bb( \frac{
\ell^{\ell-1}e^{-\ell}}{\ell !} e^{(\eps_{_i}+\log (1-\eps_{_i}))
(\ell-1)}+ 2e^{-\Omega(a\eps_{_i}^3 n)} \bb) + n e^{-\Omega(k)},
$$
where $\eps_{_0}=\eps$ and  $\eps_{_{i}} \geq ( 1+\frac{ a}{2})
\eps_{_{i-1}}, $ especially $\eps_{_{i}} \geq ( 1+\frac{ a}{2})^i
\eps\geq (1+\frac{ia}{2}) \eps$. For any   $\gd>0$, taking
$$a= \frac{1}{\log (\eps^3 n)} \and ~~ k= \frac{\log (\eps^3 n)
-2.5 \log \log (\eps^3 n)+c}{-(\eps +\log (1-\eps) ) }, $$ we have
$$
\eps_{_i}^2\sum_{\ell\geq k+1} \frac{ \ell^{\ell-1}e^{-\ell}}{\ell
!} e^{(\eps_{_i}+\log (1-\eps_{_i})) (\ell-1)} =
O\bb(\frac{e^{(\eps_{_i}+\log(1-\eps_{_i}))k}}{k^{3/2}}\bb) .$$
Using
$$ \eps_{_i} + \log (1-\eps_{_i}) \leq
(1+\mbox{$\frac{ia}{2}$})(\eps + \log (1-\eps)),
$$
we finally have
$$ \pr [ \, \exists \, v, ~ |C_v| > k ] = O\bb(
\frac{ n e^{(\eps+\log(1-\eps))k}}{ k^{3/2}\log (\eps^3 n)}\bb)+
2e^{-\Omega(\frac{\eps^3 n}{\log (\eps^3 n)})} +O(ne^{-\Omega
(\log(\eps^3 n+\gd \log\log(\eps^3 n))/ \eps^2)}).
$$
and,  for $\eps\ll 1$,
$$ \pr [ \, \exists \, v, ~ |C_v| > k ] \leq  2e^{-\Omega (c)}+ 2e^{-\Omega(\frac{\eps^3 n}{\log (\eps^3 n)})}
+ 2 e^{-\Omega(\log (\eps^3 n)/\eps^2)}.
$$

\bn For the lower bound, let $ a = 1/ \log (\eps^3 n) $ and $\rho
= a \eps^2$. We will approximate the size of the component $C_i$
of the vertex activated by the $i^{\rm th}$ free step, $i=1, ...,
0.9\rho n$. If the secondary stack becomes empty earlier, $C_i$
set to be empty. To be more precise, the following auxiliary
branching process is needed. Initially, there is one organism.
Generally, the number of children is given by the random variable
$(1-Y)X-Y$ where $X$ is a Poisson $1-(1+2a)\eps$ random variable
and $Y$ is an independent Bernoulli random variable with
$\pr[Y=1]=2a\eps^2$. In other words, the population is given by
$Z_0=1$, and $Z_j = Z_{j-1} + (1-Y_j)X_j-1-Y_j$, where $(X_j,
Y_j)$ are i.i.d random variables with the same distribution as
$(X,Y)$. The branching process ends when $Z_j \leq 0$. We will
couple the branching process with $C_i$ so  that, under certain
conditions that hold with sufficiently high probability,  $|C_i| $
is  at least the sum of $(1- Y_j)$'s over $j$ with $Z_{\ell} >0 $
for all $\ell \leq j$.

To estimate $C_i$'s, let $\gL_i$ be the cut-off value at the
beginning of the $i^{\rm th} $ free step.  Then the $i^{\rm th}$
free step will generate Poisson $\gL_i$ active clones. In a
subsequent step of the $i^{\rm th} $ branching process, an active
clone $x$ and the largest unmatched clone excluding $x$, say $y$,
are to be matched. If the cut-off value $\gL$ were $(1-\th)\gl$,
$\AA(\th)\not=0$, and $\VV(\th)$ were given at the beginning of
the step, then we may lower bound the probability that  $y$ has
not been activated. (If $\AA(\th)=0$, the branching process ends.)
Clearly, the set of vertices that have not been activated at the
beginning of the step contains $V^*:=V(\th) \sm V_\rho$. Thus, the
probability is at least as large as the probability that a clone
of a vertex in $V^*$ is larger than all of $\AA(\th)-1$ currently
active clones excluding $x$. Notice that the largest number
assigned to clones of vertices of $V^*$ is less than $(1-t)\gL$
with probability $e^{-t\gL |V^*|}$ and the corresponding density
function is $\gL |V^*| e^{-t\gL |V^*|}. $ Using $|V^*| \geq |\VV
(\th)|-\rho n$, we have
 \bean \pr [  \mbox{$y$ has not been activated}] &\geq &
 \int_{0}^1
(1-t)^{|\AA(\th)|-1} \gL |V^*| e^{-t\gL |V^*|}\, dt \\ & \geq&  1-
\gL|\AA(\th)| |V^*| \int_{0}^\infty t     e^{-t\gL |V^*|}\, dt \\
 &=& 1- \frac{|\AA(\th)|}{\gL |V^*| } \geq
 1- \frac{|\AA(\th)|}{\gL ( |\VV (\th)|-\rho n)}.\eean
Conditioned that $y$ has  not  been active, the number of new
active clones is a Poisson $\gL'$ random variable, where $\gL'$ is
the cut-off value  at the end of the step.

Provided  the  $i^{\rm th}$ free step started and
$$\frac{|\AA(\th)|}{\gL (|\VV (\th)|-\rho n)} \leq 2a\eps^2, \and
\gL'\geq 1- (1+2a)\eps, $$  the number of active clones after the
end of the step is at least $\AA(\th)+ (1-Y_j) X_j-1 - Y_j $.
Moreover, in each step $Y_j=0$, a distinct vertex is added to
eventually form the component $C_i$. Thus, \beq{zz} |C_i| \geq
\sum_{j\geq 0}  (1-Y_j) 1( Z_{\ell} >0, ~\forall
  ~\ell \leq j), \eeq
  where $Z_0=1$ and, for $j\geq 1$, $Z_j = Z_{j-1} + (1-Y_j)X_j -1 -Y_j$,
 or $Z_\ell = 1+ \sum_{j=1}^\ell  (1-Y_j)X_j -1
-Y_j. $
\old{
In particular, as $ we have
$$ \pr [ C_j \geq k]
\geq    \pr \bb[ \sum_{j= 1}^{(1+a)k}  (1-Y_j)\geq k \and 1 +
\sum_{j=1}^\ell  (1-Y_j)X_j -1 -Y_j >0~~ \mbox{for all $\ell=1,
...,(1+a)k$}\bb]. $$} 
Even if any of the conditions is not satisfied, we still define
the auxiliary branching process exactly the same way. However, the
coupling and the inequality it implies are no longer true.

Let $\tho= 1.4a\eps$. We will show that the following events occur
with large enough probability.

\mn $$ { (i)} \, \ttau_\rho \leq \tho, ~~
 { (ii)}\, \,  |\VV_\rho (\tho) |\geq 0.9\rho n,
 ~~ {(iii)} \, \,  |\VV (\tho)| \geq 0.9n,
~~ { (iv)} \,  \max_{\th:0\leq \th\leq \tho } \BB (\th) \leq 0.1
\rho n,
$$ and ${ (v)}\,$the
number of clones of vertices in $V_\rho$ is less than $1.1 \rho
n$. Notice that $(i)\cap(ii)$ implies that there are at least
$0.9\rho n$ free steps, and $(i)$ gives that the cut-off value
never becomes less than $1- (1+2a)\eps$ during the entire $0.9\rho
n$ branching processes. As $\AA(\th)$ with $\th \leq \ttau_\rho$
is at most $\BB(\th)$ plus the number of clones of vertices in
$V_\rho$,
 $(i)\cap (iii) \cap (iv) \cap (v)$ implies that
 $$\frac{|\AA(\th)|}{\gL (|\VV (\th)|-\rho n)}
 \leq \frac{1.2\rho n}{(1-\tho)\gl (0.9n -\rho n)}
 \leq 2a\eps^2.$$
Therefore, \raf{zz} gives
  \bean \pr [ |C_i| \leq k, ~ ~\forall i=1,..., 0.9\rho n]
&\leq&  \pr [ \neg(i)]+  \pr [ \neg(ii)]+  \pr [ \neg(iii)]+  \pr
[ \neg(iv)]+  \pr [ \neg(v)] \\
& & +  \pr\bb[ \sum_{j\geq 0} (1-Y_j) 1(Z_\ell>0, ~\forall \ell
\leq j) \leq k \bb]^{0.9\rho n}.\eean

First,
 Lemma \ref{secq2} gives
 $$ \pr [ \neg (i) ] = \pr[ \ttau_\rho \geq \tho] \leq
 \pr [ \BB  (\tho) \geq  - 1.1(1-\tho) \gl
e^{-\tho\gl} \rho n] + 2e^{-\Omega (  \rho n)},
$$
Since $- 1.1(1-\tho) \gl e^{-\tho\gl} \rho n \geq -1.2 \rho \gl n
= -1.2 a\eps^2 \gl n, $ and
$$ - \bt(\tho)=-(1-\tho) (1-\tho-e^{-\tho \gl})\gl n
    \geq (1-\tho)\eps\tho \gl n \geq  1.3 a\eps^2\gl n,
$$
Corollary \ref{mainlem2} gives $$ \pr [ \BB  (\tho) \geq  -
1.1(1-\tho) \gl e^{-\tho\gl} \rho n] \leq \pr [ \BB  (\tho)
-\bt(\tho) \geq 0.1a\eps^2 \gl n ] \leq 2e^{-\Omega(a\eps^3 n)}.
$$
Appealing Corollary \ref{mainlem2} and using $e^{-\tho \gl}
=1+o(1)$, we have
$$ \pr[ \neg (ii)] \leq 2e^{-\Omega(a\eps^3 n)}, \and
 \pr[ \neg (iii)] \leq 2e^{-\Omega( n)}. $$
For $(iv)$, Lemma \ref{secq2} and $\bt(\th) \leq 0 $ for $\th\geq
0$ yield
$$  \pr[ \neg (iv)] \leq \pr [ \max_{\th:0\leq \th\leq \tho} |
\BB(\th)- \bt(\th) | \geq 0.1\rho n ] \leq  2e^{-\Omega(
\frac{\rho^2 n}{\tho n} )} = 2 e^{-\Omega( a\eps^3 n)} .$$
Clearly, $ \pr [ \neg (v)] \leq 2e^{-\Omega(\rho n)} =2
e^{-\Omega( a\eps^2 n)}$ as the number is a Poisson $\rho\gl  n$
random variable.

Therefore, for $$k= \frac{\log (\eps^3 n) -2.5 \log \log (\eps^3
n)-c }{-(\eps +\log (1-\eps) )}, $$
$$ \pr [ |C_i| \leq k ~\forall i=1,..., 0.9\rho n]
\leq 2 e^{-\Omega(\frac{\eps^3 n}{\log (\eps^3 n)})} + \pr\bb[
\sum_{j\geq 0} (1-Y_j) 1(Z_\ell>0 , ~\forall \ell \leq j) \leq k
\bb]^{0.9\rho n}. $$ Hence,  it is enough to show that
 \bean \pr\bb[ \sum_{j\geq 0} (1-Y_j) 1(Z_\ell >0 , ~\forall \ell \leq
j) \leq k \bb]^{0.9\rho n}&=& \bb(1-\pr\bb[ \sum_{j\geq 0} (1-Y_j)
1(Z_\ell>0, ~\forall \ell \leq j) > k \bb]\bb)^{0.9\rho n} \\
& \leq & e^{-\Omega( e^{c})}. \eean Recalling  $Z_\ell = 1+
\sum_{j=1}^\ell (1-Y_j)X_j -1 -Y_j $, we observe that, conditioned
on $Y_1=\cdots = Y_{k+1}=0$,  $Z_\ell = 1+ \sum_{j=1}^\ell X_j
-1$, $\ell=1,..., k+1$, which is exactly the population for the
Poisson branching process with mean number of children
$1-(1+2a)\eps$. As
$$ \pr [Y_1=\cdots
= Y_{k+1}=0] = (1- 2a\eps^2)^{k+1} = (1+o(1))e^{-2ak\eps^2} =
(1+o(1))e^{-2}, $$ it follows that
$$ \pr\bb[ \sum_{j\geq 0}
(1-Y_j) 1(Z_\ell>0, ~\forall \ell \leq j) > k \bb] = \Omega\bb(
\frac{e^{((1+2a)\eps + \log (1- (1+2a)\eps))k}}{\eps^2 k^{3/2}}
\bb). $$ Using $k = \Theta (\eps^{-2}\log (\eps^3 n))$ and
$(1+2a)\eps + \log (1- (1+2a)\eps)\geq (1+3a) (\eps + \log
(1-\eps))$ for $\eps\leq 0.01$, we further have
$$ \pr\bb[ \sum_{j\geq 0}
(1-Y_j) 1(Z_\ell>0, ~\forall \ell \leq j) > k \bb] = \Omega\bb(
\frac{e^{-(1+3a)(\log (\eps^3 n)-2.5 \log\log (\eps^3 n) -c)}}{
\eps^{-1}\log^{3/2} (\eps^3 n)} \bb) = \Omega\bb(\frac{e^{c} \log
(\eps^3 n)}{ \eps^{2}n } \bb) ,
$$
and hence
$$\bb(1-\pr\bb[ \sum_{j\geq 0} (1-Y_j) 1(Z_\ell>0,
~\forall \ell \leq j) > k \bb]\bb)^{0.9\rho n} \leq e^{-\Omega(
e^{c})}.
$$

 \bn
 {\em  Inside window}: Suppose $\gl = 1+\eps$ with
 $|\eps|=O(n^{-1/3})$. For a large enough constant $K$, we may sandwich
 $G_{PC} (n,p)$ between  $G_{PC} (n,p_{_1})$ and
 $G_{PC} (n,p_{_2})$, where $p_{_1} (n-1)=1-Kn^{-1/3}$ and $p_{_2}
 (n-1)=1+Kn^{-1/3}$. Thus, $G_{PC} (n,p)$ has a component of size
 $\Omega ( K^{-2} n^{2/3} \log K )$ and no component
 of size $O(K n^{2/3})$. \qed

\section{Closing Remarks}
\label{remarks}

The Poisson $\gl$-cell is introduced  to analyze properties of
$G_{PC}(n,p)$, in which degrees are i.i.d Poisson  random
variables with mean $\gl=p(n-1)$. Then various nice properties of
Poisson random variables are used to analyze sizes  of the largest
component and the $t$-core of  $G_{PC}(n,p)$. We believe that the
approaches presented in this paper are useful to analyze problems
with similar flavors, especially problems  related to branching
processes. For example, we can easily modify the proofs of Theorem
\ref{sizeofcore} to analyze the pure literal algorithm for the
random $k$-SAT problems, $k\geq 3$. Another example may be the
Karp-Sipser algorithm to find a large matching of the random
graph. (See \cite{KS, AFP}.) In a subsequent paper, we will
analyze the structures of the  $2$-core of $G(n,p)$ and the
largest strong component of the random directed graph as well as
the pure literal algorithm for the random $2$-SAT problem.

For the random (hyper)graph with a given sequence $(d_i)$, we may
also introduce the $(d_i)$-cell, in which the vertex $v_{_i}$ has
$d_i$ clones and each clone is assigned a uniform random real
number between $0$ and the average degree
$\frac{1}{n}\sum_{i=0}^{n-1} d_i$. Though it is not possible to
use all of the nice properties of Poisson random variables any
more, we believe that the $(d_i)$-cell equipped with the cut-off
line algorithm can be used to prove stronger results for the
$t$-core problems considered in various papers including \cite{CC,
FRa, FRb, JL, MM}. The case of the random $k$-SAT problem
conditioned on given degree sequence can be similarly analyzed.

A better algorithm called the unit clause algorithm  for random
$k$-SAT problems is known. (See e.g. \cite{A1,ASb,Dsur,CF,FS}) We
believe,  the unit clause algorithm and some of its variations can
be analyzed using the Poisson cloning model equipped with the
cut-off line algorithm.

Recall  that the degrees in $G(n,p)$ has the binomial distribution
with parameters $n-1$ and $ p$. By introducing the Poisson cloning
model, we somehow first take the limit of the binomial
distribution, which is the Poisson distribution. In general, many
limiting distributions like Poisson and Gaussian ones have nice
properties. In our opinion,  this is because various small
differences are eliminated by taking the limits, and limiting
distributions have some symmetric and/or invariant properties.
Thus, it may be  natural to wonder if there is an infinite graph
that shares most properties of the random graphs $G(n,p)$ with
large enough $n$. So, in a sense, the infinite graph, if exists,
can be regarded as the limit of $G(n,p)$. An infinite graph which
Aldous \cite{DA} considered to solve the linear assignment problem
may or may not be a (primitive) version of such an infinity graph.
Though it may be impossible to construct such a graph, the
approaches taken in this paper might be useful to find one, if
any.

\bn

 \bn
 {\bf Acknowledgement.} The author thanks C. Borgs, J. Chayes,
B. Bollo\'as and Y. Peres for helpful discussions.

\bn

\newcommand{\bam}{\emph{Bull. Amer. Math. Soc.}}
\newcommand{\pmd}{\emph{Publ. Math. Debrecen} }
\newcommand{\pmi}{\emph{Publ. Math. Inst. Hung. Acad. Sci.}}
\newcommand{\ama}{\emph{Acta Math.  Acad. Sci. Hung.}}
\newcommand{\ssm}{\emph{Stud. Sci. Math. Hung.}}
\newcommand{\rsa}{\emph{Random Structures and Algorithms}}
\newcommand{\tam}{\emph{Trans. Amer. Mat. Soc.}}
\newcommand{\com}{\emph{Combinatorica}}
\newcommand{\mpc}{\emph{Math. Proc. Camb. Phil. Soc.}}
\newcommand{\dm}{\emph{Discrete Math.}}
\newcommand{\spa}{\emph{Stoch. Proc. Appl.}}
\newcommand{\focs}{\emph{Proc. 22nd  Ann. IEEE Symp. Found. Comp.}}

\bn

\bn
 {\Large \bf Appendix: The proof of Lemma \ref{gmain}}

\mn For the proof of Lemma \ref{gmain}, observing  that
$$ |V(\th)|=\sum_{i=1}^m |C_i| 1( D_i (\th) \in P )
= \frac{n}{m} \sum_{i=1}^m  1( D_i (\th) \in P ),  $$ Lemma
\ref{sobr} yields the following corollary.

\begin{cor}\label{gvth} For $\th$ in the range $\tho\leq \th\leq 1$
and with the same hypotheses as in the main lemma,
$$ \pr [ |V_P (\th)| - p(\th)  n| \geq \dd ]
\leq 2e^{-\Omega (\min\{\dd , \frac{ \dd^2}{(1-\tho) n}\})}.$$
\end{cor}
As in Example \ref{ex5}, Lemma \ref{gcb2}  yields a concentration
inequality for all of $V(\th)$'s:

\begin{lem} With the same hypotheses as in the main lemma,
$$ \pr \bb[\max_{\th:\tho\leq \th \leq 1} \bb| |V(\th)| - p(\th) n \bb| \leq \dd \bb]
\leq 2e^{-\Omega (\min\{ \dd, \frac{ \dd^2}{(1-\tho) n}\})}.$$
\end{lem}

\pf   Observing that, for $\tho\leq \th\leq 1$,
$$  \sfrac{p(\tho)}{p(\th)}\bb| |V (\th)| - p(\th) n)\bb|
\leq \bb| |V (\tho)|- p(\tho)n \bb| + \bb| |V(\tho)| -
 \sfrac{p(\tho)}{p(\th)}|V(\th)|\bb|, $$
we  set $\Gamma (\th) = | |V (\th)| - p(\th) n|,$   $$
 \psi =  \sfrac{1}{p(\tho)}\Gamma  (\tho),   \and \psi_\th = \sfrac{p(\th)}{p(\tho)}
 \bb| |V (\tho)| -   \sfrac{p(\tho)}{p(\th)} |V (\th)|\bb|. $$
Clearly,  $ \Gamma  (\th ) \leq \psi + \psi_\th$. Corollary
\ref{gvth} gives \beq{gpsi1} \pr [ \psi \geq \dd/2] \leq
2e^{-\Omega(\min\{ \dd , \frac{ \dd^2}{(1-\tho) n}\})}. \eeq
 Suppose $\{X_{\th'}:=V(\th') \}_{\th\leq \th'\leq 1}$ is given, especially
$V (\th)$ is given. Then, since  $P$ is increasing, we may write
$|V(\tho)|$ as
$$ |V (\tho)|=\sum_{i: C_i \sub V(\th)}  |C_i|  1( D_i (\tho) \in
P) ,
$$ with
$$ \pr [  D_i (\tho) \in
P |\{X_{\th'} \}_{\th'\leq \th} ]= \pr [ D_i (\tho) \in P | C_i
\sub V (\th) ] = \frac{p(\tho)}{p(\th)}=:p(\tho, \th),$$  for $i$
with $C_i \sub V (\th)$. Lemma \ref{sobr} then gives
 \beq{gpsi2} \pr [ \psi_{\th} \geq \dd/2 ]
\leq 2e^{-\Omega(
     \min \{ p(\tho, \th) \dd , \frac{p(\tho, \th)\dd^2}{(1-p(\tho,\th)) |V
     (\th)|}\} )}\leq 2e^{-\Omega(\min\{p(\tho\!) \dd , \frac{p(\tho \!)
\dd^2}{(1-p(\tho\!)) n}\})} \leq 2e^{-\Omega(\min\{\dd , \frac{
\dd^2}{(1-\tho) n}\})}. \eeq Lemma \ref{gcb2}  together with
\raf{gpsi1} and \raf{gpsi2} yields the desired inequality.

\qed

\bn

We now estimate $M(\th)$. First, since $\sum_{v\in V} \bd_v(\th)$
is a Poisson random variable with mean $(1-\th)\gl n $,
  \beq{gfirstsum} \pr\bbr  \bb|\sum_{v\in V} \bd_v(\th) - (1-\th)\gl n \bb| \geq \dd/2\bb]
\leq 2 e^{-\min\{ \dd, \frac{\dd^2}{(1-\th) n}\}}. \eeq For the
second sum, observe that
$$ \sum_{v\in V} d_v (\th) 1(v\not\in V (\th))
=  \sum_{i=1}^m  \bb( \sum_{v\in C_i}  d_v(\th) \bb) 1(D_i(\th)
\not \in P)$$ is a sum of $m$ i.i.d random variables with
$$ E\bb[ \bb( \sum_{v\in C_i}  d_v(\th) \bb) 1(D_i(\th)
\not \in P)\bb]=E\bb[ \sum_{v\in C_i}  d_v(\th) - \bb( \sum_{v\in
C_i} d_v (\th)\bb) 1(D_i (\th) \in P ) \bb] = (\th\gl- q(\th))
|C_i|.
$$
  Moreover, since $P$ is an
increasing property and $\sum_{v\in C_i}  d_v(\th) $ is a Poisson
$\th\gl |C_i|$ random variable, FKG inequality (see e.g.  Chapter
6 of \cite{AS}) gives
$$
E\bb[ \bb( \bb( \sum_{v\in C_i}  d_v(\th) \bb) 1(D_i(\th) \not \in
P)\bb)^i \bb] \leq E\bb[ \bb( \sum_{v\in C_i}  d_v(\th) \bb)^i\bb]
\pr [D_i(\th) \not \in P ] = O(1-p(\th)),
$$
for all fixed $i$, e.g. $i=1,2,3$. Thus one may take $\xi_{_0}=1$
and $a_i, b_i=\Theta(1-p(\th))$ to satisfy all the conditions to
apply the generalized Chernoff bound and to obtain
 $$ \pr \bbr  \bb|\sum_{v\in V} d_v (\th) 1(v\not\in V (\th))
   - (\th\gl-q(\th)) n \bb| \geq \dd \bb] \leq
2e^{-\Omega(\min\{ \dd, \frac{\dd^2}{(1-p(\th))n}\})},
 $$
provided $|C_1|=O(1)$. This together with \raf{gfirstsum} implies
that if $|C_1|=O(1)$  and $1-p(\th) = O(1-\th)$, then
   \beq{gsmth}
 \pr \bbr  \bb|M(\th)
   - (\gl-q(\th)) n \bb| \geq \dd \bb] \leq
2e^{-\Omega(\min\{ \dd, \frac{\dd^2}{(1-\th)n} \})}.
 \eeq

We proof a concentration result for all of $M(\th)$ that together
with the cut-off lemma implies Lemma \ref{gmain} follows.

\begin{lem} \label{gamth} With the same hypotheses as in the main lemma,
$$ \pr \bb[\max_{\th:\tho\leq \th \leq 1} \bb| M (\th) - (\gl-q(\th)) n \bb| \leq \dd \bb]
\leq 2e^{-\Omega (\min\{ \dd, \frac{\dd^2}{(1-\tho) n}\})}.$$
\end{lem}

\pf Clearly,
 $$| M(\th) -
 (\gl-q(\th)) n|\leq |M(\tho)-  (\gl-q(\tho)) n|+
 |M(\tho)-M (\th) - (q(\th)-q(\tho)) n|. $$
Let $\Gamma (\th) =| M (\th) -
 (\gl-q(\th)) n|$,
 $$ \psi = \Gamma (\tho),~~ \psi_{\th} =|M (\tho)-M (\th) - (q(\th)-q(\tho))
 n|, $$
and $\Phi_\th$ is the event $|V (\th) - p(\th)n| \leq
\frac{p(\tho)\dd}{4\gl}$. Then, \raf{gsmth} gives
$$ \pr [ \psi \geq \dd/2 ] \leq 2e^{-\Omega(\min \{\dd,
\frac{\dd^2}{(1-\tho) n}\})}. $$ For $\psi_\th$, suppose
$\{X_{\th'} :=(V(\th'), M (\th'))\}_{\th\leq \th' \leq 1}$ is
given. Using $$ M(\th) = \sum_{v\in V} \bd_v (\th) + d_v (\th ) 1
(v\not \in V (\th))=  \sum_{v\in V}   d_v (1) -  d_v(\th)1 (v \in
V (\th)) ,$$  we obtain
 \bean M( \tho) - M (\th) &=& \sum_{v\in  V} d_v (\th)
 1(v\in  V (\th))
- d_v (\tho) 1(v\in V (\tho)) \\ &=&  \sum_{i=1}^m \bb(\sum_{v\in
C_i} d_v (\th)\bb) 1(C_i \sub V (\th)) - \bb( \sum_{v\in C_i} d_v
(\tho) \bb) 1(C_i \sub  V (\tho)).\eean
 Once $V(\th)$ is given, the
distributions of $$Y_i :=\bb(\sum_{v\in C_i} d_v (\th)\bb)1(C_i
\sub V (\th)) - \bb( \sum_{v\in C_i} d_v (\tho) \bb) 1(C_i \sub V
(\tho))$$  depend on neither $\{M(\th')\}_{\th\leq \th'\leq 1} $
nor  $\{V(\th')\}_{\th<  \th'\leq 1} $ and hence, for $C_i \sub V
(\th)$,
 \bean
 E\bb[  Y_i
\bb| \{ X_{\th'}\}_{\th\leq \th' \leq 1}\bb] & =& E\bb[\sum_{v\in
C_i} d_v (\th) - \bb( \sum_{v\in C_i}
d_v (\tho) \bb) 1(v\in V (\tho))\bb| C_i \sub V (\th) \bb] \\
&= &
 \frac{q(\th) - q(\tho)}{p(\th)} |C_i|.
\eean  
If $C_i \not\sub V (\th)$, $Y_i=0$   since $P$ is increasing.

Also, for $C_i \sub V(\th)$, we may write
 $$ Y_i = \sum_{v\in C_i}
d_v (\th) - d_v (\tho) + d_v (\tho) 1(C_i \not\sub V (\tho))
$$
and \bean
 E\bbr Y_i^2 \bb| C_i \sub V (\th)\bb]
&\leq &
 2 E\bbr \bb( \sum_{v\in C_i} d_v (\th) - d_v (\tho) \bb)^2 \bb| C_i \sub V (\th) \bb]
 \\
 & &  +  2E\bbr
  \bb( \sum_{v\in C_i}
d_v (\tho) 1(C_i \not\sub V (\tho))\bb)^2 \bb| C_i \sub V (\th)
\bb].
 \eean
For $j=1,2,3$,
$$
 E\bbr \bb( \sum_{v\in C_i} d_v (\th) - d_v (\tho) \bb)^j \bb|
C_i \sub V (\th) \bb]
  \leq
p(\th)^{-1} E\bbr \bb( \sum_{v\in C_i} d_v (\th) - d_v (\tho)
\bb)^2j \bb] = O(\th-\tho) = O(1-\tho) $$ for $p(\th)\geq p(\tho)
=\Omega(1)$ and $\sum_{v\in C_i} d_v (\th) - d_v (\tho) $ is a
Poisson  random variable with mean $(\th-\tho)\gl
|C_i|=O(\th-\tho)$. For the second term, FKG inequality gives
 \bean
E\bbr
  \bb( \sum_{v\in C_i}
d_v (\tho) 1(C_i \not\sub V (\tho))\bb)^j \bb| C_i \sub V (\th)
\bb] &\leq & p(\th)^{-1} E\bbr
  \bb( \sum_{v\in C_i}
d_v (\tho) 1(C_i \not\sub V (\tho))\bb)^j\bb] \\
 &\leq&
 p(\tho)^{-1} E\bbr
  \bb( \sum_{v\in C_i}
d_v (\tho) \bb)^j\bb]E [1(C_i \not\sub V (\tho))]\\
&=& O(1-p(\tho)) = O(1-\tho),
 \eean
 for $j=1,2,3$.
Therefore,
   $$ E\bbr \bb( Y_i - E[Y_i] \bb)^2 \bb| \{ X_{\th'}\}_{\th\leq \th' \leq 1}\bb] \leq
 E\bb[  Y_i^2 \bb|  \{ X_{\th'}\}_{\th\leq \th' \leq 1}\bb]= O(1-\tho).
$$ Similarly,
for $\xi$ in the range $|\xi|\leq \xi_{_0}=1$, it is not hard to
show
 $$ \bb| E\bbr (Y_i-E[Y_i])^3 e^{\xi(Y_i-E[Y_i])}  \bb| \{ X_{\th'}\}_{\th\leq \th' \leq 1}\bb]
 \bb|
 = O(1-\tho). $$
Applying the generalized Chernoff bound, we have
$$ \pr \bbr \bb| \sum_{j=1}^m  Y_i -
\frac{q(\th)-q(\tho)}{p(\th)} | V (\th) | \bb| \geq \dd/4  \bb| \{
X_{\th'}\}_{\th\leq \th' \leq 1}\bb] \leq 2e^{-\Omega(\min \{\dd,
\frac{\dd^2}{1-\tho}\})}.
$$

Finally, as    the event $\Phi_\th$ guarantees
$$ \frac{q(\th)-q(\tho)}{p(\th)}  \bb| | V (\th) | -
p(\th)n\bb| \leq \dd/4  $$ for  $p(\tho) \leq p(\th)$ and $q(\th)
\leq \gl$, we  have
  \bean  & & 1(\Phi_\th) \pr \bbr \bb| \sum_{j=1}^m  Y_i -
(q(\th)-q(\tho))n \bb| \geq \dd/2 \bb| \{ X_{\th'}\}_{\th\leq \th'
\leq 1}\bb]  \\ & & \leq \pr \bbr \bb| \sum_{j=1}^m Y_i -
\frac{q(\th)-q(\tho)}{p(\th)} | V (\th) | \bb| \geq \dd/4 \bb|
\{ X_{\th'}\}_{\th\leq \th' \leq 1}\bb]\\
& & \leq 2e^{-\Omega(\min \{\dd, \frac{\dd^2}{1-\tho}\})}.  \eean
Lemma \ref{gcb2}  yields the desired inequality.

\qed

\end{document}